\documentclass[letterpaper,10pt,twoside]{article}
\usepackage[utf8]{inputenc}
 
\usepackage[OT1]{fontenc}

\usepackage[letterpaper,margin=2cm,
  headheight=12pt,headsep=14pt,footskip=25pt]{geometry}
\usepackage{amsmath,amssymb,amsthm,mathtools,mathrsfs}
\usepackage[cal=rsfso]{mathalpha}
\usepackage{graphicx,booktabs,array}
\usepackage[font=small,labelfont=normalfont,labelsep=colon]{caption}
\usepackage{microtype}
\usepackage{xcolor}
\usepackage{xurl}
\usepackage{balance}
\usepackage{fancyhdr}
\usepackage{hyperref}
\mathtoolsset{showonlyrefs}
\definecolor{referenceblue}{rgb}{0.26,0.32,0.82}
\hypersetup{
  colorlinks=true,
  linkcolor=red,
  citecolor=referenceblue,
  urlcolor=referenceblue,
  pdfauthor={Alexandre Autran}
}

\allowdisplaybreaks[1]
\fancypagestyle{plain}{\fancyhf{}\fancyfoot[C]{\normalfont\small\thepage}}

\makeatletter
\renewcommand\section{\@startsection{section}{1}{\z@}
  {-3.5ex plus -1ex minus -0.2ex}{2.3ex plus 0.2ex}
  {\normalfont\large\bfseries}}
\renewcommand\subsection{\@startsection{subsection}{2}{\z@}
  {-3.25ex plus -1ex minus -0.2ex}{2ex plus 0.2ex}
  {\normalfont\normalsize\bfseries}}
\renewcommand\maketitle{\begingroup
  \centering
  \hrule height 4pt
  \vskip 0.25in
  {\normalfont\fontsize{14.4}{16}\selectfont\bfseries\@title\par}
  \vskip 0.25in
  \hrule height 1pt
  \vskip 0.47in
  {\normalfont\normalsize\@author\par}
  \endgroup
  \vskip 0.35in}
\makeatother
\renewenvironment{abstract}{\begin{center}\normalfont\large\bfseries Abstract\end{center}
  \vspace{7pt}
  \begin{list}{}{\setlength{\leftmargin}{0pt}\setlength{\rightmargin}{0pt}\setlength{\topsep}{0pt}\setlength{\parsep}{5pt}\setlength{\itemsep}{0pt}}
  \item\relax\normalfont\normalsize\ignorespaces
}{\end{list}\par\vspace{8pt}}


\newtheorem{theorem}{Theorem}[section]
\newtheorem{proposition}[theorem]{Proposition}
\newtheorem{lemma}[theorem]{Lemma}
\newtheorem{corollary}[theorem]{Corollary}
\theoremstyle{remark}\newtheorem{remark}[theorem]{Remark}
\newcommand{\R}{\mathbf R}
\newcommand{\Uad}{\mathcal U_{\mathrm{ad}}}
\newcommand{\Uk}{\mathcal U_{\mathrm{ad},h,\tau}}
\newcommand{\ip}[2]{\langle #1,#2\rangle_h}
\newcommand{\ipt}[2]{\langle #1,#2\rangle_{h,\tau}}
\newcommand{\normh}[1]{\lVert #1\rVert_h}
\newcommand{\normt}[1]{\lVert #1\rVert_{h,\tau}}
\newcommand{\Wht}{W_{h,\tau}}
\newcommand{\one}{\mathbf 1}
\newcommand{\diag}{\operatorname{diag}}
\newcommand{\proj}{\Pi_{[0,U_{\max}]}}
\newcommand{\argmin}{\operatorname*{arg\,min}}
\newcommand{\Ih}{\mathcal I_h}
\newcommand{\Qh}{\mathcal Q_h}
\newcommand{\Qht}{\mathcal Q_{h,\tau}}
\newcommand{\avg}{\mathcal A_{h,\tau}}
\title{Quadrature-Induced Local Minima in Structure-Preserving Lotka--Volterra Control}
\author{\textbf{Alexandre Autran}\\
Department of Mathematics
\\
\'{E}cole normale sup\'{e}rieure de Rennes, France\\
\href{mailto:alexandre.autran@ens-rennes.fr}{alexandre.autran@ens-rennes.fr}
}

\newcommand{\Vh}{V_h}
\newcommand{\Y}{\mathcal Y}
\newcommand{\mapdef}[5]{#1:\left|\begin{array}{@{}r@{\;}c@{\;}l@{}}
    #2&\longrightarrow&#3\\[2pt]
    #4&\longmapsto&\displaystyle #5
  \end{array}\right.}
\newcommand{\byeq}[2]{\mathrel{\overset{\mathclap{\text{\eqref{#1}}}}{#2}}}
\newcommand{\Yht}{\mathcal Y_{h,\tau}}

\newsavebox{\mapdefbox}
\makeatletter
\newcommand{\compactmapdef}[5]{\begingroup
  \mapdeffallback{#1}{#2}{#3}{#4}{#5}
  \endgroup}
\newcommand{\mapdeffallback}[5]{
  \sbox{\mapdefbox}{$\displaystyle
    \begin{array}{@{}l@{\;}c@{\;}c@{\;}l@{}}
      #1:&#2&\longrightarrow&#3\\[2pt]
      &#4&\longmapsto&\displaystyle #5
    \end{array}$}
  \ifdim\wd\mapdefbox<\dimexpr\linewidth-1em\relax
    \usebox{\mapdefbox}
  \else
    \if@twocolumn
      \begin{gathered}
        #1:#2\longrightarrow#3\\[2pt]
        #4\longmapsto\displaystyle #5
      \end{gathered}
    \else
      \fullwidthmapdef{#1}{#2}{#3}{#4}{#5}
    \fi
  \fi}
\makeatother

\begin{document}
\raggedbottom
\twocolumn[{\maketitle}]
\thispagestyle{plain}
\setlength{\parindent}{0pt}
\setlength{\parskip}{5pt plus 1pt minus 1pt}
\let\fullwidthmapdef\mapdef
\let\mapdef\compactmapdef
\begin{abstract}
We study how the choice of harvest quadrature affects optimal control in a positivity-preserving Lotka--Volterra scheme. With harvesting alone, counting the prey removed at each step gives a strictly convex cost. Using the prey density at the start of the step overestimates the catch and can create a strict local minimum with cost above the global minimum. We certify the global minimizers in a two-step example by rational root isolation. A second example shows that, with both species present and positive diffusion, a strict local minimizer for the old-state rule admits a descent direction under the flux rule.

For nonnegative $H^1$ initial data and positive diffusivities in one space dimension, the two rules give costs, weighted gradients and Hessians that differ by $O(\tau)$ at the same control, uniformly in the spatial mesh. Here $\tau$ is the time step. If the old-state Hessian is uniformly coercive near a stationary control, the flux objective has a local minimizer within $O(\tau)$ of it. With uniformly bounded tracking targets, a large enough effort penalty makes both objectives strongly convex for all sufficiently small time steps. We also derive exact discrete adjoints and Hessian actions, prove a common variational limit with an IMEX discretization, and compare stationary controls across grids. The same harvest accounting applies to weighted production--destruction systems.
\end{abstract}
{\small\noindent\textbf{Keywords:}
reaction--diffusion, optimal harvesting, positive scheme, exact discrete adjoint, reduced Hessian, variational convergence, forward--backward sweep.
\par}\vspace{4pt}
\setlength{\parskip}{4pt plus 1pt minus 1pt}
\setlength{\abovedisplayskip}{9pt plus 2pt minus 4pt}
\setlength{\belowdisplayskip}{9pt plus 2pt minus 4pt}
\setlength{\abovedisplayshortskip}{0pt plus 3pt}
\setlength{\belowdisplayshortskip}{5pt plus 3pt minus 3pt}
\section{Introduction}
Harvesting removes prey and leaves less food for predators. Diffusion spreads both effects in space \cite{AutranStage2026,Murray2002,Hrinca2002,Apreutesei2010,Leung1995}. We write $v_1(x,t)$ and $v_2(x,t)$ for predator and prey densities, and $u(x,t)$ for fishing effort. The objective balances population tracking and the cost of effort against the value of the catch. Discretizing it requires a choice about how to count the prey removed during each time step \cite{BurkardtGunzburgerPeterson2002,Campos2015,LiuWang2019,ShakyaSinha2020,Harzer2025}.

Consider a single population with harvesting alone. Let $y^n$ be its density at the start of step $n$, $u^n$ the effort, $q$ the catchability and $\tau$ the step length. The positive update gives
\[
 y^{n+1}=\frac{y^n}{1+\tau q u^n},
 \qquad y^n-y^{n+1}=\tau q u^n y^{n+1}.
\]
The quantity on the right is exactly the catch over the step. If we use the starting density to count the catch, the excess is
\[
 \tau q u^n y^n-(y^n-y^{n+1})
 =\tau^2q^2(u^n)^2y^{n+1}.
\]
This extra reward depends on how effort is spread over time. Even with only two steps, it can favor putting all the effort into one step or create a strict local minimum above the global minimum. We derive the corresponding cost identity in \eqref{eq:pure-defect-new}.

For the predator--prey system, we compare three discrete problems. The first, I, uses implicit--explicit (IMEX) Euler and counts the catch using the starting prey density. The state stays nonnegative under a time-step restriction \cite{Hundsdorfer2003,Quarteroni2007}. The other two problems, $P_0$ and $P_1$, use the same Patankar scheme, which preserves nonnegativity for every time step \cite{Burchard2003}. The old-state rule $P_0$ uses the starting density. The flux rule $P_1$ uses the end-of-step density that appears in the population balance. Comparing I with $P_0$ isolates the choice of state scheme. Comparing $P_0$ with $P_1$ isolates the harvest quadrature. Proposition~\ref{prop:coupled-local-minimum} gives a two-species example in which changing this rule removes a strict local minimum, for any positive diffusivities on a fixed spatial grid.

With diffusion in one dimension, we bound the differences between the costs, gradients and Hessians of $P_0$ and $P_1$ at a common control. The bounds are uniform in the spatial mesh. Under a local coercivity assumption, they also bound the distance between the corresponding stationary controls. All three discrete problems have the same weak-$*$ variational limit. Controls whose costs approach the discrete global minimum have subsequences converging strongly in $L^2$ to continuous global minimizers. Each strict $L^2$ local minimizer of the continuous problem can also be approximated by discrete local minimizers.

We derive exact discrete adjoints and Hessian-vector products \cite{Giles2000,Hinze2009}. For the flux rule, even the last harvest depends on the final prey state, so its derivative enters the terminal adjoint condition. We use these derivatives to compare four optimizers and to examine how temporal switching and curvature change with the grid and the initial control. For relaxed forward--backward sweep, strict complementarity identifies the active bounds locally. Once they are identified, the smallest and largest eigenvalues of the positive definite free Hessian determine the range of relaxation parameters that gives a contraction.

This analysis builds on work on discrete adjoints and second derivatives \cite{Hager2000,Herty2013,SanzSerna2016,Ito2021} and on positivity, conservation and accuracy for Patankar schemes \cite{Burchard2003,KopeczMeister2018,Torlo2022}. Cacace et al.\ use a modified Patankar semi-Lagrangian scheme to solve production--destruction control problems by dynamic programming \cite{Cacace2026}. Here we fix the state scheme and ask how the choice of harvest reward changes the controls.
\section{Continuous problem}
Let $\Omega\subset\R^d$ be a bounded Lipschitz domain, with $d\in\mathbf N=\{1,2,\ldots\}$, and let $T>0$ be the final time. On the space--time domain $Q\coloneqq \Omega\times(0,T)$, write $v_1,v_2:Q\to\R$ for the densities and $u:Q\to\R$ for effort. In the following system, the evolution equations hold in $Q$, the boundary conditions on $\partial\Omega\times(0,T)$, and the initial conditions in $\Omega$ for $i=1,2$:
\begin{equation}\label{eq:pde}
\left\{\begin{aligned}
\partial_t v_1-\mu_1\Delta v_1&=v_1(-m+\alpha v_2),\\
\partial_t v_2-\mu_2\Delta v_2&=v_2(r-\beta v_1)-quv_2,\\
\partial_n v_1=\partial_n v_2&=0,\\
v_i(0)&=v_i^0.
\end{aligned}\right.
\end{equation}
Here $\Delta$ is the spatial Laplacian and $\partial_n$ denotes the outward normal derivative. We write $v_i(t)$ for the spatial profile $v_i(\cdot,t)$ and $v_i^0$ for its initial value. The model has positive diffusivities $\mu_i$, mortality $m$, growth rate $r$, interaction coefficients $\alpha,\beta$, catchability $q$, and effort bound $U_{\max}$. We assume $\mu_i,m,r,\alpha,\beta,q,U_{\max}>0$. The weights $\gamma_i\ge0$, $\lambda>0$, and $\eta\ge0$ determine the penalties for tracking and effort, and the reward for harvest. For a domain $D$, $L^p(D)$, $1\le p<\infty$, consists of measurable functions with finite norm
\[
\|f\|_{L^p(D)}\coloneqq\left(\int_D|f|^p\right)^{1/p}.
\]
In $L^\infty(D)$ we use the essential supremum norm. We write $|D|$ for the measure of $D$.
We constrain effort to the \emph{admissible set}
\begin{equation}\label{eq:admissible}
\Uad\coloneqq \{u\in L^\infty(Q):0\le u\le U_{\max}\text{ a.e.}\}.
\end{equation}
For target profiles $v_i^d\in L^2(Q)$, write $v\coloneqq (v_1,v_2)$ and define the \emph{cost functional}
\begin{equation}\label{eq:cost}
\mapdef{J}{L^2(Q)^2\times L^2(Q)}{\R}{(v,u)}{
\begin{aligned}[t]
&\frac12\sum_{i=1}^2\gamma_i\|v_i-v_i^d\|_{L^2(Q)}^2\\
&\quad+\frac\lambda2\|u\|_{L^2(Q)}^2-\eta q\int_Q uv_2.
\end{aligned}}
\end{equation}
\subsection{Spaces, bounds, and existence}
All spaces are real. The Sobolev space $H^1(D)$ consists of $L^2(D)$ functions whose first weak derivatives also lie in $L^2(D)$, with
\[
\|f\|_{H^1(D)}^2=\|f\|_{L^2(D)}^2+\|\nabla f\|_{L^2(D)}^2.
\]
Here $\nabla f$ is the weak spatial gradient.
For a Banach space $X$, $L^p(0,T;X)$ denotes strongly measurable $X$-valued functions whose norm lies in scalar $L^p(0,T)$. The space $H^1(0,T;X)$ also requires a weak time derivative in $L^2(0,T;X)$. The continuous dual $X^*$ consists of continuous linear functionals, with pairing $\langle\ell,x\rangle_{X^*,X}\coloneqq\ell(x)$. We write $X^m$ for the product of $m$ copies of $X$. Set
\begin{equation}\label{eq:spaces}
\begin{aligned}
\mathcal W&\coloneqq L^2(0,T;H^1(\Omega))\cap H^1(0,T;H^1(\Omega)^*),\\
\Y&\coloneqq \mathcal W\cap L^\infty(Q),
\end{aligned}
\end{equation}
with norm
\begin{equation*}
\|y\|_{\Y}\coloneqq\|y\|_{L^2(0,T;H^1)}
+\|\partial_ty\|_{L^2(0,T;H^1{}^*)}+\|y\|_{L^\infty(Q)}.
\end{equation*}
The \emph{control-to-state map} and \emph{reduced objective} are
\[
 \mapdef{S}{\Uad}{\Y^2}{u}{v\text{ solving \eqref{eq:pde}}}
\]
and
\[
 \mapdef{j}{\Uad}{\R}{u}{J(S(u),u)}.
\]
We write $S_i(u)$ for the $i$th component of $S(u)$. The embedding $\mathcal W\hookrightarrow C([0,T];L^2(\Omega))$ gives each state an $L^2(\Omega)$-continuous representative. We use the standard product norm on products of Hilbert spaces and the sum norm on $\Y^2$. Write $(a,b)_D\coloneqq \int_Dab$, using the Euclidean dot product for vectors.
We use $\to$ for strong convergence and $\rightharpoonup$ for weak convergence. In $L^\infty(D)$, $f_k\rightharpoonup^* f$ means $\int_D f_k\phi\to\int_D f\phi$ for every $\phi\in L^1(D)$.
For $\phi\in H^1(\Omega)$ and almost every $t\in(0,T)$, the weak equation is
\begin{equation}\label{eq:weakstate}
\langle\partial_tv_i,\phi\rangle_{H^1{}^*,H^1}
+\mu_i(\nabla v_i,\nabla\phi)_\Omega=(f_i(v,u),\phi)_\Omega,
\end{equation}
where the \emph{reaction map} $f\coloneqq (f_1,f_2)$ is
\[
 \mapdef{f}{\R^2\times\R}{\R^2}{(v,u)}{
 \binom{v_1(-m+\alpha v_2)}{v_2(r-\beta v_1-qu)}}.
\]
This weak equation includes the homogeneous Neumann conditions. We omit the domain in $\|\cdot\|_p$ when the context fixes it, and allow constants $C$ to change between estimates. Each result states what these constants depend on. For $b_\varepsilon>0$, as $\varepsilon\to0$, $O(b_\varepsilon)$ denotes a quantity bounded in the relevant norm by $Cb_\varepsilon$, and $o(b_\varepsilon)$ one whose norm divided by $b_\varepsilon$ tends to zero.
Set $V_i^0\coloneqq \|v_i^0\|_{L^\infty(\Omega)}$.
\begin{theorem}[State bounds]\label{th:continuous}
For $u\in\Uad$ and $0\le v_i^0\in L^\infty(\Omega)$, \eqref{eq:pde} has a unique weak solution in $\Y^2$. For almost every $(x,t)\in\Omega\times(0,T)$, we have
\begin{equation}\label{eq:preybound}
0\le v_2(x,t)\le V_2^0 e^{rt},
\end{equation}
and
\begin{equation}\label{eq:predbound}
0\le v_1(x,t)\le V_1^0\exp\!\left(-mt+\frac{\alpha V_2^0}{r}(e^{rt}-1)\right).
\end{equation}
Define the weighted population by
\begin{equation}
\mathcal M(t)\coloneqq \beta\int_\Omega v_1(t)+\alpha\int_\Omega v_2(t).
\end{equation}
For almost every $t\in(0,T)$, it satisfies
\begin{equation}\label{eq:balancec}
\mathcal M'(t)=-\beta m\int_\Omega v_1+\alpha r\int_\Omega v_2
-\alpha q\int_\Omega uv_2.
\end{equation}
Thus, for every $t\in[0,T]$,
\begin{equation*}
\mathcal M(t)\le e^{rt}\mathcal M(0).
\end{equation*}
If $\mu_1=\mu_2$, then for almost every $(x,t)\in\Omega\times(0,T)$, we also have
\begin{equation}\label{eq:equaldiff}
\beta v_1+\alpha v_2\le e^{rt}\|\beta v_1^0+\alpha v_2^0\|_\infty.
\end{equation}
\end{theorem}
\begin{theorem}[Existence of an optimal control]\label{th:exist}
Under the assumptions above, $\min_{u\in\Uad}j(u)$ has a solution.
\end{theorem}
Appendix~\ref{app:continuous} proves both results.
\subsection{First and second derivatives}
For a differentiable map $F$, $DF(u)[d]$ is its derivative at $u$ applied to direction $d$, and $D^2F(u)[d,e]$ its second derivative in directions $d,e$. We also use $F'(u)$ and $F''(u)$. Subscripts specify the variable being differentiated. Derivatives between Banach spaces are understood in the Fr\'echet sense.
The state map has a $C^2$ extension $S:L^\infty(Q)\to\Y^2$ (Appendix~\ref{app:continuous}). The same formula extends $j$ to this space. With the adjoint $p\coloneqq(p_1,p_2)\in\mathcal W^2$ defined by \eqref{eq:adjc} in Appendix~\ref{app:continuous-derivative-formulas}, the \emph{gradient map}
\[
 \mapdef{g}{L^\infty(Q)}{L^2(Q)}{u}{\lambda u-qS_2(u)(\eta+p_2(u))}
\]
gives
\begin{equation}\label{eq:gradc}
\mapdef{j'(u)}{L^\infty(Q)}{\R}{\xi}{\int_Q g(u)\xi}.
\end{equation}
The \emph{scalar projection}
\begin{equation*}
\mapdef{\proj}{\R}{[0,U_{\max}]}{s}{\min\{U_{\max},\max\{0,s\}\}}
\end{equation*}
acts pointwise on functions and componentwise on arrays. At a local minimizer $u^*$, for every $u\in\Uad$,
\begin{equation*}
\int_Q g(u^*)(u-u^*)\ge0.
\end{equation*}
Equivalently, for almost every $(x,t)\in\Omega\times(0,T)$,
\begin{equation}\label{eq:projectc}
u^*=\proj\!\left(\frac q\lambda S_2(u^*)(\eta+p_2(u^*))\right).
\end{equation}
For second-order sufficient conditions in PDE-constrained optimization, see \cite{Troeltzsch2010,CasasChrysafinosMateos2025}.
\section{Two state schemes and three discrete objectives}
\subsection{Spatial operator and quadrature}\label{sec:spatial-operator}
For the convergence analysis and experiments, take $\Omega=(0,1)$. For $M,N\in\mathbf N$, set $h\coloneqq 1/M$ and $\tau\coloneqq T/N$. For $0\le j\le M$ and $0\le n\le N$, define $x_j\coloneqq jh$ and $t_n\coloneqq n\tau$, respectively.
We identify arrays with their coordinate vectors. The entry $v_{i,j}^n$ is the discrete density of species $i$ at node $x_j$ and time $t_n$. $v_i^n$ collects its nodal values, and $u_j^n$ is the fishing effort. State trajectories lie in $\Yht\coloneqq \R^{2(N+1)(M+1)}$, the adjoints $(p_i^n)_{1\le n\le N}$ in $\R^{2N(M+1)}$, and the controls in $\R^{N(M+1)}$. Define the \emph{discrete Neumann Laplacian} by
\begin{equation}\label{eq:laplacian}
\mapdef{L_h}{\R^{M+1}}{\R^{M+1}}{y}{
 h^{-2}\begin{pmatrix}
 2(y_1-y_0)\\
 (y_{j-1}-2y_j+y_{j+1})_{1\le j<M}\\
 2(y_{M-1}-y_M)
 \end{pmatrix}}.
\end{equation}
The \emph{dual cells} are $C_0\coloneqq [0,h/2)$, $C_M\coloneqq (1-h/2,1]$, and for $1\le j<M$,
\begin{equation*}
C_j\coloneqq (x_j-h/2,x_j+h/2).
\end{equation*}
Their lengths are $w_0\coloneqq w_M\coloneqq h/2$, $w_j\coloneqq h$ otherwise. For $a,b\in\R^{M+1}$ and $a,b\in\R^{N(M+1)}$, respectively,
\begin{equation}\label{eq:inner}
\begin{aligned}
\ip a b&\coloneqq \sum_{j=0}^M w_ja_jb_j,\qquad\text{and}\\
\ipt a b&\coloneqq \tau\sum_{n=0}^{N-1}\ip{a^n}{b^n}.
\end{aligned}
\end{equation}
These inner products give the norms $\normh a^2\coloneqq \ip aa$ and $\normt a^2\coloneqq \ipt aa$. For arrays, $\|\cdot\|_\infty$ denotes the largest absolute entry and $\|\cdot\|_2$ the Euclidean norm. Set
\begin{equation*}
\begin{aligned}
W_h&\coloneqq \diag(w_0,\ldots,w_M)\in\R^{(M+1)\times(M+1)},\\
\Wht&\coloneqq \tau I_N\otimes W_h\in\R^{N(M+1)\times N(M+1)}.
\end{aligned}
\end{equation*}
Here $\diag(a)$ is the diagonal matrix with entries from $a$, $I_s$ the $s\times s$ identity ($I$ when its size is clear), $\otimes$ the Kronecker product, $\one\in\R^{M+1}$ the vector of ones, the superscript $T$ matrix transposition, and $\circ$ componentwise multiplication. Inequalities between arrays or matrices are componentwise. Then $W_hL_h=L_h^TW_h$, $L_h\one=0$, and
\begin{equation*}
\ip{-L_hy}{y}=h^{-1}\sum_{j=0}^{M-1}(y_{j+1}-y_j)^2.
\end{equation*}
The \emph{control box} is $$\Uk\coloneqq [0,U_{\max}]^{N(M+1)}\subset\R^{N(M+1)}.$$
Let $\Vh\subset H^1(0,1)$ be the space of continuous functions that are affine on each $[x_j,x_{j+1}]$, $j=0,\ldots,M-1$. For $x\in[x_j,x_{j+1}]$, define the \emph{affine reconstruction}
\[
 \mapdef{\Ih}{\R^{M+1}}{\Vh}{y}{
 \left[x\mapsto\frac{x_{j+1}-x}{h}y_j+\frac{x-x_j}{h}y_{j+1}\right]}
\]
and the \emph{piecewise-constant reconstruction}
\[
 \mapdef{\Qh}{\R^{M+1}}{L^2(0,1)}{y}{\sum_{j=0}^M y_j\mathbf1_{C_j}}.
\]
Here $\mathbf1_E$ is the indicator of $E$, equal to one on $E$ and zero outside. Define
\begin{equation}\label{eq:diffusion_form}
\mapdef{a_h}{\R^{M+1}\times\R^{M+1}}{\R}{(y,z)}{
 \begin{aligned}[t]
 &h^{-1}\sum_{j=0}^{M-1}(y_{j+1}-y_j)(z_{j+1}-z_j)\\
 &\qquad=\ip{-L_hy}{z}
 \end{aligned}},
\end{equation}
with seminorm
\begin{equation*}
|y|_{1,h}\coloneqq a_h(y,y)^{1/2}=\|\partial_x\Ih y\|_{L^2(0,1)}.
\end{equation*}
The \emph{space--time reconstruction} is
\begin{equation}\label{eq:reconstruction}
 \mapdef{\Qht}{\R^{N(M+1)}}{L^2(Q)}{u}{
 \sum_{n=0}^{N-1}\sum_{j=0}^M u_j^n\mathbf1_{C_j\times[t_n,t_{n+1})}}.
\end{equation}
It is an isometry: $\normt u=\|\Qht u\|_{L^2(Q)}$. The \emph{averaging map} is
\begin{equation}\label{eq:averages}
 \mapdef{\avg}{L^2(Q)}{\R^{N(M+1)}}{f}{
 \left(\frac1{\tau w_j}\int_{t_n}^{t_{n+1}}\!\int_{C_j}f(x,t)\,\mathrm{d}x\,\mathrm{d}t\right)_{\substack{0\le n<N\\0\le j\le M}}}.
\end{equation}
Averaging preserves the control box. It also gives $\Qht\avg f\xrightarrow[h,\tau\to0]{L^2(Q)}f$. Choose target arrays $v_i^{d,n}\in\R^{M+1}$ whose reconstructions converge strongly to $v_i^d$ in $L^2(Q)$. Cell averages provide such arrays. In the experiments, we instead sample the smooth stationary targets at the grid nodes. We also write $v_i^0$ for the initial grid vector and specify its discretization when needed.
\subsection{IMEX Euler step}
For $n=0,\ldots,N-1$, set $A_i\coloneqq I-\tau\mu_iL_h\in\R^{(M+1)\times(M+1)}$ and solve
\begin{equation}\label{eq:imex}
\left\{\begin{aligned}
A_1v_1^{n+1}&=v_1^n+\tau v_1^n\circ(-m+\alpha v_2^n),\\
A_2v_2^{n+1}&=v_2^n+\tau v_2^n\circ(r-\beta v_1^n-qu^n).
\end{aligned}\right.
\end{equation}
Since $A_i^{-1}\ge0$ componentwise, positivity at each step follows from
\begin{equation}\label{eq:predator-positivity}
\min_{0\le j\le M}\{1+\tau(-m+\alpha v_{2,j}^n)\}\ge0,
\end{equation}
and
\begin{equation}\label{eq:chi}
\min_{0\le j\le M}\{1+\tau(r-\beta v_{1,j}^n-qu_j^n)\}\ge0.
\end{equation}
\subsection{Positive transfer-implicit step}
Define the loss matrices
\begin{equation}\label{eq:BC}
\begin{aligned}
B&\coloneqq I-\tau\mu_1L_h+\tau mI,\\
C^n&\coloneqq I-\tau\mu_2L_h+\tau\diag(\beta v_1^n+qu^n).
\end{aligned}
\end{equation}
Solve first for the prey, then for the predator:
\begin{equation}\label{eq:positive}
\left\{\begin{aligned}
C^n v_2^{n+1}&=(1+\tau r)v_2^n,\\
Bv_1^{n+1}&=v_1^n+\tau\alpha v_1^n\circ v_2^{n+1}.
\end{aligned}\right.
\end{equation}
Both equations use the transfer $v_1^n\circ v_2^{n+1}$, which cancels in the weighted balance.
\begin{theorem}[Positivity and exact balance]\label{th:positive}
For $u\in\Uk$, nonnegative initial vectors, and any $\tau>0$, \eqref{eq:positive} has a unique nonnegative trajectory. With
\begin{equation}
\mathcal M_h^n\coloneqq \beta\ip{v_1^n}{\one}+\alpha\ip{v_2^n}{\one},
\end{equation}
it satisfies
\begin{equation}\label{eq:balancep}
\begin{aligned}
\frac{\mathcal M_h^{n+1}-\mathcal M_h^n}{\tau}
={}&-\beta m\ip{v_1^{n+1}}\one+\alpha r\ip{v_2^n}\one\\
&-\alpha q\ip{u^n}{v_2^{n+1}},
\end{aligned}
\end{equation}
and $$\mathcal M_h^n\le(1+\tau r)^n\mathcal M_h^0\le e^{rt_n}\mathcal M_h^0.$$
\end{theorem}
See Appendix~\ref{app:discrete}.
\begin{proposition}[Bounds uniform in the control and mesh]\label{prop:uniform}
Let $K_2\coloneqq \|v_2^0\|_\infty e^{rT}$, $K_1\coloneqq \|v_1^0\|_\infty e^{\alpha K_2T}$, and 
\[
\tau_0\coloneqq \min\{m^{-1},(\beta K_1+qU_{\max})^{-1}\}.
\]
For every $\tau>0$ with $N\tau=T$, the positive scheme satisfies $0\le v_i^n\le K_i\one$. For $\tau\le\tau_0$, IMEX satisfies the same bounds independently of $h$. If $\mu_1=\mu_2$, both nonnegative schemes satisfy
\begin{equation}\label{eq:equaldiscrete}
\beta v_1^n+\alpha v_2^n\le
(1+\tau r)^n\|\beta v_1^0+\alpha v_2^0\|_\infty\one.
\end{equation}
\end{proposition}
Appendix~\ref{app:discrete} proves these bounds.
\begin{remark}[Other spatial operators]
Positivity, balance, and adjoint identities hold for $L\in\R^{s\times s}$ and positive diagonal $W_h\in\R^{s\times s}$ satisfying $L_{ij}\ge0$ for $i\ne j$, $L\one=0$, $W_hL=L^TW_h$. The shifted loss matrices are nonsingular M-matrices. We use the one-dimensional grid for the convergence analysis.
\end{remark}
\subsection{Harvest quadrature is part of the discrete problem}
The switch $\zeta\in\{0,1\}$ selects old-state reward when $\zeta=0$ and new-state reward when $\zeta=1$. Define
\begin{equation}\label{eq:discretecost}
\mapdef{J_{h,\tau}^\zeta}{\Yht\times\R^{N(M+1)}}{\R}{(v,u)}{
 \begin{aligned}[t]
 &\frac\tau2\sum_{n=0}^{N-1}\sum_{i=1}^2\gamma_i\normh{v_i^n-v_i^{d,n}}^2
 +\frac\lambda2\normt u^2\\
 &\quad-\eta q\tau\sum_{n=0}^{N-1}\ip{u^n}{(1-\zeta)v_2^n+\zeta v_2^{n+1}}.
 \end{aligned}}
\end{equation}
The state maps are
\[
 \mapdef{S^{\rm I}_{h,\tau}}{\Uk}{\Yht}{u}{v\text{ solving \eqref{eq:imex}}}
\]
and
\[
 \mapdef{S^+_{h,\tau}}{\Uk}{\Yht}{u}{v\text{ solving \eqref{eq:positive}}}.
\]
Their reduced costs are
\[
 \mapdef{j^{\mathrm I}_{h,\tau}}{\Uk}{\R}{u}{J_{h,\tau}^0(S^{\rm I}_{h,\tau}(u),u)}
\]
and
\[
 \mapdef{j^{+,\zeta}_{h,\tau}}{\Uk}{\R}{u}{J_{h,\tau}^\zeta(S^+_{h,\tau}(u),u)}.
\]
The two costs for the positive scheme share one trajectory, but $P_1$ rewards exactly the removal flux $q\tau\ip{u^n}{v_2^{n+1}}$ in \eqref{eq:balancep}.
\section{Exact discrete adjoints and Hessian actions}
For any of the three discrete costs $j$, the state map and cost are smooth on an open neighborhood $\mathcal O_{h,\tau}\subset\R^{N(M+1)}$ of $\Uk$ (Appendix~\ref{app:derivatives}). Write $\nabla_Ej$ for the Euclidean gradient. The \emph{weighted gradient} is
\[
 \mapdef{g}{\mathcal O_{h,\tau}}{\R^{N(M+1)}}{u}{\Wht^{-1}\nabla_Ej(u)}
\]
and the \emph{weighted Hessian action} is
\[
 \mapdef{H(u)}{\R^{N(M+1)}}{\R^{N(M+1)}}{d}{Dg(u)d},
\]
where $Dg(u)d$ is the directional derivative of $g$.
Equivalently,
\begin{equation}\label{eq:euclidean}
(\nabla_E j)_{n,j}=\tau w_jg_{n,j},\qquad\nabla_E j=\Wht g.
\end{equation}
Write $\mathcal L(X)$ for bounded linear maps from a normed space $X$ into itself. In $\R^{N(M+1)}$, every linear map is bounded. The state and adjoint maps are
\[
\mapdef{S}{\mathcal O_{h,\tau}}{\Yht}{u}{
  S^{\rm I}_{h,\tau}(u)\text{ or }S^+_{h,\tau}(u)
}
\]
and
\[
\mapdef{p}{\mathcal O_{h,\tau}}{\R^{2N(M+1)}}{u}{
  \bigl(p_i^n(u)\bigr)_{\substack{i=1,2\\1\le n\le N}}
}.
\]
The adjoints solve \eqref{eq:adjI} for IMEX and
\eqref{eq:terminalP}--\eqref{eq:adjP} for the positive scheme.
Their derivatives are
\[
 \mapdef{DS(u)}{\R^{N(M+1)}}{\Yht}{d}{z}
\]
and
\[
 \mapdef{Dp(u)}{\R^{N(M+1)}}{\R^{2N(M+1)}}{d}{\xi}.
\]
Algorithms that use Euclidean gradients take $\Wht g$ as input. Appendix~\ref{app:imex-derivatives} gives the IMEX adjoint, gradient and Hessian action.
\subsection{Positive scheme: terminal and backward equations}
For $P_\zeta$, the terminal conditions are
\begin{equation}\label{eq:terminalP}
p_1^N=0,\qquad C^{N-1}p_2^N=-\tau\zeta\eta q u^{N-1}.
\end{equation}
For $k=N-1,\ldots,1$, solve $p_1^k$ before $p_2^k$:
\begin{equation}\label{eq:adjP}
\left\{\begin{aligned}
Bp_1^k={}&p_1^{k+1}+\tau\big[\gamma_1(v_1^k-v_1^{d,k})\\
&+v_2^{k+1}\circ(\alpha p_1^{k+1}-\beta p_2^{k+1})\big],\\
C^{k-1}p_2^k={}&(1+\tau r)p_2^{k+1}\\
&+\tau\big[\gamma_2(v_2^k-v_2^{d,k})+\alpha v_1^{k-1}\circ p_1^k\\
&-\eta q((1-\zeta)u^k+\zeta u^{k-1})\big].
\end{aligned}\right.
\end{equation}
The gradient is
\begin{equation}\label{eq:gradP}
\begin{split}
g^{+,\zeta,n}={}&\lambda u^n-\eta q[(1-\zeta)v_2^n+\zeta v_2^{n+1}]\\
&-qv_2^{n+1}\circ p_2^{n+1}.
\end{split}
\end{equation}
For $\varphi\in\R^{M+1}$, differentiating the terminal harvest term gives
\[
D_{v_2^N}J_{h,\tau}^{\zeta}(v,u)[\varphi]
=-\tau\zeta\eta q\ip{u^{N-1}}{\varphi},
\]
hence \eqref{eq:terminalP}. The initial state is fixed, so the adjoint stops at $n=1$. Appendix~\ref{app:positive-derivatives} gives the positive tangent, tangent-adjoint and Hessian formulas.
\begin{proposition}[Exactness and weighted symmetry]\label{prop:symmetry}
For each of I, $P_0$, and $P_1$, the formulas in this section and Appendix~\ref{app:derivatives} satisfy
\begin{equation}\label{eq:exactness}
Dj(u)d=\ipt{g(u)}d,
\end{equation}
and
\begin{equation*}
D^2j(u)[d,e]=\ipt{H(u)d}e=\ipt d{H(u)e}.
\end{equation*}
\end{proposition}
With the states and adjoints stored, one tangent-state sweep and one tangent-adjoint sweep give a Hessian action. Both use tridiagonal spatial solves. The total cost is $O(NM)$. Appendix~\ref{app:derivatives} derives the formulas.
\section{Convergence of the discrete control problems}\label{sec:variational}
Throughout this section, $\Omega=(0,1)$ and
\begin{equation}\label{eq:dataV}
0\le v_i^0\in H^1(0,1)\cap L^\infty(0,1),\quad v_i^d\in L^2(Q).
\end{equation}
We interpolate the initial data at the nodes and choose target reconstructions that converge strongly in $L^2(Q)$. Let $h_k,\allowbreak\tau_k\xrightarrow[k\to\infty]{}0$ independently. For IMEX, the step restriction $\tau_k\le\tau_0$ holds for all sufficiently large $k$ by Proposition~\ref{prop:uniform}. Reconstructions indexed by $k$ use $(h_k,\tau_k)$.
\subsection{Uniform state estimates and compactness}
On $[t_n,t_{n+1})$, define the \emph{forward time difference}
\begin{equation}\label{eq:time_reconstructions}
D_\tau v_i^n\coloneqq (v_i^{n+1}-v_i^n)/\tau\in\R^{M+1},
\end{equation}
and the \emph{old-state and new-state reconstructions}
\begin{equation*}
V_i^-\coloneqq \Qh v_i^n,\qquad\text{and}\qquad V_i^+\coloneqq \Qh v_i^{n+1},
\end{equation*}
with \emph{linear time interpolation}
\begin{equation*}
\widehat V_i(x,t)\coloneqq \frac{t_{n+1}-t}{\tau}\Ih v_i^n(x)
+\frac{t-t_n}{\tau}\Ih v_i^{n+1}(x).
\end{equation*}
On each grid, $V_i^\pm\in L^\infty(Q)$ and $\widehat V_i\in H^1(0,T;\Vh)$.
\begin{lemma}[Uniform discrete energy and time differences]\label{lem:energy}
For either scheme, a constant $C$ independent of $h$, $\tau$, and $u\in\Uk$ gives
\begin{equation}\label{eq:energy}
\max_{0\le n\le N}\|\Ih v_i^n\|_{H^1(0,1)}^2
+\tau\sum_{n=0}^{N-1}\normh{D_\tau v_i^n}^2\le C.
\end{equation}
In particular,
\begin{equation}\label{eq:jump}
\|V_i^+-V_i^-\|_{L^2(Q)}\le C\tau,
\end{equation}
and
\begin{equation}\label{eq:reconstructionerror}
\|\Qh v_i^n-\Ih v_i^n\|_{L^2(0,1)}\le Ch\|\partial_x\Ih v_i^n\|_{L^2}.
\end{equation}
\end{lemma}
\begin{theorem}[State consistency for varying controls]\label{th:statecompact}
Let $u_k\in \mathcal U_{\mathrm{ad},h_k,\tau_k}$ and set $U_k\coloneqq \Qht u_k$. Suppose that in $L^\infty(Q)$ we have $U_k\mathrel{\underset{k\to\infty}{\rightharpoonup^{*}}}u$. Then, for either state scheme, both reconstructions converge strongly in $L^2(Q)$:
\begin{equation}\label{eq:statecompact}
V_{i,k}^-\xrightarrow[k\to\infty]{}S_i(u),\qquad
V_{i,k}^+\xrightarrow[k\to\infty]{}S_i(u).
\end{equation}
The linear interpolants also converge strongly in $L^2(Q)$: $\widehat V_{i,k}\xrightarrow[k\to\infty]{}S_i(u)$. This limit is the nonnegative solution of \eqref{eq:pde} with initial data $v_i^0$.
\end{theorem}
Appendix~\ref{app:convergence} proves that the sequence
$(\widehat V_{i,k})_k$ is bounded in
$L^\infty(0,T;H^1(0,1))\cap H^1(0,T;L^2(0,1))$, 
then identifies the limit using compactness \cite{Simon1987} and products with one weakly convergent factor and one strongly convergent factor.
\subsection{The two harvest objectives have the same limit}
\begin{proposition}[Uniform harvest-quadrature defect]\label{prop:quadrature}
For the positive state scheme and data \eqref{eq:dataV},
\begin{equation}\label{eq:quadraturedefect}
\sup_{u\in\Uk}|j^{+,1}_{h,\tau}(u)-j^{+,0}_{h,\tau}(u)|\le C\tau.
\end{equation}
The same bound holds for the difference between their global minimum values.
\end{proposition}
See Appendix~\ref{app:proof-prop-quadrature}.
\begin{remark}[Control comparison under quadratic growth]\label{rem:quadrature-controls}
Write $j_\zeta\coloneqq j^{+,\zeta}_{h,\tau}$, $\delta\coloneqq \sup_{u\in\Uk}|j_1(u)-j_0(u)|$. For $u_0\in\Uk$, $r>0$, set
\[
K_r\coloneqq \{u\in\Uk:\normt{u-u_0}\le r\}.
\]
Suppose that for some $c>0$ and every $u\in K_r$,
\[
j_0(u)\ge j_0(u_0)+c\normt{u-u_0}^2.
\]
If $u_1\in K_r$ and
$j_1(u_1)\le j_1(u_0)+\varepsilon$, with $\varepsilon\ge0$, then
\[
\normt{u_1-u_0}\le\sqrt{\frac{2\delta+\varepsilon}{c}}.
\]
If $2\delta<cr^2$, minimizing $j_1$ on $K_r$ gives a local minimizer on $\Uk$ (Appendix~\ref{app:proof-growth-comparison}). Since $\delta\le C\tau$ by Proposition~\ref{prop:quadrature}, a mesh-independent positive lower bound for $c$ and $\varepsilon=O(\tau)$ give an $O(\tau^{1/2})$ bound for $u_1\in K_r$. \end{remark}
\subsection{Variational convergence and strong control convergence}
For I, $P_0$, or $P_1$, denote the reduced cost on grid $k$ by $j_k$ and define
\begin{equation}\label{eq:extended}
 \mapdef{\mathcal F_k}{\Uad}{\R\cup\{+\infty\}}{U}{
 \begin{cases}
 j_k(u),&U=\mathcal Q_{h_k,\tau_k}u,\quad u\in \mathcal U_{\mathrm{ad},h_k,\tau_k},\\
 +\infty,&\text{otherwise}.
 \end{cases}}
\end{equation}
We write $\argmin_K f$ for the set of points in $K$ where $f$ attains its minimum.
\begin{theorem}[Variational consistency and global minimizers]\label{th:gamma}
Assume \eqref{eq:dataV}. For each of I, $P_0$, and $P_1$:
\begin{enumerate}
\item If, in $L^\infty(Q)$, $U_k\mathrel{\underset{k\to\infty}{\rightharpoonup^{*}}}u$, then
\begin{equation}\label{eq:liminf}
j(u)\le\liminf_{k\to\infty}\mathcal F_k(U_k).
\end{equation}
\item For every $u\in\Uad$, $u_k\coloneqq \mathcal A_{h_k,\tau_k}u$ is a recovery sequence:
\begin{equation}\label{eq:recovery}
\begin{gathered}
\mathcal Q_{h_k,\tau_k}u_k\xrightarrow[k\to\infty]{L^2(Q)}u,\\
j_k(u_k)\xrightarrow[k\to\infty]{}j(u).
\end{gathered}
\end{equation}
\item The minima converge:
\begin{equation}\label{eq:minima}
\min_{\mathcal U_{\mathrm{ad},h_k,\tau_k}}j_k\xrightarrow[k\to\infty]{}\min_{\Uad}j.
\end{equation}
If $\varepsilon_k\to0$ and
\begin{equation}\label{eq:epsglobal}
j_k(u_k)\le\min j_k+\varepsilon_k,
\end{equation}
then, along a subsequence, the reconstructed controls converge strongly in $L^2(Q)$ to some $u^*\in\argmin_{\Uad}j$:
\begin{equation}\label{eq:strongcontrols}
\mathcal Q_{h_k,\tau_k}u_k\xrightarrow[k\to\infty]{}u^*.
\end{equation}
The states converge as in \eqref{eq:statecompact}.
\end{enumerate}
\end{theorem}
The liminf bound and recovery-sequence property in the first two statements define sequential $\Gamma$-convergence in weak-$*$ $L^\infty(Q)$. Appendix~\ref{app:convergence} proves the theorem.
\begin{corollary}[Comparison of global minimizers]\label{cor:unique}
If the continuous global minimizer is unique, every sequence satisfying \eqref{eq:epsglobal} with $\varepsilon_k\to0$ converges strongly to it, so reconstructed controls from different schemes approach each other in $L^2(Q)$. When there are several continuous minimizers, the controls converge to their set.
\end{corollary}
\begin{corollary}[Approximation of a strict local minimizer]\label{cor:local}
Under Theorem~\ref{th:gamma}, let $u^*\in\Uad$ be a strict $L^2(Q)$ local minimizer: for some $r>0$ and every $u\in\Uad$ with $0<\|u-u^*\|_{L^2(Q)}\le r$, $j(u)>j(u^*).$ For each scheme, there are local minimizers on all sufficiently fine grids whose reconstructions converge strongly to $u^*$ in $L^2(Q)$.
\end{corollary}
Appendix~\ref{app:convergence} constructs these controls by minimizing over the closed $L^2(Q)$ ball of radius $r$ around $u^*$.
With the regularity and sampling required by Proposition~\ref{prop:accuracy}, either scheme has $O(\tau+h^2)$ state error (Appendix~\ref{app:smooth-accuracy}).
\section{What the harvest quadrature changes}
\label{sec:quadrature-mechanism-new}
For $\zeta\in\{0,1\}$, write $j_\zeta\coloneqq j^{+,\zeta}_{h,\tau}$. The \emph{quadrature defect} is
\[
 \mapdef{D}{\Uk}{\R}{u}{j_1(u)-j_0(u)}.
\]
Both costs use the same state, so
\begin{equation}\label{eq:defect-exact-new}
 D(u)\byeq{eq:discretecost}{=}-\eta q\tau\sum_{n=0}^{N-1}
       \langle u^n,v_2^{n+1}-v_2^n\rangle_h.
\end{equation}
Where the curvature is small, even a small cost difference can change the local minima.
\subsection{A reduced model}
On $\R^N$, use the weighted inner product and norm
\[
\langle\xi,\chi\rangle_\tau\coloneqq \tau\sum_{n=0}^{N-1}\xi^n\chi^n,\qquad
\|\xi\|_\tau^2\coloneqq \langle\xi,\xi\rangle_\tau.
\]
For $u\in [0,U_{\max}]^N$ and initial population $y^0>0$, let $y^n$ denote the population at time $t_n$ in the pure-harvesting model. For $0\le n<N$,
\begin{equation}\label{eq:pure-harvest-new}
y^{n+1}=\frac{y^n}{1+\tau q u^n}.
\end{equation}
The trajectory map and costs are
\[
 \mapdef{y}{[0,U_{\max}]^N}{(0,\infty)^{N+1}}{u}{
 \left(y^0\prod_{n=0}^{k-1}(1+\tau q u^n)^{-1}\right)_{k=0}^N},
\]
\begin{equation}\label{eq:pure-cost0-new}
\mapdef{j_0}{[0,U_{\max}]^N}{\R}{u}{
  \frac{\lambda\tau}{2}\sum_{n=0}^{N-1}(u^n)^2
  -\eta\tau q\sum_{n=0}^{N-1}u^ny^n
}
\end{equation}
and
\begin{equation}\label{eq:pure-cost1-new}
\mapdef{j_1}{[0,U_{\max}]^N}{\R}{u}{
  \frac{\lambda\tau}{2}\sum_{n=0}^{N-1}(u^n)^2
  -\eta\tau q\sum_{n=0}^{N-1}u^ny^{n+1}
}.
\end{equation}
Taking the empty product as $1$ and summing the population removed at each step gives
\begin{equation}\label{eq:pure-telescope-new}
\begin{aligned}
 j_1(u)&=\frac{\lambda\tau}{2}\sum_{n=0}^{N-1}(u^n)^2
 -\eta\tau q\sum_{n=0}^{N-1}u^ny^{n+1}\\
 &\byeq{eq:pure-harvest-new}{=}
   \frac{\lambda\tau}{2}\sum_{n=0}^{N-1}(u^n)^2
 -\eta\sum_{n=0}^{N-1}(y^n-y^{n+1})\\
 &=\frac{\lambda\tau}{2}\sum_{n=0}^{N-1}(u^n)^2-\eta(y^0-y^N).
\end{aligned}
\end{equation}
Old-state quadrature adds a reward term:
\begin{equation}\label{eq:pure-defect-new}
\begin{aligned}
 j_0(u)&=j_1(u)-\eta\tau q\sum_{n=0}^{N-1}u^n(y^n-y^{n+1})\\
 &\byeq{eq:pure-harvest-new}{=}j_1(u)
 -\eta\tau^2q^2\sum_{n=0}^{N-1}(u^n)^2y^{n+1}.
\end{aligned}
\end{equation}
This term rewards population already removed (Appendix~\ref{app:proof-removal-identities}).
For $u\in [0,U_{\max}]^N$ and $\zeta\in\{0,1\}$, the weighted Hessian action is
\[
 \mapdef{H_\zeta(u)}{\R^N}{\R^N}{\xi}{
 \tau^{-1}\nabla_E^2j_\zeta(u)\xi},
\]
where $\nabla_E^2$ is the Euclidean Hessian of the cost's smooth rational extension.
\begin{proposition}[Convexity and a quadrature-induced negative mode]
\label{prop:pure-curvature-new}
For every $N$, $\lambda>0$, and $\eta\ge0$, $j_1$ is strictly convex on $[0,U_{\max}]^N$. Its unique minimizer is constant in time. If it lies strictly inside the control box, its value $s$ is the unique solution of
\begin{equation}\label{eq:pure-root-new}
 \lambda s=\eta q y^0(1+\tau q s)^{-N-1}.
\end{equation}
The operators $H_\zeta(u)\in\mathcal L(\R^N)$ satisfy
\[
D^2j_\zeta(u)[\xi,\chi]=\langle H_\zeta(u)\xi,\chi\rangle_\tau
\]
and
\begin{equation}\label{eq:pure-H1-new}
\begin{aligned}
H_1(0)&=\lambda I+\eta\tau q^2y^0(\mathbf 1\mathbf 1^T+I),\\
H_0(0)&=\lambda I+\eta\tau q^2y^0(\mathbf 1\mathbf 1^T-I),
\end{aligned}
\end{equation}
where $I\in\R^{N\times N}$ and $\mathbf1\in\R^N$.
For every nonzero direction whose components sum to zero,
the Rayleigh quotient is $\lambda+\eta\tau q^2y^0$ for $H_1(0)$
and $\lambda-\eta\tau q^2y^0$ for $H_0(0)$. For $N\ge2$ and $\lambda<\eta\tau q^2y^0$, $j_0$ is nonconvex.
\end{proposition}
The proof is in Appendix~\ref{app:proof-prop-pure-curvature-new}.
For even $N$, the direction $(1,-1,\ldots,1,-1)$ shifts effort between adjacent steps while keeping the total fixed. At zero effort, the old-state curvature along this direction changes sign at $\lambda=\eta\tau q^2y^0$. This gives a mechanism for temporal oscillations.
The refinement estimate is
\[
\begin{aligned}
&\|H_1(u)-H_0(u)\|_{\mathcal L(\R^N)}\\
&\quad\le\eta\tau q^2y^0\bigl(2+4qU_{\max}T+2q^2U_{\max}^2T^2\bigr).
\end{aligned}
\]
Here the operator norm is induced by $\|\cdot\|_\tau$. Once the bound falls below $\lambda$, $j_0$ is strictly convex (Appendix~\ref{app:proof-reduced-refinement}).
\begin{proposition}[A quadrature-induced local minimum that is not global]
\label{prop:quadrature-induced-local-minimum}
For $N=2$, set $x=\tau q u^0$, $z=\tau q u^1$ and $c=\eta\tau q^2 y^0/\lambda$. If $1<c\le \tau q U_{\max}$, then $(x,z)=(0,c)$ is a strict constrained local minimizer of $j_0$ with cost strictly above the global minimum. The flux cost has a unique global minimizer.
\end{proposition}
The proof is in Appendix~\ref{app:proof-prop-quadrature-induced-local-minimum}.
Both the flux and continuous costs are strictly convex. For $u\in L^\infty(0,T)$ with $0\le u\le U_{\max}$ a.e., the continuous cost is
\[
j_{\rm c}(u)=\tfrac\lambda2\int_0^T u^2\,\mathrm{d}t-\eta y^0
+\eta y^0\exp\left(-q\int_0^T u\,\mathrm{d}t\right).
\]
Its smooth extension to $L^2(0,T)$ satisfies, for $\xi\in L^2(0,T)$,
\[
\begin{aligned}
D^2j_{\rm c}(u)[\xi,\xi]
&=\lambda\|\xi\|_{L^2(0,T)}^2\\
&\quad+\eta q^2y^0e^{-q\int_0^T u\,\mathrm{d}t}
       \Bigl(\int_0^T\xi\,\mathrm{d}t\Bigr)^2\\
&\ge\lambda\|\xi\|_{L^2(0,T)}^2.
\end{aligned}
\]
Quadrature-induced local minima can therefore occur even with a positive state step.
\subsection{Optimal allocation at fixed total effort}
For two time steps with $u^0+u^1=E$, $0<E\le U_{\max}$, set $a\coloneqq u^0\in[0,E]$, $u^1=E-a$. The restricted costs are
\[
 \mapdef{f_0}{[0,E]}{\R}{a}{
 \begin{aligned}[t]
 &\frac{\lambda\tau}{2}\big[a^2+(E-a)^2\big]-\eta\tau qy^0E\\
 &\qquad+\eta\tau^2q^2y^0\frac{a(E-a)}{1+\tau qa}
 \end{aligned}}
\]
and
\[
 \mapdef{f_1}{[0,E]}{\R}{a}{
 \begin{aligned}[t]
 &\frac{\lambda\tau}{2}\big[a^2+(E-a)^2\big]-\eta y^0\\
 &\qquad+\frac{\eta y^0}{(1+\tau qa)(1+\tau q(E-a))}
 \end{aligned}}.
\]
\begin{proposition}[Quadrature-induced extreme allocations]\label{prop:budget}
For $\lambda>0$, $f_1$ has the unique global minimizer $a=E/2$. If
\[
 \lambda\tau<\frac{\eta\tau^2q^2y^0}{(1+\tau qE)^2},
\]
then $f_0$ is strictly concave on $[0,E]$, with exactly two global minimizers $a=0,E$. The exact continuous removal model, with effort constant on each interval and the same total effort, has the unique minimizer $a=E/2$.
\end{proposition}
See Appendix~\ref{app:proof-prop-budget}.
\subsection{Harvest quadrature in the coupled system}
Return to $j_\zeta\coloneqq j^{+,\zeta}_{h,\tau}$ and $D\coloneqq j_1-j_0$. Then
\begin{equation}\label{eq:coupled-defect-new}
\begin{aligned}
D(u)&\byeq{eq:defect-exact-new}{=}
 -\eta q\tau\sum_{n=0}^{N-1}\langle u^n,v_2^{n+1}-v_2^n\rangle_h\\
&\byeq{eq:positive}{=}
 \eta q\tau^2\sum_{n=0}^{N-1}\bigl[
 q\langle (u^n)^2,v_2^{n+1}\rangle_h\\
&\qquad +\beta\langle u^n v_1^n,v_2^{n+1}\rangle_h
 -r\langle u^n,v_2^n\rangle_h\\
&\qquad -\mu_2\langle u^n,L_hv_2^{n+1}\rangle_h\bigr].
\end{aligned}
\end{equation}
Growth, predation and diffusion determine the sign: the defect can be positive or negative (Appendix~\ref{app:proof-coupled-sign}). Differentiating it also requires the state's response to $u$.
The next example starts from a spatially constant predator--prey equilibrium. The effort may vary in space at either step.
\begin{proposition}[A local minimizer in the coupled system]
\label{prop:coupled-local-minimum}
On any of the spatial grids defined in Section~\ref{sec:spatial-operator}, take $T=N=2$, $\mu_1,\mu_2>0$ and
\[
\begin{aligned}
 &m=\alpha=\beta=q=\lambda=1,\quad r=\tfrac12,\quad \eta=2,\\
 &U_{\max}=3,\quad \gamma_1=\gamma_2=0,\quad
 (v_1^0,v_2^0)=(\tfrac12\one,\one).
\end{aligned}
\]
The control $(u^0,u^1)=(0,2\one)$ is a strict constrained local minimizer of $j_0$ on the full control box. It has cost $-2$. The control $(2\one,\one)$ has lower cost $-33/14$. At $(0,2\one)$, decreasing the second-step effort uniformly gives a strict descent direction for $j_1$.
\end{proposition}
The proof is in Appendix~\ref{app:proof-coupled-sign}.
\section{Uniform derivative comparison and control stability}
\label{sec:uniform-derivatives-new}
On $\R^{N(M+1)}$, equipped with $\|\cdot\|_{h,\tau}$, use the \emph{operator norm}
\[
 \|A\|_{\mathcal L(\R^{N(M+1)})}
 \coloneqq \sup_{a\in\R^{N(M+1)}\setminus\{0\}}\frac{\|Aa\|_{h,\tau}}{\|a\|_{h,\tau}}.
\]
For $w=(w^n)_{n=0}^N\in\R^{(N+1)(M+1)}$, define
\[
\begin{aligned}
w^-&\coloneqq (w^n)_{n=0}^{N-1}\in\R^{N(M+1)},\\
w^+&\coloneqq (w^{n+1})_{n=0}^{N-1}\in\R^{N(M+1)}.
\end{aligned}
\]
Set
\[
D_\tau w^n\coloneqq \frac{w^{n+1}-w^n}{\tau},\quad
\|D_\tau w\|_{h,\tau}^2=\tau\sum_{n=0}^{N-1}\|D_\tau w^n\|_h^2.
\]
\begin{theorem}[A uniform $C^2$ quadrature estimate in one dimension]
\label{th:quadrature-C2-new}
On $\Omega=(0,1)$, fix $T$, the model parameters with $\mu_1,\mu_2>0$, and $U_{\max}$. For $i=1,2$, assume $0\le v_i^0\in\R^{M+1}$ and
\[
\sup_{\substack{h=1/M\\M\in\mathbf N}}
\bigl(\|v_i^0\|_\infty+\|v_i^0\|_h+|v_i^0|_{1,h}\bigr)
<\infty.
\]
For \eqref{eq:positive} and identical tracking terms in $j_0,j_1$, there exist $\bar\tau>0$, $C<\infty$, independent of $h$, $N$, and $u\in\Uk$, such that for $\tau\le\bar\tau$,
\begin{equation}\label{eq:C0-new}
|j_1(u)-j_0(u)|\le C\tau,
\end{equation}
and
\begin{equation}\label{eq:C1-new}
\|g_1(u)-g_0(u)\|_{h,\tau}\le C\tau,
\end{equation}
and
\begin{equation}\label{eq:C2-new}
\|H_1(u)-H_0(u)\|_{\mathcal L(\R^{N(M+1)})}\le C\tau.
\end{equation}
\end{theorem}
The proof in Appendix~\ref{app:proof-th-quadrature-C2-new} uses the one-dimensional $H^1$ embedding to bound the second variation. Remark~\ref{rem:explicit-constants-new} in Appendix~\ref{app:proof-explicit-constants} defines the constants $C_z,C_s$ and the step bound $\bar\tau$ used below.
\begin{corollary}[Global uniqueness for a sufficiently large effort penalty]
\label{cor:effort-dominance-new}
Under Theorem~\ref{th:quadrature-C2-new}, also assume
$\|v_i^d\|_{h,\tau}\le D_i$ uniformly. Define
\begin{equation}\label{eq:effort-bound-new}
B\coloneqq C_s\sum_{i=1}^2\gamma_i(K_i\sqrt T+D_i)
+2\eta qC_z+\eta qU_{\max}\sqrt T C_s.
\end{equation}
If $\lambda>B$ and $\tau\le\bar\tau$, both objectives for the positive scheme are strongly convex on $\Uk$, with modulus at least $\lambda-B$ uniformly in the mesh. Each has a unique global minimizer.
\end{corollary}
The bound $B$ is independent of $\lambda$ and $u$ (Appendix~\ref{app:proof-cor-effort-dominance-new}).
\begin{corollary}[Control stability under uniform coercivity]
\label{cor:control-stability-new}
Under Theorem~\ref{th:quadrature-C2-new}, take $\tau\le\bar\tau$, a common constant $C$ in \eqref{eq:C1-new}--\eqref{eq:C2-new}, and $u_0\in\Uk$ such that, for every $u\in\Uk$, $\langle g_0(u_0),u-u_0\rangle_{h,\tau}\ge0$. For $r>0$, set
\[
K_r\coloneqq \{u\in\Uk:\|u-u_0\|_{h,\tau}\le r\}.
\]
Suppose that for some $\kappa>0$, every $u\in K_r$ and every $a\in\R^{N(M+1)}$,
\begin{equation}\label{eq:uniform-coercivity-new}
\langle H_0(u)a,a\rangle_{h,\tau}\ge\kappa\|a\|_{h,\tau}^2.
\end{equation}
If $C\tau<\kappa\min\{1,r\}$, the flux problem has a unique stationary point $u_1$ in $K_r$, which is a strict local minimizer and satisfies
\begin{equation}\label{eq:O-tau-stability-new}
 \|u_1-u_0\|_{h,\tau}\le\frac{C\tau}{\kappa}.
\end{equation}
Stationarity means the box variational inequality, and the estimate is uniform in the mesh for fixed $r,\kappa$. With coercivity on the whole box, the minimizers are unique and global, with the same bound.
\end{corollary}
See Appendix~\ref{app:proof-cor-control-stability-new}.
At a common control and on a common subspace, corresponding ordered Hessian eigenvalues differ by at most $C\tau$ (Appendix~\ref{app:proof-spectral-perturbation}).
\section{Second-order conditions and local sweep theory}
Fix a grid and one discrete objective.
\subsection{Stationarity and the critical cone}
Define the \emph{stationarity residual} by
\begin{equation}\label{eq:residual}
 \mapdef{R}{\Uk}{\R^{N(M+1)}}{u}{u-\proj(u-g(u))}.
\end{equation}
Then $R(u^*)=0$ is equivalent to the following inequality for every $u\in\Uk$:
\begin{equation*}
\ipt{g(u^*)}{u-u^*}\ge0.
\end{equation*}
Write the control array as a vector with coordinates $i=1,\ldots,N(M+1)$. Coordinates at their bounds form the \emph{active sets}:
\begin{equation*}
A_0\coloneqq \{i:u_i^*=0\},\quad A_U\coloneqq \{i:u_i^*=U_{\max}\},
\end{equation*}
The other coordinates form the \emph{free set}:
\begin{equation*}
F\coloneqq \{i:0<u_i^*<U_{\max}\},\quad A\coloneqq A_0\cup A_U.
\end{equation*}
The \emph{tangent cone} and \emph{critical cone} in $\R^{N(M+1)}$ are
\begin{equation*}
T(u^*)\coloneqq\{d:d_i\ge0\text{ on }A_0,\ d_i\le0\text{ on }A_U\},
\end{equation*}
and
\begin{equation*}
\mathcal C(u^*)\coloneqq \{d\in T(u^*):\ipt{g(u^*)}d=0\}.
\end{equation*}
\emph{Strict complementarity} consists of the following two conditions, the first for every $i\in A_0$ and the second for every $j\in A_U$:
\begin{equation}\label{eq:strict}
g_i(u^*)>0,\qquad g_j(u^*)<0,
\end{equation}
and, at a stationary control, gives $$\mathcal C(u^*)=\{d\in\R^{N(M+1)}:d_i=0\text{ on }A\}.$$
\begin{theorem}[Fixed-grid second-order conditions]\label{th:second}
At a local minimizer, for every $d\in\mathcal C(u^*)$,
\begin{equation}
\ipt{H(u^*)d}d\ge0.
\end{equation}
Conversely, if $u^*$ is stationary and, for some $\kappa>0$ and every $d\in\mathcal C(u^*)$,
\begin{equation}\label{eq:coercive}
\ipt{H(u^*)d}d\ge\kappa\normt d^2,
\end{equation}
then, for some $c>0$ and all nearby feasible $u$, we have $j(u)\ge j(u^*)+c\normt{u-u^*}^2$.
\end{theorem}
See Appendix~\ref{app:local}.
\subsection{Relaxation inside the projection}
For $\theta>0$, the \emph{projected sweep} is
\begin{equation}\label{eq:fbs}
 \mapdef{\Phi_\theta}{\Uk}{\Uk}{u}{
 \proj\!\left(u-\frac\theta\lambda g(u)\right)}.
\end{equation}
Extend $\Phi_\theta$ to $\mathcal O_{h,\tau}$ by the same formula. For IMEX, the unprojected update is
\[
 u^n-\frac{g^n(u)}{\lambda}
 \byeq{eq:gradI}{=}\frac{q}{\lambda}v_2^n\circ(\eta+p_2^{n+1}).
\]
FBS is projected gradient descent with step $\theta/\lambda$, with relaxation applied before projection. Reversing this order changes the update on active coordinates. The choice $\theta\le1$ gives under-relaxation.
For finite index sets, vertical bars denote the number of elements. For $F=\{i_1<\cdots<i_{|F|}\}\ne\varnothing$, define insertion and projection by
\[
 \mapdef{E_F}{\R^{|F|}}{\R^{N(M+1)}}{x}{\sum_{k=1}^{|F|} x_ke_{i_k}},
\]
and
\[
 \mapdef{P_F}{\R^{N(M+1)}}{\R^{N(M+1)}}{d}{E_FE_F^Td},
\]
where $(e_i)_{i=1}^{N(M+1)}$ is the coordinate basis of $\R^{N(M+1)}$. Set
\begin{equation}\label{eq:free_types}
W_F\coloneqq E_F^T\Wht E_F\in\R^{|F|\times |F|},
\end{equation}
and
\begin{equation*}
H_F\coloneqq P_FH(u^*)P_F\big|_F\simeq E_F^TH(u^*)E_F\in\R^{|F|\times |F|}.
\end{equation*}
Here $\big|_F$ means restriction to the free-coordinate subspace. Use $\|x\|_{W_F}\coloneqq (x^TW_Fx)^{1/2}$ and set $P_F\coloneqq 0$ if $F=\varnothing$.
\begin{theorem}[Exact local FBS criterion]\label{th:fbs}
Suppose $u^*$ is stationary and strictly complementary. Then $\Phi_\theta$ is differentiable at $u^*$ and
\begin{equation}\label{eq:fullderivative}
D\Phi_\theta(u^*)=P_F\left(I-\frac\theta\lambda H(u^*)\right).
\end{equation}
Near $u^*$, one step fixes all active coordinates at their bounds \cite{BurkeMore1988}. If $F=\varnothing$, the map is locally constant. Otherwise, on the free affine space,
\begin{equation}\label{eq:Richardson}
D\Phi_\theta(u^*)\big|_F=I_F-\frac\theta\lambda H_F.
\end{equation}
If $H_F$ is positive definite, with smallest and largest weighted eigenvalues $0<\mu_{\min}\le\mu_{\max}$, the linear iteration factor on the free coordinates is
\begin{equation}\label{eq:rho}
\rho_\theta\coloneqq \max\left\{\left|1-\frac{\theta\mu_{\min}}\lambda\right|,
\left|1-\frac{\theta\mu_{\max}}\lambda\right|\right\}.
\end{equation}
It is below one exactly when
\begin{equation}\label{eq:threshold}
0<\theta<\theta_{\max}\coloneqq \frac{2\lambda}{\mu_{\max}}.
\end{equation}
For $0<\theta<\theta_{\max}$, convergence is locally $q$-linear once the active set is identified. The optimal positive relaxation is
\begin{equation}\label{eq:optimalrelax}
\theta_R\coloneqq \frac{2\lambda}{\mu_{\min}+\mu_{\max}},\qquad
\rho_{\theta_R}= \frac{\mu_{\max}-\mu_{\min}}{\mu_{\max}+\mu_{\min}}.
\end{equation}
For under-relaxation, use $\theta_{\rm opt}\coloneqq \min\{1,\theta_R\}$ in \eqref{eq:rho}. For $\theta>\theta_{\max}$, the fixed point is linearly unstable.
\end{theorem}
Here $q$-linear convergence means that each error norm is eventually bounded by a fixed factor smaller than one times the preceding error norm. Locally, FBS is Richardson iteration on the free coordinates (Appendix~\ref{app:local}). Global FBS convergence is studied in \cite{McAsey2012}.
\subsection{Matrix-free spectral diagnostics}
Set
\begin{equation}\label{eq:symmetricH}
\widetilde H_F\coloneqq W_F^{1/2}H_FW_F^{-1/2},
\end{equation}
which is Euclidean symmetric by weighted symmetry. For $y\in\R^{|F|}$, its action is
\begin{equation*}
\mapdef{\widetilde H_F}{\R^{|F|}}{\R^{|F|}}{y}{W_F^{1/2}E_F^TH(u^*)E_FW_F^{-1/2}y}.
\end{equation*}
For a matrix $A$, $\sigma(A)$ denotes its set of eigenvalues and $\|A\|_2$ the operator norm induced by the Euclidean vector norm. Positive diagonal matrix powers, such as $W_F^{1/2}$, are taken entrywise on the diagonal. We estimate the smallest and largest eigenvalues with Lanczos and mark the estimates with a hat. For nonempty $F$ and an approximate eigenpair $(\widehat\mu,y)$ with $\|y\|_2=1$, the \emph{Ritz residual} bounds the distance to the spectrum $\sigma(\widetilde H_F)$:
\[
\min_{\mu\in\sigma(\widetilde H_F)}|\widehat\mu-\mu|
\le\|\widetilde H_Fy-\widehat\mu y\|_2.
\]
See Appendix~\ref{app:proof-ritz-residual}.
\section{Algorithms and reproducibility}\label{sec:algorithms}
We accept a run when
\begin{equation}\label{eq:stop}
\|R(u)\|_\infty\le10^{-6}.
\end{equation}
The four methods are FBS, Anderson-accelerated FBS, spectral projected gradient (SPG), and the limited-memory BFGS method with bound constraints (L-BFGS-B). Appendix~\ref{app:algorithm-details} gives their iteration formulas and parameters. To compare their work across hardware, we count state and adjoint solves. The code used for the numerical experiments is available in~\cite{AutranCode}.
\section{Numerical results}\label{sec:numerical}
\subsection{Reference data and comparison protocol}
We use $\Omega\coloneqq(0,1)$, $(T,M,N)\coloneqq(4,80,120)$,
$(m,r,\alpha,\beta,q)\coloneqq(1,1.2,1,1,1)$,
$\mu_1\coloneqq\mu_2\coloneqq0.002$ and
\[
(U_{\max},\gamma_1,\gamma_2,\lambda,\eta)
\coloneqq(2,0.5,0.5,0.05,0.8),
\]
with initial profiles
\[
\begin{aligned}
v_1^0(x)&\coloneqq1.3e^{-(x-0.30)^2/0.02}+0.05,\\
v_2^0(x)&\coloneqq1.4e^{-(x-0.70)^2/0.02}+0.05,
\end{aligned}
\]
and targets
\[
\begin{aligned}
v_1^d(x)&\coloneqq r/\beta+0.15\cos(2\pi x),\\
v_2^d(x)&\coloneqq m/\alpha-0.15\cos(2\pi x).
\end{aligned}
\]
Unless stated otherwise, all runs start from zero. For baseline validation, we use IMEX with $M=80$, $N=1920$. For each tested $N\mid1920$, copy each control value onto its corresponding finer time intervals. The common objective is
\begin{equation*}
j^{\rm I}_{h,T/1920}\!\left(\bigl(u^{\lfloor kN/1920\rfloor}\bigr)_{k=0}^{1919}\right).
\end{equation*}
Here $\lfloor s\rfloor$ is the largest integer not exceeding $s$.
The derivative and Hessian tests recover the expected Taylor orders. Appendix~\ref{app:numerical-checks} gives these results alongside the weighted symmetry defects and positivity tests.
\subsection{How harvest quadrature changes the computed control}
\begin{figure*}[t]\centering
\includegraphics[width=.325\textwidth]{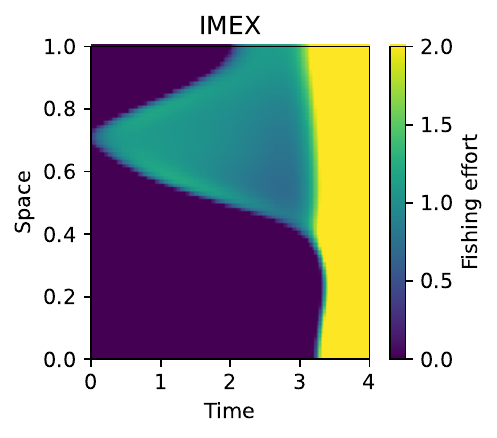}\hfill
\includegraphics[width=.325\textwidth]{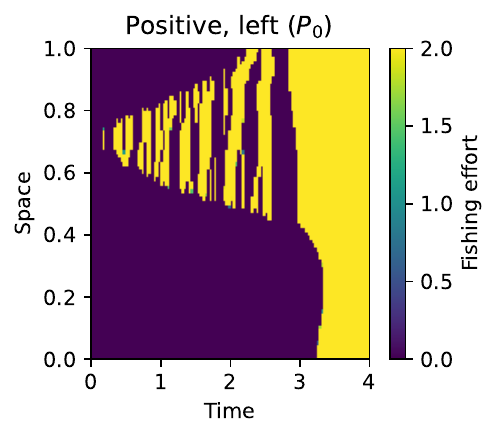}\hfill
\includegraphics[width=.325\textwidth]{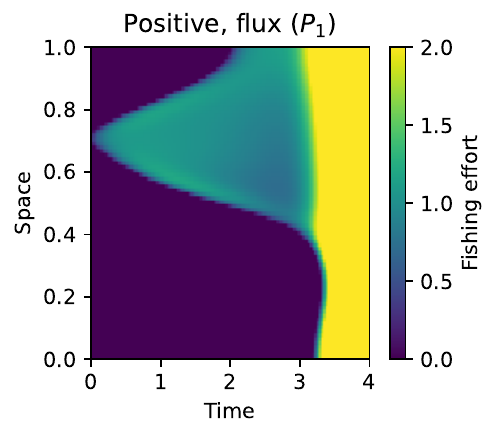}
\caption{Reference controls on a common color scale: I (left), $P_0$ (center), $P_1$ (right). The positive $P_0$ scheme produces temporal switching.}\label{fig:controls}
\end{figure*}

Table~\ref{tab:reference} and Fig.~\ref{fig:controls} show the controls reached from zero. A coordinate $i$ is counted as free when $10^{-8}<u_i<U_{\max}-10^{-8}$. Of the 9720 coordinates, 18 are free for $P_0$, compared with 3127 for I and 3076 for $P_1$. To measure switching, define the \emph{temporal variation} by
\begin{equation*}
{\rm TV}_t(u)\coloneqq\sum_{n=0}^{N-2}\sum_{j=0}^M w_j|u_j^{n+1}-u_j^n|.
\end{equation*}
$P_0$ switches sharply in time and has the lowest own-grid cost. On the common validation grid, I has the lowest cost, followed by $P_1$, then $P_0$.
\begin{table*}[t]\centering\small
\caption{Zero-start reference controls. Common validation uses the same IMEX objective with $N=1920$.}\label{tab:reference}
\begin{tabular}{lrrrrrr}
\toprule
Problem & Own-grid objective & Common validation & $\|R\|_\infty$ & $|F|$ & ${\rm TV}_t$ & $\min v_i$ \\
\midrule
I & -1.5491 & -1.553 & $9.95\times10^{-8}$ & 3127 & 2.252 & 0.02198 \\
$P_0$ & -1.5949 & -1.517 & $6.25\times10^{-7}$ & 18 & 19.11 & 0.02306 \\
$P_1$ & -1.521 & -1.5495 & $1.32\times10^{-7}$ & 3076 & 2.232 & 0.02305 \\
\bottomrule
\end{tabular}
\end{table*}
\subsection{Free Hessians and the local threshold}
\begin{figure*}[t]\centering
\includegraphics[width=.485\textwidth]{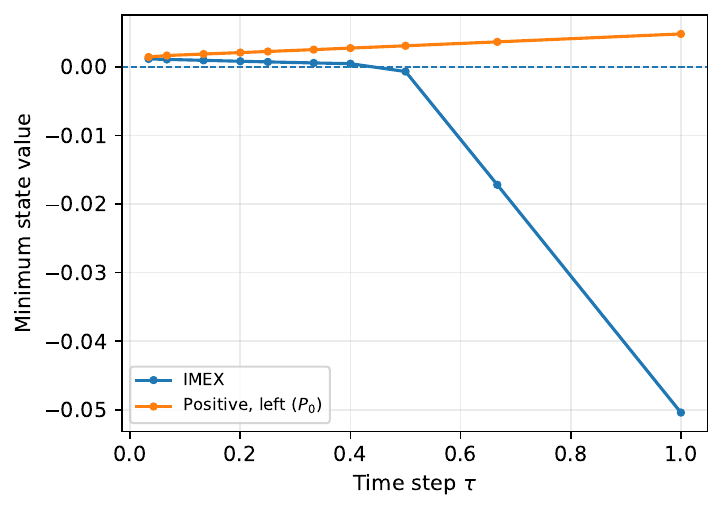}\hfill
\includegraphics[width=.485\textwidth]{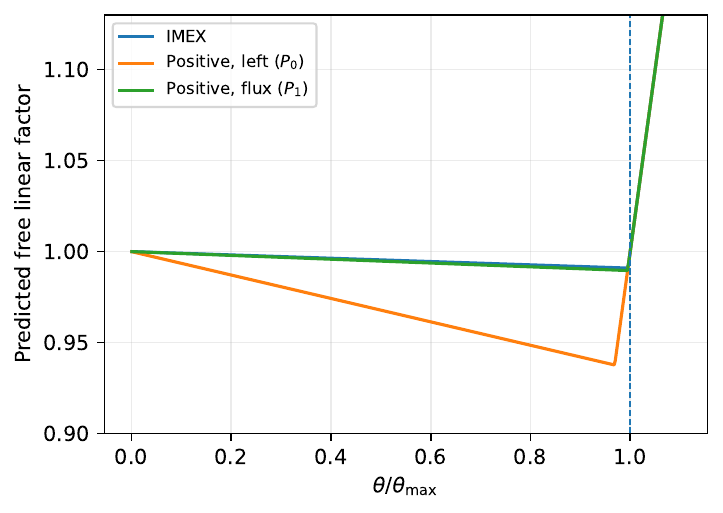}
\caption{Minimum state at fixed maximal effort (left), predicted free linear factor against normalized relaxation (right). }\label{fig:positivityspectrum}
\end{figure*}

The matrix-free Ritz estimates in Table~\ref{tab:spectra} put the relaxation threshold near $0.01$ for I and $P_1$. IMEX has unrelaxed free modes as low as about $-207.3$. For $P_0$, the threshold exceeds one, but the small $\widehat\mu_{\min}$ gives $\widehat\rho_\theta\approx0.9995$ at the shared $\theta=0.0093$ by \eqref{eq:rho}. Thus the local sweep is predicted to converge slowly near this stationary control (Fig.~\ref{fig:positivityspectrum}, right).
\begin{table*}[t]\centering\small
\caption{Ritz estimates of free spectra and relaxation thresholds. The last column gives the Euclidean residual for the smallest computed eigenpair of $\widetilde H_F$.}\label{tab:spectra}
\begin{tabular}{lrrrrrr}
\toprule
Problem & $\widehat\mu_{\min}$ & $\widehat\mu_{\max}$ & $\widehat\theta_{\max}$ & $\widehat\theta_{\rm opt}$ & $\widehat\rho_{\rm opt}$ & Smallest-mode residual \\
\midrule
I & 0.04737 & 10.42 & 0.009600 & 0.009557 & 0.9909 & $4.27\times10^{-8}$ \\
$P_0$ & 0.002752 & 0.08558 & 1.168 & 1 & 0.9450 & $1.70\times10^{-17}$ \\
$P_1$ & 0.05016 & 9.588 & 0.01043 & 0.01038 & 0.9896 & $4.92\times10^{-8}$ \\
\bottomrule
\end{tabular}
\end{table*}
For the largest computed mode $d$, the local directional test compares
\begin{equation*}
\frac{\Phi_\theta(u+\varepsilon d)-\Phi_\theta(u)}{\varepsilon}
\quad\text{with }d.
\end{equation*}
At $\theta=0.95\widehat\theta_{\max}$ and $1.05\widehat\theta_{\max}$, the predicted factors are $-0.9$ and $-1.1$, with magnitudes $0.9$ and $1.1$. All three problems match these predictions within $3\times10^{-4}$. For $P_0$, both tests use extrapolated steps, $\theta>1$.
\subsection{Algorithm comparison and rejected runs}
All four methods pass the acceptance test for I and $P_1$, with L-BFGS-B requiring the fewest state and adjoint solves. The two sweeps reach the iteration limit for $P_0$. SPG and L-BFGS-B reach stationary controls, but with different costs. Appendix~\ref{app:algorithm-comparison} gives the parameters, final residuals and solve counts for every run.

\begin{figure*}[t]\centering
\includegraphics[width=\textwidth]{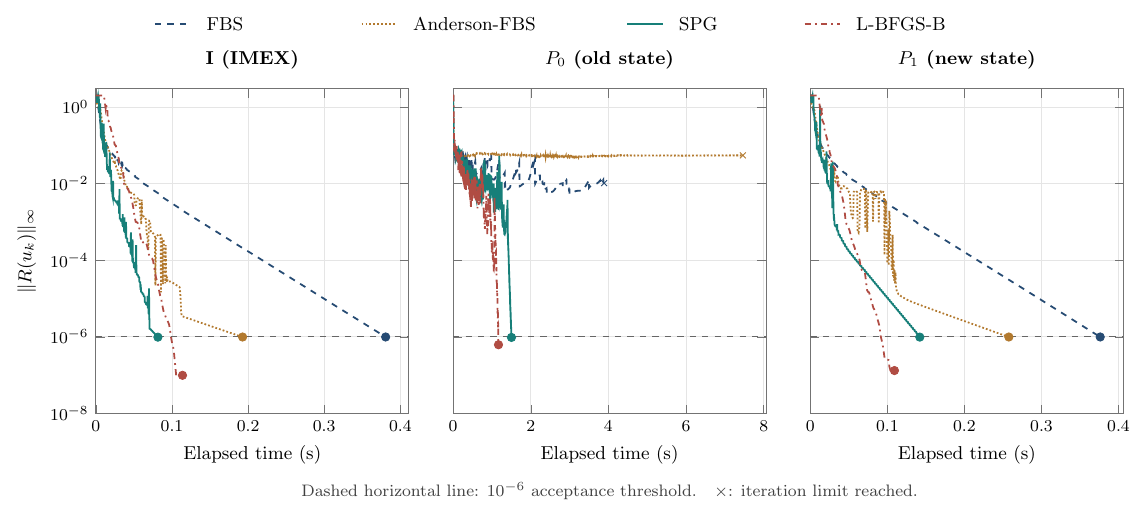}
\caption{Stationarity residual against elapsed time for the reference problem with $M=80$, $N=120$ and zero initial control. Each curve shows the run with median total runtime among five sequential repetitions. Horizontal scales differ between panels. The dashed horizontal line marks the $10^{-6}$ acceptance threshold. Crosses mark runs that reached the iteration limit without satisfying this threshold.}\label{fig:algorithm-timing}
\end{figure*}

Figure~\ref{fig:algorithm-timing} compares the residual histories against elapsed time on the reference grid $M=80$, $N=120$, starting from zero. We repeat each method--problem pair five times, with all runs executed sequentially and thread limits set to one. A monotonic clock measures each solver call, including final diagnostics. JIT compilation is warmed up beforehand. We perform the independent residual check after timing. Each curve is the run with median total runtime, rather than an average of trajectories.

Table~\ref{tab:algorithm-timing} reports total runtimes under each method's stopping rules from Appendix~\ref{app:algorithm-details}. These durations need not equal the first time the residual crosses the common $10^{-6}$ acceptance threshold. The runs used ARM64 macOS~27.2.1, Python~3.12.14, NumPy~2.3.5, SciPy~1.17.0 and Numba~0.65.1. SPG has the shortest median runtime for I. L-BFGS-B is fastest for $P_1$ and among the accepted $P_0$ runs. The two accepted $P_0$ methods reach different objective values, as reported in Appendix~\ref{app:algorithm-comparison}.
\begin{table*}[t]\centering\small
\caption{Total runtimes in seconds over five repetitions, shown as median [minimum, maximum]. A dagger marks runs that reached 12000 iterations without passing the acceptance test.}\label{tab:algorithm-timing}
\begin{tabular}{lccc}
\toprule
Method & I & $P_0$ & $P_1$ \\
\midrule
FBS & 0.3808 [0.3803, 0.4438] & 3.888 [3.878, 4.574]$^{\dagger}$ & 0.3764 [0.3716, 0.4446] \\
Anderson-FBS & 0.1929 [0.1917, 0.212] & 7.466 [7.445, 8.13]$^{\dagger}$ & 0.2577 [0.2562, 0.2833] \\
SPG & 0.08162 [0.08018, 0.09635] & 1.503 [1.497, 1.743] & 0.1421 [0.141, 0.167] \\
L-BFGS-B & 0.1139 [0.1105, 0.1185] & 1.168 [1.165, 1.253] & 0.1091 [0.1076, 0.1141] \\
\bottomrule
\end{tabular}
\end{table*}

\subsection{Control differences under temporal refinement}
\begin{figure*}[t]\centering
\includegraphics[width=.485\textwidth]{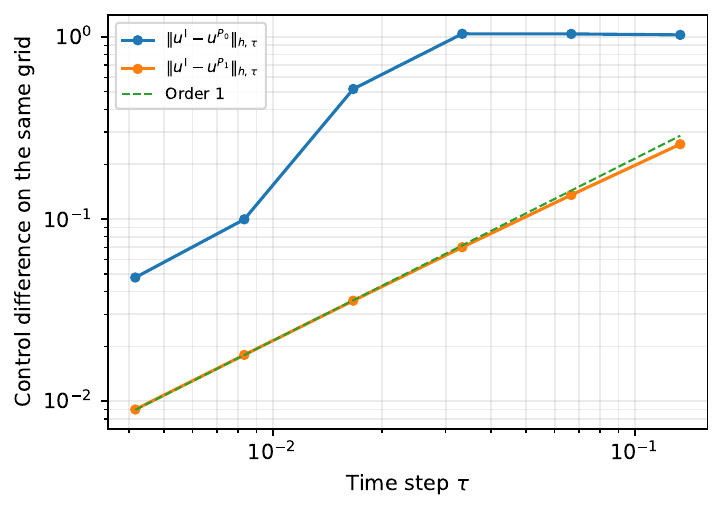}\hfill
\includegraphics[width=.485\textwidth]{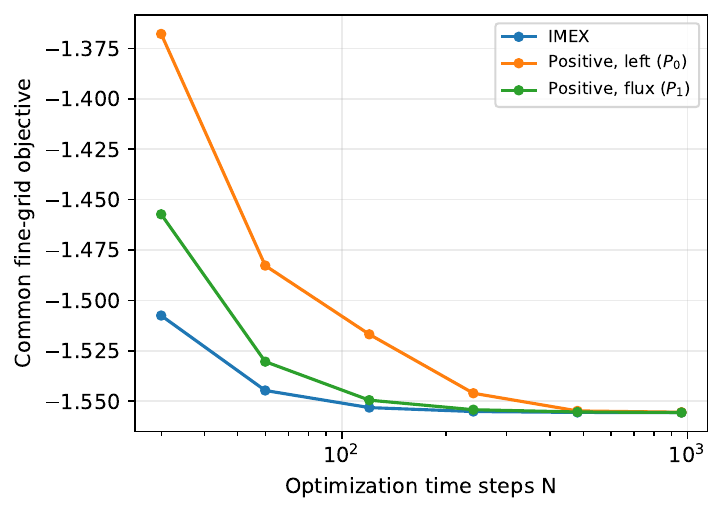}
\caption{Same-grid control differences from I (left), common $N=1920$ validation (right). The differences shrink on the finer time grids.}\label{fig:refinement}
\end{figure*}

With $M=80$ fixed, we refine through $N=30,\allowbreak60,\allowbreak120,\allowbreak240,\allowbreak480,\allowbreak960$. The first run starts from zero at $N=30$. Later runs start from the previous control, copied to the two corresponding finer intervals. At $N=120$, we replace that control with the reference control computed from zero before continuing the refinement. For candidates $u_N^{\rm I},u_N^{P_\zeta}\in\R^{N(M+1)}$ on the same grid, set
\begin{equation*}
D_\zeta(N)\coloneqq\normt{u_N^{\rm I}-u_N^{P_\zeta}}.
\end{equation*}
All controls pass the residual check. The $P_1$ distance decreases with ratios tending to two. The $P_0$ distance follows a different path: it stays near one through $N=120$, then drops on finer grids (Table~\ref{tab:grid} and Fig.~\ref{fig:refinement}).
\begin{table*}[t]\centering\small
\caption{Temporal refinement at $M=80$. All candidates pass the independent residual check.}\label{tab:grid}
\begin{tabular}{rrrrrrr}
\toprule
$N$ & $j^{\rm I}$ & $j^{P_0}$ & $j^{P_1}$ & $D_0$ & $D_1$ & $\|R^{P_0}\|_\infty$ \\
\midrule
30 & -1.532 & -1.743 & -1.423 & 1.025 & 0.2574 & $0$ \\
60 & -1.543 & -1.651 & -1.487 & 1.038 & 0.1355 & $1.30\times10^{-8}$ \\
120 & -1.5491 & -1.5949 & -1.521 & 1.039 & 0.07003 & $6.25\times10^{-7}$ \\
240 & -1.5525 & -1.566 & -1.538 & 0.5173 & 0.03564 & $1.33\times10^{-7}$ \\
480 & -1.554 & -1.561 & -1.547 & 0.09960 & 0.01796 & $1.17\times10^{-7}$ \\
960 & -1.555 & -1.558 & -1.5516 & 0.04775 & 0.00901 & $2.04\times10^{-7}$ \\
\bottomrule
\end{tabular}
\end{table*}
For state self-convergence at a fixed control, Appendix~\ref{app:state-self-convergence} separates time and space refinement.
\subsection{Coarse-grid optimization and multistart}
At $N=4$, the accepted IMEX control gives negative states. Both positive problems keep their states nonnegative, yet their common validation costs fall on opposite sides of I: higher for $P_0$, lower for $P_1$ (Table~\ref{tab:coarse}). The ordering changes with $N$.
\begin{table*}[!t]\centering\small
\caption{Coarse-grid candidates and common fine-grid validation. All candidates meet the residual tolerance.}\label{tab:coarse}
\begin{tabular}{rlrrrr}
\toprule
$N$ & Problem & $\min v_i$ & Own-grid objective & Common validation & $\|R\|_\infty$ \\
\midrule
4 & I & -0.4039 & -1.654 & 0.4478 & $6.41\times10^{-8}$ \\
4 & $P_0$ & 0.02653 & -1.852 & 1.255 & $1.28\times10^{-9}$ \\
4 & $P_1$ & 0.02710 & -0.7576 & -0.2562 & $2.86\times10^{-8}$ \\
8 & I & 0.01021 & -1.523 & -0.7813 & $5.57\times10^{-8}$ \\
8 & $P_0$ & 0.02723 & -2.049 & 0.1542 & $5.85\times10^{-10}$ \\
8 & $P_1$ & 0.02903 & -1.113 & -0.7170 & $4.67\times10^{-8}$ \\
20 & I & 0.01911 & -1.525 & -1.440 & $1.07\times10^{-7}$ \\
20 & $P_0$ & 0.02541 & -1.825 & -1.158 & $1.45\times10^{-16}$ \\
20 & $P_1$ & 0.02543 & -1.362 & -1.347 & $1.02\times10^{-7}$ \\
\bottomrule
\end{tabular}
\par\medskip
\caption{Results from five initial guesses. Distances are weighted $L^2$ distances to the zero-start control.}\label{tab:multistart}
\begin{tabular}{lrrrr}
\toprule
Problem & Accepted starts & Objective spread & Max. distance to zero start & Largest residual \\
\midrule
I & 5/5 & $2.66\times10^{-15}$ & $3.08\times10^{-7}$ & $4.26\times10^{-7}$ \\
$P_0$ & 5/5 & $8.02\times10^{-4}$ & $1.42\times10^{0}$ & $9.82\times10^{-7}$ \\
$P_1$ & 5/5 & $1.78\times10^{-15}$ & $2.00\times10^{-7}$ & $3.19\times10^{-7}$ \\
\bottomrule
\end{tabular}
\end{table*}
The five initial guesses are $0$, $1$, $U_{\max}$, independent uniform samples on $[0,U_{\max}]$, and $U_{\max}\mathbf1_{[0,T/2)}(t)$. Every control passes the residual check. For I and $P_1$, the costs differ only by rounding, and weighted distances to the zero-start control stay below $3.1\times10^{-7}$. The $P_0$ runs reach distinct controls, with a spread of about $8.02\times10^{-4}$ in cost and 1.425 in distance (Table~\ref{tab:multistart}).
\raggedbottom
\section{Targeted tests of the quadrature mechanism}\label{sec:targetedtests}
The next tests vary one parameter at a time. They also reverse the direction of continuation and separate time refinement from space refinement.
\subsection{Exact global verification in the removal model}\label{sec:reduced-certificate}
Take $\tau=q=y^0=\lambda=1$, $\eta=2$, $U_{\max}=3$, and $(x,z)\coloneqq(u^0,u^1)\in[0,3]^2$. The two costs are
\[
 \mapdef{F_0}{[0,3]^2}{\R}{(x,z)}{
 \tfrac12(x^2+z^2)-2\left(x+\frac z{1+x}\right)},
\]
and
\[
 \mapdef{F_1}{[0,3]^2}{\R}{(x,z)}{
 \tfrac12(x^2+z^2)+\frac2{(1+x)(1+z)}-2}.
\]
The old-state global minimizer is near
\[
(1.822,\,0.7087),
\] with cost $-2.235$. The extra strict local minimum $(0,2)$ has cost $-2$ (Appendix~\ref{app:proof-global-certificate}). The unique flux minimizer is $(t,t)$, $t\approx0.5437$, with cost $-0.8651$.
For $T=2$ and $K_c\coloneqq\{u\in L^2(0,2):0\le u\le3\text{ a.e.}\}$, the common exact continuous cost is
\[
 \mapdef{j_c}{K_c}{\R}{u}{
 \frac12\int_0^2u(t)^2\,\mathrm{d}t-2\left(1-e^{-\int_0^2u(t)\,\mathrm{d}t}\right)}.
\]
The continuous cost has the unique minimizer $u\equiv s$, $s\approx0.6011$, with value $-1.038$. For two-step controls, $j_c(x,z)=\frac12(x^2+z^2)-2+2e^{-(x+z)}$: it takes values $0.07021$ and $-1.030$ at the certified old-state and flux global minimizers. The quadrature bias is visible in these values.

\begin{figure*}[t]\centering
\includegraphics[width=.65\textwidth]{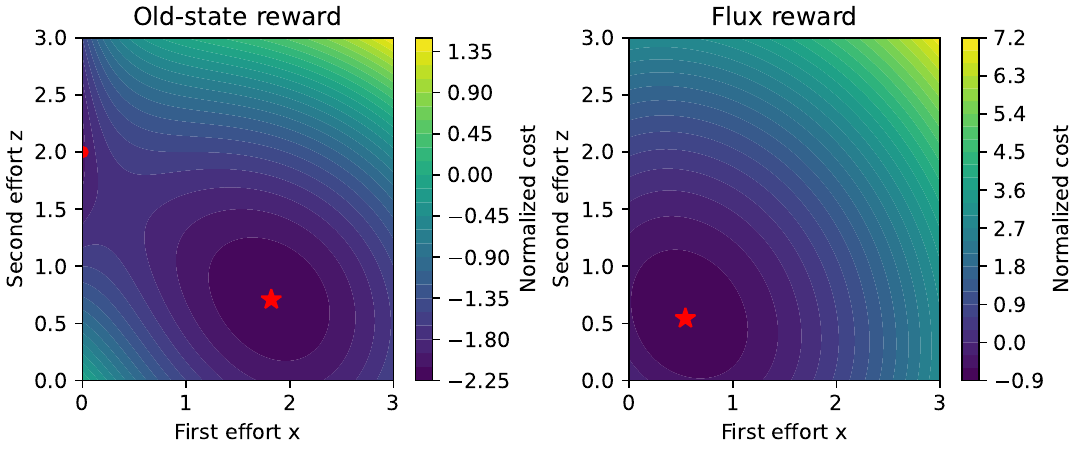}\hfill
\includegraphics[width=.33\textwidth]{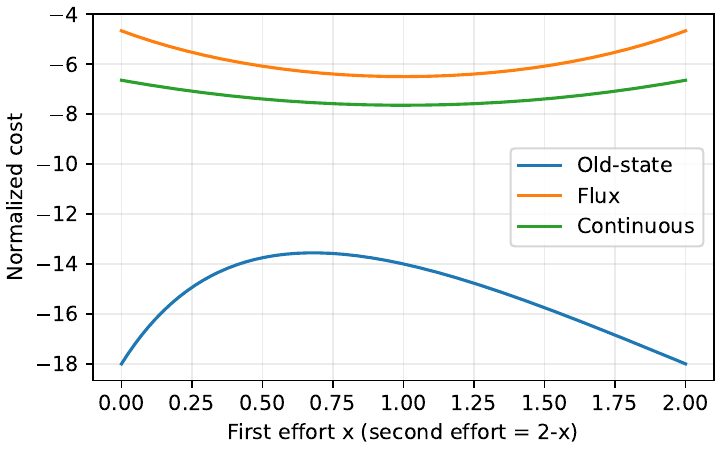}
\caption{Pure removal. Left: two-step costs at $\eta=2$. Red markers identify old-state global and extra local minima and the unique flux minimum. Right: at $\eta=10$ and fixed total effort, old-state quadrature favors extreme allocations, while the flux and continuous costs favor equal allocation.}\label{fig:reducedglobal}
\end{figure*}

With $\eta=10$, $u^0+u^1=2$ and the other coefficients unchanged, Proposition~\ref{prop:budget} gives old-state minimizers $(0,2),(2,0)$, while the flux and continuous costs share the unique allocation $(1,1)$ (Fig.~\ref{fig:reducedglobal}).
\subsection{A common reference and acceptance protocol}
For $u\in\Uk$, spatial interpolation gives $U_{\rm ref}(t)\in \R^{M_{\rm ref}+1}$, constant on each original time interval. Set $h_{\rm ref}\coloneqq1/M_{\rm ref}$, $M_{\rm ref}\coloneqq80$, and integrate
\begin{align*}
 \dot V_1&=\mu_1L_{h_{\rm ref}}V_1+(-m\mathbf1+\alpha V_2)\circ V_1,\\
 \dot V_2&=\mu_2L_{h_{\rm ref}}V_2+(r\mathbf1-\beta V_1-qU_{\rm ref})\circ V_2,
\end{align*}
where $V_i:[0,T]\to \R^{M_{\rm ref}+1}$ and $\dot V_i$ denotes its time derivative. We sample the analytic initial and target profiles at the reference nodes, then integrate the four cost terms alongside the state equations. DOP853 restarts at every control jump.
With relative and absolute tolerances of $10^{-9}$ and $10^{-11}$, this gives a continuous-time reference on the shared spatial grid. We check the solver against the exact exponential solution for zero predators and piecewise constant harvesting.
Of $94$ final targeted controls, $93$ satisfy \eqref{eq:stop}. We retry failed runs with L-BFGS memory eight instead of 38.
\subsection{Curvature at common controls}
\begin{figure*}[t]
\centering
\includegraphics[width=\textwidth]{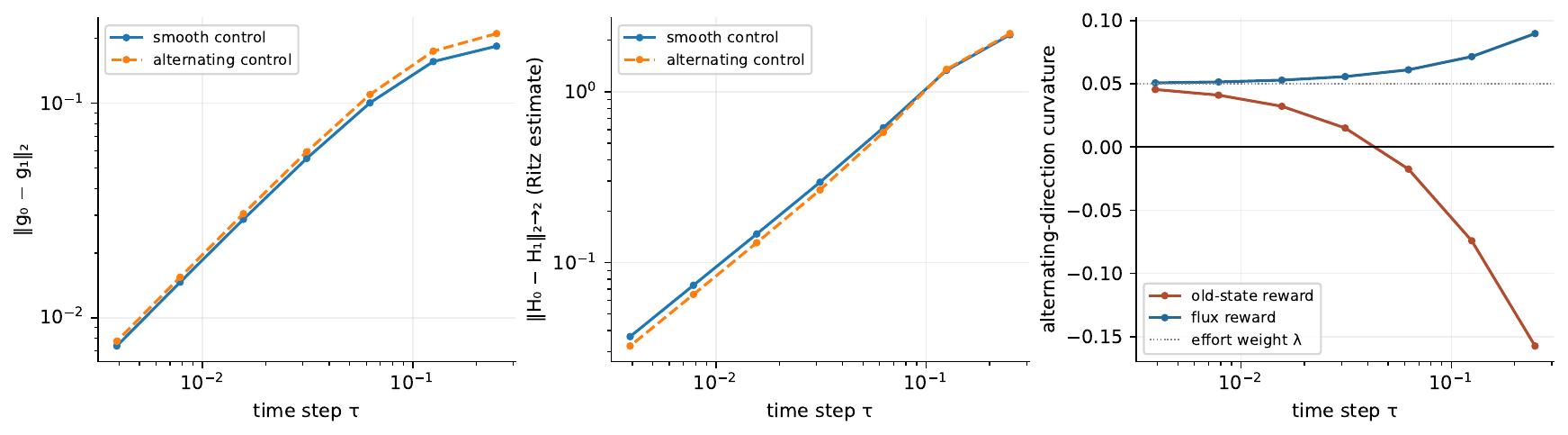}
\caption{Quadrature perturbations at common controls, $M=32$. Left and middle: weighted gradient differences and largest absolute Ritz values of the mass-scaled Hessian difference. Right: normalized alternating-direction curvature at the smooth control.}
\label{fig:researchcurvature}
\end{figure*}

Both positive problems produce the same state at a fixed control. On $M=32$, take
\[
u^n(x)\coloneqq0.7+0.2\cos(2\pi x)\cos(2\pi t_n/T).
\]
The alternating direction $d\in\R^{N(M+1)}$, with entries $d_j^n\coloneqq(-1)^n/\sqrt T$, has $\normt d=1$ because $\sum_{j=0}^M w_j=1$. Exact Hessian actions give
\begin{center}
\begin{tabular}{rcc}
\toprule
$N$ & $\ipt{d}{H_0d}$ & $\ipt{d}{H_1d}$\\
\midrule
16 & $-0.1575$ & $0.08976$\\
32 & $-0.07423$ & $0.07153$\\
64 & $-0.01750$ & $0.06108$\\
128 & $0.01502$ & $0.05571$\\
1024 & $0.04551$ & $0.05074$\\
\bottomrule
\end{tabular}
\end{center}
Along this direction, $P_0$ has negative curvature on coarse grids, while $P_1$ has positive curvature (Fig.~\ref{fig:researchcurvature}, right). Refinement brings both towards $\lambda=0.05$.
At these controls, measure $\normt{g_1(u)-g_0(u)}$ and set
\[
\begin{aligned}
B(u)&\coloneqq \Wht^{1/2}(H_0(u)-H_1(u))\Wht^{-1/2}\\
&\in\R^{N(M+1)\times N(M+1)}.
\end{aligned}
\]
For symmetric $B(u)$, let $\widehat\nu$ be the Ritz value with largest magnitude. We estimate the norm $$\|B(u)\|_2=\|H_1(u)-H_0(u)\|_{\mathcal L(\R^{N(M+1)})}$$ by $|\widehat\nu|$ and, for $\|y\|_2=1$, record the residual $\|B(u)y-\widehat\nu y\|_2$. The tests cover $M=8,16,32$, $N=16,32,\ldots,1024$, using both the smooth control and the spatially constant alternating control $u^n=0.7+0.2(-1)^n$.

At $M=32$, a fit over $N=256,512,1024$ gives gradient-difference rates $0.981$, $0.988$ and Hessian-estimate rates $0.997$, $1.004$, respectively. At $N=1024$, division by $\tau$ gives gradient differences $1.888$, $1.987$ and Hessian estimates $9.444$, $8.320$. The measured rates agree with Theorem~\ref{th:quadrature-C2-new}.
\subsection{Continuation and parameter changes}
\begin{figure*}[t]
\centering
\includegraphics[width=\textwidth]{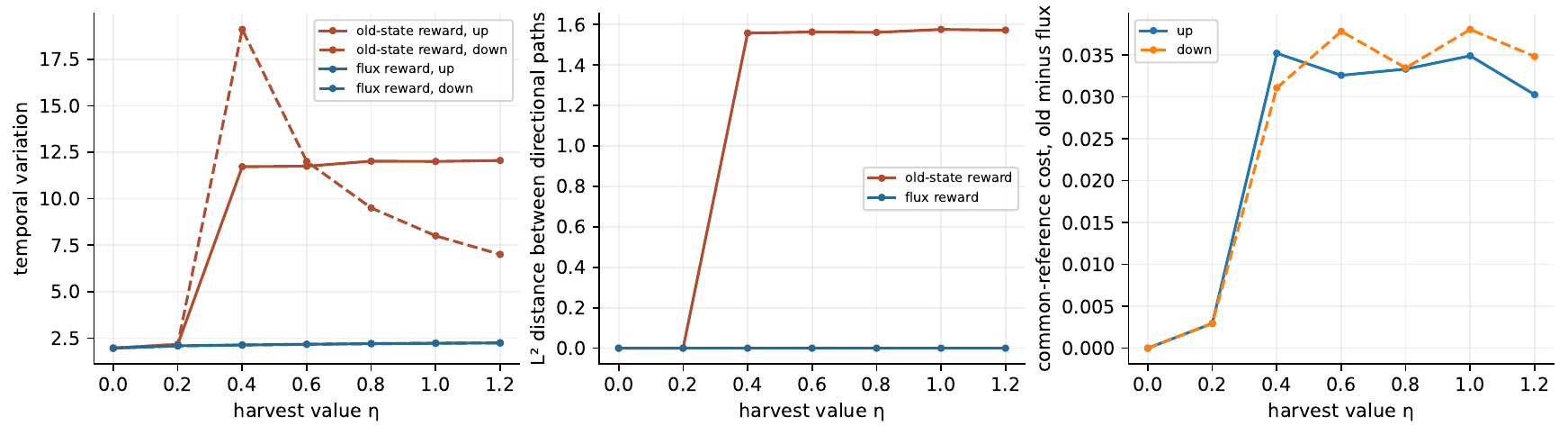}
\caption{Opposite continuation directions on $M=20,N=60$, each starting from zero: temporal variation, same-parameter control distance and common-reference cost difference.}
\label{fig:researchpaths}
\end{figure*}

On $M=20,N=60$, we increase $\eta$ through $0,0.2,\ldots,1.2$, then run a separate sequence in reverse. Each sequence starts from zero and passes its accepted control to the next run. The $P_0$ controls from opposite directions have $L^2$ distances of $1.558$--$1.576$ at $\eta=0.4,\ldots,1.2$, but agree within about $1.1\times10^{-6}$ at $\eta=0,0.2$. All $P_1$ distances stay below $4.9\times10^{-7}$. For the old-state rule, the direction of continuation changes the control reached (Fig.~\ref{fig:researchpaths}).
At $\eta=0.8$, we start from zero, one, the upper bound, a seeded uniform random control, and the upper bound until $0.45T$ followed by zero. All ten runs pass the residual test. The $P_0$ runs reach distinct controls. Their costs range from $-1.653$ to $-1.648$ and their temporal variations from $6.501$ to $19.61$. The five $P_1$ runs give numerically indistinguishable controls. Their costs agree within $9\times10^{-16}$ at $-1.484$, and their temporal variations are about $2.199$.
Separate series vary $\eta=0,0.2,0.4,0.8,1.2$, $\lambda=0.02,0.05,0.1,0.2$ and $T=1,2,4,6$. In the horizon series, $N=15T$ keeps $\tau=1/15$ fixed. The rewards agree at $\eta=0$. With positive $\eta$ and a zero initial guess, the $P_1$ control has a lower common-reference cost than the $P_0$ control in every case tested in these series. The differences range from $0.000502$ to $0.07664$. At $T=6$, temporal variation is $24.57$ for $P_0$ and $2.341$ for $P_1$.
\subsection{Separate temporal and spatial refinement}
\begin{figure*}[t]
\centering
\includegraphics[width=\textwidth]{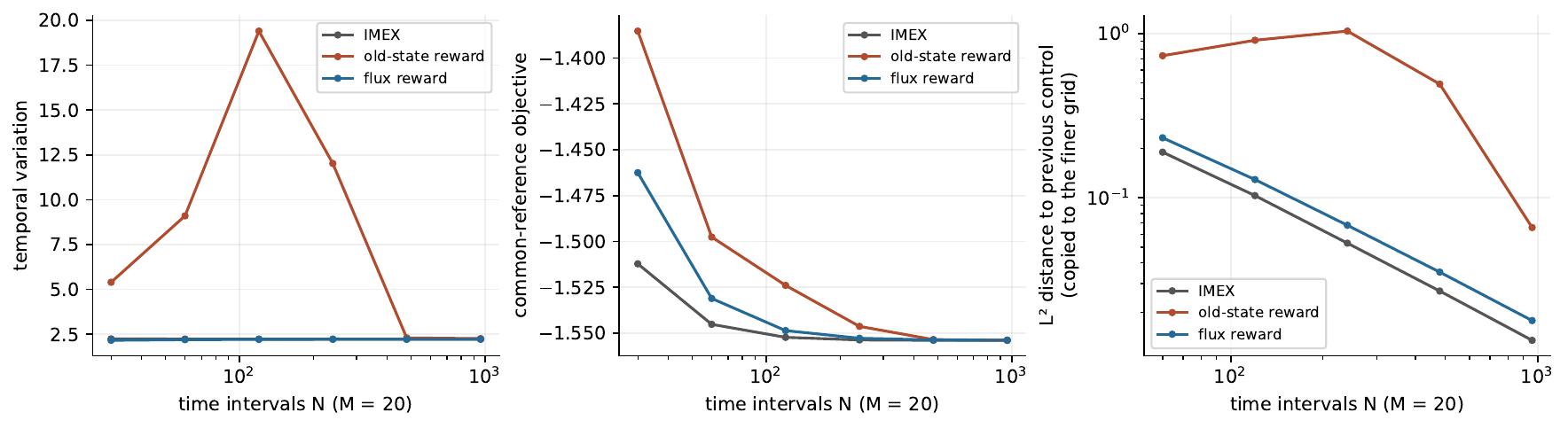}
\caption{Temporal refinement at $M=20$, with each run starting from the previous accepted control. Right: weighted $L^2$ distance between consecutive controls, after copying each coarse value to its two finer intervals. The old-state distance first grows, then drops sharply for $N\ge480$. All controls pass the residual test.}
\label{fig:researchrefinement}
\end{figure*}

With $M=20$ fixed, we refine through $N=30,\allowbreak60,\allowbreak120,\allowbreak240,\allowbreak480,\allowbreak960$, starting each run from the previous accepted control (Fig.~\ref{fig:researchrefinement}). The $P_0$ temporal variations are
\[
5.400,\quad9.100,\quad19.40,\quad12.02,\quad2.302,\quad2.261,
\]
while the variation for $P_1$ moves gradually from $2.164$ to $2.225$. On these six grids, the $L^2$ distance between the $P_0$ and $P_1$ controls is $0.949,\allowbreak 1.011,\allowbreak 1.029,\allowbreak 0.513,\allowbreak 0.0963,\allowbreak 0.0455$, respectively. The distance peaks at $N=120$ and then decreases. At $N=960$, the costs on the common reference grid ($M_{\rm ref}=80$) are approximately $-1.55389$ for IMEX, $-1.55388$ for $P_0$ and $-1.55384$ for $P_1$. On this grid, the $P_0$ cost is about $4\times10^{-5}$ below the $P_1$ cost.

For spatial refinement we fix $N=120$ and use $M=10,20,40,80$. The first and last $P_0$ controls have temporal variations $16.48$ and $27.87$. The last is a rejected candidate. At $M=80$, restarting with L-BFGS memory eight reduces its residual from $0.01886$ to $1.161\times10^{-6}$, just above the acceptance threshold. All $P_1$ controls pass, with temporal variations between $2.221$ and $2.232$. Thus temporal switching remains visible in the computed $P_0$ controls at this time step.

For the three candidates at $M=80,N=120$, tightening the reference solver's relative tolerance from $10^{-8}$ to $10^{-10}$ changes the costs by at most about $8.3\times10^{-12}$. Increasing $M_{\rm ref}$ from $80$ to $160$ changes them by about $2.22\times10^{-4}$ for IMEX, $1.03\times10^{-3}$ for the rejected $P_0$ candidate and $2.48\times10^{-4}$ for $P_1$. Reference-grid sensitivity of the cost differences at $M=20,N=960$ remains untested.
\section[Extension to production--destruction networks]{Extension to\\production--destruction networks}
\label{sec:pd-network}
A single flux can account for both donor loss and receiver gain in a network \cite{BlanesIserlesMacnamara2022}. With rates fixed within each step, this shared-flux scheme preserves positivity and records harvest exactly. Let $s,d\in\mathbf N$ count species and spatial nodes, with indices $1\le i,j\le s$, $1\le a\le d$. For positive weights $c_i,w_a$, set
\[
y\coloneqq(y_{i,a})\in\R^{sd},\quad\rho\coloneqq(c_iw_a)\in(0,\infty)^{sd}.
\]
Define the \emph{weighted mass} by
\[
 \mapdef{\mathcal M}{\R^{sd}}{\R}{y}{\rho^Ty}.
\]
Here $|y|$ is the componentwise absolute value, and $\R_+^{sd}\coloneqq[0,\infty)^{sd}$ is the set of nonnegative vectors. Let $\mathcal V\subset\R^{d_u}$, $d_u\in\mathbf N$, contain the admissible controls. At step $0\le n<N$, take $y^n\in\R_+^{sd}$, $u^n\in\mathcal V$ and finite rates $k_{ij,a}^n\ge0$ for transfer from $j$ to $i\ne j$, evaluated at $(y^n,u^n)$. The donor loses and the receiver gains, respectively $k_{ij,a}^ny_{j,a}^{n+1}$ and $(c_j/c_i)k_{ij,a}^ny_{j,a}^{n+1}$. Weighting by $c_jw_a$ and $c_iw_a$ makes the loss and gain equal. Define $K^n\in\R^{sd\times sd}$ as follows. For $i\ne j$, set
\begin{equation}\label{eq:pd-transfer}
K_{(i,a),(j,a)}^n\coloneqq(c_j/c_i)k_{ij,a}^n,
\end{equation}
and
\begin{equation*}
 K_{(j,a),(j,a)}^n\coloneqq-\sum_{\substack{i=1\\i\ne j}}^s k_{ij,a}^n,
\end{equation*}
with zero entries between different spatial nodes, so $\rho^TK^n=0$. Let $L_h\in\R^{d\times d}$ have nonnegative off-diagonal entries and $w^TL_h=0$, $w\coloneqq(w_a)\in(0,\infty)^d$. For $\mu_i\ge0$, set
\[
 G^n\coloneqq\operatorname{blockdiag}_{i=1}^s(\mu_iL_h)+K^n.
\]
Here $\operatorname{blockdiag}$ places the listed matrices on the diagonal and zeros elsewhere. Then $G^n$ is \emph{Metzler} (nonnegative off-diagonal entries), with $\rho^TG^n=0$. For mortality, harvesting and production $m^n,\allowbreak h^n,\allowbreak p^n\in\R_+^{sd}$, use
\begin{equation}\label{eq:pd-step}
 A^ny^{n+1}=y^n+\tau p^n,
\end{equation}
where
\begin{equation*}
A^n\coloneqq I_{sd}-\tau G^n+\tau\operatorname{diag}(m^n+h^n)\in\R^{sd\times sd}.
\end{equation*}
The rates may depend on $(y^n,u^n)$. We require them to remain finite even when a state component is zero.
\begin{theorem}[Positive transfer network and harvest accounting]
\label{th:pd-network}
For every $\tau>0$, nonnegative $y^n,p^n,m^n,h^n$ and finite nonnegative transfer rates, \eqref{eq:pd-step} has a unique nonnegative solution. It satisfies the exact mass balance
\begin{equation}\label{eq:pd-mass}
\begin{aligned}
&\mathcal M(y^{n+1})-\mathcal M(y^n)\\
&\quad=\tau\rho^Tp^n-\tau\rho^T\operatorname{diag}(m^n+h^n)y^{n+1}.
\end{aligned}
\end{equation}
With $A^n$ fixed at these rates, for every $z\in\R^{sd}$,
\begin{equation}\label{eq:pd-resolvent}
 \rho^T|(A^n)^{-1}z|\le\rho^T|z|.
\end{equation}
If $\rho^Tp^n\le\sigma\mathcal M(y^n)+b_n$, with $\sigma,b_n\ge0$, then, for $t_n=n\tau$,
\begin{equation}\label{eq:pd-growth}
\mathcal M(y^n)\le e^{\sigma t_n}
\left(\mathcal M(y^0)+\tau\sum_{r=0}^{n-1}b_r\right).
\end{equation}
If $p^n=0$ at every step, then
\[
 \tau\sum_{n=0}^{N-1}\rho^T\operatorname{diag}(h^n)y^{n+1}\le\mathcal M(y^0).
\]
\end{theorem}
The proof is in Appendix~\ref{app:proof-th-pd-network}.
For Lotka--Volterra, $s=2$, $d=M+1$, $c_1=\beta$, $c_2=\alpha$, and $k_{12,a}^n\coloneqq\beta v_{1,a}^n$ and $k_{21,a}^n\coloneqq0$. We also set $m^n\coloneqq(m\mathbf1,0)$, $h^n\coloneqq(0,qu^n)$, and $p^n\coloneqq(0,rv_2^n)$.
These choices turn \eqref{eq:pd-step} into \eqref{eq:positive}. An acyclic transfer graph can be ordered into a block triangular system. Cycles require a coupled solve.
\subsection{Reward accounting}
For economic weights $\ell\in\R_+^{sd}$ and $\zeta\in\{0,1\}$, define
\begin{equation}\label{eq:pd-reward}
\mapdef{\mathcal R_\zeta}{\mathcal V^N}{\R}{u}{
 \tau\sum_{n=0}^{N-1}\ell^T\operatorname{diag}(h^n)\cdot\big((1-\zeta)y^n+\zeta y^{n+1}\big)}.
\end{equation}
For $\ell=\rho$, $\mathcal R_1$ is the controlled mass loss in \eqref{eq:pd-mass}. Other weights assign different prices to those removals. For a common cost $\mathcal C:\mathcal V^N\to\R$, set
\[
 \mapdef{j_\zeta}{\mathcal V^N}{\R}{u}{\mathcal C(u)-\eta\mathcal R_\zeta(u)}.
\]
Then
\begin{equation}\label{eq:pd-defect}
\begin{aligned}
 j_0(u)-j_1(u)&=\eta\bigl(\mathcal R_1(u)-\mathcal R_0(u)\bigr)\\
 &\byeq{eq:pd-reward}{=}
 \eta\tau\sum_{n=0}^{N-1}\ell^T\operatorname{diag}(h^n)(y^{n+1}-y^n).
\end{aligned}
\end{equation}
Production and incoming transfers can change its sign. For one donor, $y,p,h\ge0$, $y^+=(y+\tau p)/(1+\tau h)$,
\begin{equation}\label{eq:pd-sign}
\tau hy-\tau hy^+=\tau h\left(y-\frac{y+\tau p}{1+\tau h}\right)=\frac{\tau^2h(hy-p)}{1+\tau h}.
\end{equation}
When $p=0$, $\tau hy^+=y-y^+\le y$, whereas $\tau hy>y$ if $\tau h>1$ and $y>0$. For $h>0$, this difference is negative when $p>hy$.
We test three species in a directed cycle at a fixed control, with positive spatial diffusion and an initially empty receiver. For $\tau=0.001,0.1,2,20$, the positivity, mass and resolvent contraction checks pass to roundoff, with mass errors below $1.8\times10^{-15}$. Without production, eight steps at $\tau=20$ harvest about $2.969$ from an initial mass of $3.844$. The old-state rule rewards $46.31$.
If the coefficients are $C^2$ on an open neighborhood of the computed state--control pairs, the fixed-grid state map is $C^2$ near the computed control. For a $C^2$ cost, Appendix~\ref{app:network-derivative-formulas} gives the exact tangent, adjoint and Hessian formulas. Although the continuous problem has no terminal cost, the final discrete removal contributes to the terminal adjoint.
\subsection{Relation to existing work}
Patankar positivity and conservative transfers \cite{Burchard2003,KopeczMeister2018} extend to multistep methods \cite{Izzo2025}. Accuracy and oscillations are studied in \cite{Torlo2022,IzginOffner2023}. Cacace, Oliviero and Pezzella \cite{Cacace2026} use dynamic programming for production--destruction control with a modified Patankar semi-Lagrangian Hamilton--Jacobi--Bellman approximation. Our two reward rules share a positive state solver, so removal quadrature accounts for their curvature difference. We obtain derivatives directly from the discrete problem \cite{Hager2000,Herty2013,SanzSerna2016,Ito2021}.
\raggedbottom
\section{Conclusion}
The positive state scheme gives a direct way to count the catch. In the model with harvesting alone, summing its removal flux gives the total population loss and a strictly convex reduced cost. Using the starting density adds a reward proportional to $(u^n)^2y^{n+1}$ at each step. This can favor extreme allocations and introduce a strict local minimum with a higher cost than another feasible control. The two-species example shows that the same change in local optimality can occur with positive diffusion and spatially varying admissible controls.

For $P_0$ and $P_1$ in one space dimension, the costs and their first two derivatives at the same control differ by $O(\tau)$ uniformly in the spatial mesh. A uniform positive lower bound on the old-state Hessian throughout a common neighborhood gives the same order for the distance between the corresponding stationary controls. With uniformly bounded tracking targets, a sufficiently large effort penalty and a sufficiently small time step make each objective strongly convex and give it a unique global minimizer.

On coarse time grids, the old-state controls in our tests switched more sharply and were more sensitive to the initial guess. In the refinement study at $M=20$, the distance between the two controls first grew and then fell as the time grid was refined. At $N=960$, the three reference costs differed only in the fifth decimal place on the chosen validation grid. The population balance also gives a harvest rule for weighted production--destruction networks, where the same transfer flux records the loss from one species and the gain to another.
\section*{Acknowledgment}
The author used AI tools to edit the language, organize the manuscript, and help write and review the code used in the numerical experiments. These tools also assisted with formalizing and checking proofs in \(\mathrm{L}\exists\forall\mathrm{N}\) (Lean).

\balance
\balance
\setlength{\splitmaxdepth}{0pt}

\onecolumn
\nobalance
\normalsize
\setlength{\abovedisplayskip}{10pt plus 2pt minus 2pt}
\setlength{\belowdisplayskip}{10pt plus 2pt minus 2pt}
\setlength{\abovedisplayshortskip}{10pt plus 2pt minus 2pt}
\setlength{\belowdisplayshortskip}{10pt plus 2pt minus 2pt}
\setlength{\parindent}{0pt}
\setlength{\parskip}{5pt plus 1pt minus 1pt}
\let\mapdef\compactmapdef
\appendix
\section{Continuous proofs}\label{app:continuous}
\subsection{State bounds and well-posedness}
We use $\|\cdot\|_2$ for the $L^2(\Omega)$ norm and the Euclidean product norm for vectors. To truncate the reactions, fix $R>0$ and set
\[
 \mapdef{T_R}{\R}{[0,R]}{s}{\min\{R,\max\{0,s\}\}}.
\]
The truncated reactions are
\[
f_{R,1}(v,u)=-mv_1+\alpha T_R(v_1)T_R(v_2),\qquad
f_{R,2}(v,u)=(r-qu)v_2-\beta T_R(v_1)T_R(v_2).
\]
Since $T_R$ is $1$-Lipschitz and $0\le T_R\le R$,
\begin{equation*}
|T_R(a_1)T_R(a_2)-T_R(b_1)T_R(b_2)|
\le R|a_1-b_1|+R|a_2-b_2|,
\end{equation*}
and 
\begin{equation*}
|f_R(a,u)-f_R(b,u)|\le C_R|a-b|,\qquad |f_R(a,u)|\le C_R(1+|a|).
\end{equation*}
Here $C_R$ is uniform for $0\le u\le U_{\max}$. To construct the solution, take an $L^2(\Omega)$-orthonormal Neumann eigenbasis $(e_\ell)_{\ell\ge1}$ and let $P_K$ be the projection onto its first $K$ elements. Write
\[
v_{i,K}(t)=\sum_{\ell=1}^K c_{i,\ell}(t)e_\ell,\qquad
v_{i,K}(0)=P_Kv_i^0.
\]
The Galerkin equations, for $1\le\ell\le K$, are
\[
c_{i,\ell}'(t)+\mu_i\sum_{k=1}^K c_{i,k}(t)(\nabla e_k,\nabla e_\ell)_\Omega
=(f_{R,i}(v_K,u),e_\ell)_\Omega.
\]
Their coefficients are measurable in time and globally Lipschitz in the coefficient vector. Multiplying by $c_{i,\ell}$ and summing gives
\begin{equation*}
\frac12\frac{\mathrm{d}}{\mathrm{d}t}\sum_{i=1}^2\|v_{i,K}\|_2^2
+\sum_{i=1}^2\mu_i\|\nabla v_{i,K}\|_2^2
=\sum_{i=1}^2(f_{R,i}(v_K,u),v_{i,K})_\Omega \le C_R\left(1+\sum_{i=1}^2\|v_{i,K}\|_2^2\right).
\end{equation*}
Gronwall and integration in time yield bounds independent of $K$:
\[
\sup_{0\le t\le T}\sum_{i=1}^{2}\|v_{i,K}(t)\|_2^2
+\sum_{i=1}^{2}\int_0^T\|\nabla v_{i,K}\|_2^2\,\mathrm{d}t
\le C_R.
\]
The projections $P_K$ are bounded in $H^1(\Omega)$. For $\phi\in H^1(\Omega)$, the Galerkin equation therefore gives
\[
|\langle\partial_tv_{i,K},\phi\rangle|
\le\bigl(\mu_i\|\nabla v_{i,K}\|_2+\|f_{R,i}(v_K,u)\|_2\bigr)\|\phi\|_{H^1},
\]
and hence $\|\partial_tv_{i,K}\|_{L^2(0,T;H^1{}^*)}\le C_R$. Compactness gives a subsequence with
\[
v_K\xrightarrow[K\to\infty]{L^2(Q)^2}v,
\qquad
\|f_R(v_K,u)-f_R(v,u)\|_{L^2(Q)^2}
\le C_R\|v_K-v\|_{L^2(Q)^2}\xrightarrow[K\to\infty]{}0.
\]
Passing to the limit in the time-integrated Galerkin equations gives the weak Neumann problem and the initial trace $v_i(0)=v_i^0$ \cite{Amann1984,Evans2010}.
Set $s^-\coloneqq\max\{-s,0\}$ and test component $i$ with $-v_i^-\in H^1(\Omega)$. The interaction terms vanish because $T_R(v_i)=0$ where $v_i<0$. For almost every $t\in(0,T)$,
\begin{equation*}
\begin{aligned}
&\frac12\frac{\mathrm{d}}{\mathrm{d}t}\sum_{i=1}^2\|v_i^-\|_2^2 +\sum_{i=1}^2\mu_i\|\nabla v_i^-\|_2^2 \le C\sum_{i=1}^2\|v_i^-\|_2^2.
\end{aligned}
\end{equation*}
The two reaction terms on the negative sets are
\[
\int_\Omega f_{R,1}(v,u)(-v_1^-)=-m\|v_1^-\|_2^2,\qquad
\int_\Omega f_{R,2}(v,u)(-v_2^-)
=\int_\Omega(r-qu)(v_2^-)^2\le r\|v_2^-\|_2^2.
\]
Thus
\[
\sum_{i=1}^{2}\|v_i^-(t)\|_2^2
\le e^{2rt}\sum_{i=1}^{2}\|(v_i^0)^-\|_2^2=0.
\]
Since $0\le T_R(v_i)\le v_i$, in $\Omega\times(0,T)$ we have
\begin{equation*}
\partial_tv_2-\mu_2\Delta v_2\le rv_2.
\end{equation*}
For the comparison, set $b_2(t)=V_2^0e^{rt}$ and test the equation for $v_2-b_2$ with $(v_2-b_2)^+$. Since $b_2'=rb_2$ and $(v_2^0-b_2(0))^+=0$,
\[
\frac12\frac{\mathrm{d}}{\mathrm{d}t}\|(v_2-b_2)^+\|_2^2
+\mu_2\|\nabla(v_2-b_2)^+\|_2^2
\le r\|(v_2-b_2)^+\|_2^2.
\]
Gronwall gives, for almost every $(x,t)\in\Omega\times(0,T)$,
\begin{equation*}
v_2(x,t)\le V_2^0e^{rt}.
\end{equation*}
In $\Omega\times(0,T)$, we then have
\begin{equation*}
\partial_tv_1-\mu_1\Delta v_1\le(-m+\alpha V_2^0e^{rt})v_1.
\end{equation*}
Compare $v_1$ with the solution of the ODE on $(0,T)$, with initial value $V_1^0$:
\begin{equation*}
y'(t)=(-m+\alpha V_2^0e^{rt})y(t),\qquad y(0)=V_1^0.
\end{equation*}
Its solution is
\[
y(t)=V_1^0\exp\!\left(\int_0^t[-m+\alpha V_2^0e^{rs}]\,\mathrm{d}s\right) =V_1^0\exp\!\left(-mt+\frac{\alpha V_2^0}{r}(e^{rt}-1)\right).
\]
Testing the inequality for $v_1-y$ with $(v_1-y)^+$ gives
\[
\frac12\frac{\mathrm{d}}{\mathrm{d}t}\|(v_1-y)^+\|_2^2
+\mu_1\|\nabla(v_1-y)^+\|_2^2
\le(-m+\alpha V_2^0e^{rt})\|(v_1-y)^+\|_2^2.
\]
The initial positive part is zero, so $v_1\le y$. Once $R$ exceeds both bounds on $[0,T]$, the truncation disappears: the solution satisfies \eqref{eq:pde}. To prove uniqueness, take two solutions $v,\widetilde v\in\Y^2$ and set $e\coloneqq v-\widetilde v\in\mathcal W^2$. Then
\begin{equation}\label{eq:pfa_product_difference}
v_1v_2-\widetilde v_1\widetilde v_2=e_1v_2+\widetilde v_1e_2.
\end{equation}
The difference equations are
\[
\left\{\begin{aligned}
\partial_te_1-\mu_1\Delta e_1&=(-m+\alpha v_2)e_1+\alpha\widetilde v_1e_2,\\
\partial_te_2-\mu_2\Delta e_2&=-\beta v_2e_1+(r-qu-\beta\widetilde v_1)e_2.
\end{aligned}\right.
\]
All coefficients are bounded. In particular,
\[
\left|\int_\Omega(\alpha\widetilde v_1-\beta v_2)e_1e_2\right|
\le\frac12\|\alpha\widetilde v_1-\beta v_2\|_\infty
(\|e_1\|_2^2+\|e_2\|_2^2),
\]
so $|f(v,u)-f(\widetilde v,u)|\le C|e|$. For almost every $t\in(0,T)$, testing the difference equation with $e$ yields
\begin{equation*}
\begin{aligned}
&\frac12\frac{\mathrm{d}}{\mathrm{d}t}\|e\|_2^2 +\sum_{i=1}^2\mu_i\|\nabla e_i\|_2^2 \le C\|e\|_2^2.
\end{aligned}
\end{equation*}
Since $e(0)=0$, Gronwall gives $e=0$.
Testing with $\phi=1\in H^1(\Omega)$ eliminates diffusion. For almost every $t\in(0,T)$,
\begin{equation*}
\mathcal M'
\byeq{eq:weakstate}{=}\beta\int_\Omega(-mv_1+\alpha v_1v_2)
 +\alpha\int_\Omega(rv_2-\beta v_1v_2-quv_2) =\underbrace{-\beta m\int_\Omega v_1}_{\le0}
+\underbrace{\alpha r\int_\Omega v_2}_{\le r\mathcal M}
+\underbrace{(-\alpha) q\int_\Omega uv_2}_{\le0}
\le r\mathcal M.
\end{equation*}
Equivalently,
\[
\frac{\mathrm{d}}{\mathrm{d}t}\bigl(e^{-rt}\mathcal M(t)\bigr)
=-e^{-rt}\left[\beta(m+r)\int_\Omega v_1+\alpha q\int_\Omega uv_2\right]\le0.
\]
Integration from $0$ to $t$ gives $\mathcal M(t)\le e^{rt}\mathcal M(0)$. For $\mu_1=\mu_2=\mu$, set $w\coloneqq\beta v_1+\alpha v_2\in\Y$. In $\Omega\times(0,T)$,
\begin{equation*}
\partial_tw-\mu\Delta w
=-\beta mv_1+\alpha rv_2-\alpha quv_2\le rw,
\end{equation*}
and comparison with $e^{rt}\|w(0)\|_\infty$ proves \eqref{eq:equaldiff}.
\subsection{Existence of a continuous minimizer}
The bound $j(u)\ge-\eta qU_{\max}K_2|Q|$, with $K_2\coloneqq V_2^0e^{rT}$, gives a finite infimum. Take a minimizing sequence $(u_k)\subset\Uad$. Because $Q$ has finite measure and $L^1(Q)$ is separable, we can extract a subsequence such that, in $L^\infty(Q)$, $u_k\mathrel{\underset{k\to\infty}{\rightharpoonup^{*}}}u^*$. For $0\le\phi\in L^1(Q)$,
\begin{equation*}
0\le\int_Q u^*\phi=\lim_k\int_Q u_k\phi\le U_{\max}\int_Q\phi,
\end{equation*}
so $u^*\in\Uad$. Set $v_k\coloneqq S(u_k)\in\Y^2$, $K_1\coloneqq V_1^0\max_{0\le t\le T}\exp\!\left(-mt+\frac{\alpha V_2^0}{r}(e^{rt}-1)\right)$. Theorem~\ref{th:continuous} gives $\|v_{i,k}\|_{L^\infty(Q)}\le K_i$ uniformly in $k$. For almost every $t\in(0,T)$, testing by $v_{i,k}$ gives
\begin{equation*}
\frac12\frac{\mathrm{d}}{\mathrm{d}t}\|v_{i,k}\|_{L^2(\Omega)}^2
+\mu_i\|\nabla v_{i,k}\|_{L^2(\Omega)}^2
=(f_i(v_k,u_k),v_{i,k})_\Omega.
\end{equation*}
The right side is uniformly bounded, so integration in time, together with
\begin{equation*}
\|\partial_tv_{i,k}\|_{H^1{}^*}
\le\mu_i\|\nabla v_{i,k}\|_{L^2(\Omega)}+\|f_i(v_k,u_k)\|_{L^2(\Omega)}
\end{equation*}
gives
\begin{equation*}
\|v_{i,k}\|_{L^2(0,T;H^1)}
+\|\partial_tv_{i,k}\|_{L^2(0,T;H^1{}^*)}\le C.
\end{equation*}
The compact embedding $H^1(\Omega)\hookrightarrow L^2(\Omega)$ and time compactness \cite{Simon1987} give a subsequence converging almost everywhere in $Q$, with
\begin{equation*}
v_{i,k}\xrightharpoonup[k\to\infty]{\mathcal W}v_i^*,
\qquad
v_{i,k}\xrightarrow[k\to\infty]{L^2(Q)}v_i^*.
\end{equation*}
Thus $0\le v_i^*\le K_i$ and $v^*\in\Y^2$. Then
\begin{align*}
\|v_{1,k}v_{2,k}-v_1^*v_2^*\|_2
\byeq{eq:pfa_product_difference}{\le}K_1\|v_{2,k}-v_2^*\|_2+K_2\|v_{1,k}-v_1^*\|_2\xrightarrow[k\to\infty]{}0.
\end{align*}
For $\phi\in L^2(Q)$,
\begin{align*}
\int_Q(u_kv_{2,k}-u^*v_2^*)\phi
&=\int_Q u_k(v_{2,k}-v_2^*)\phi+\int_Q(u_k-u^*)v_2^*\phi.
\end{align*}
Then
\begin{equation*}
\left|\int_Q u_k(v_{2,k}-v_2^*)\phi\right|
\le U_{\max}\|v_{2,k}-v_2^*\|_2\|\phi\|_2\xrightarrow[k\to\infty]{}0.
\end{equation*}
The second term tends to zero by weak-$*$ convergence, since $v_2^*\phi\in L^1(Q)$. Hence $u_kv_{2,k}\xrightharpoonup[k\to\infty]{L^2(Q)}u^*v_2^*$. Test the weak equations with $\phi(x)\chi(t)$, where $\phi\in H^1(\Omega)$, $\chi\in C^1([0,T])$ and $\chi(T)=0$. We pass to the limit in every term, including the fixed initial trace. Uniqueness identifies $v^*=S(u^*)$.
For each tracking term, strong convergence gives
\begin{equation*}
\bigl|\|v_{i,k}-v_i^d\|_2^2-\|v_i^*-v_i^d\|_2^2\bigr|
=\left|\int_Q(v_{i,k}-v_i^*)(v_{i,k}+v_i^*-2v_i^d)\right| \le\|v_{i,k}-v_i^*\|_2
\bigl(\|v_{i,k}-v_i^d\|_2+\|v_i^*-v_i^d\|_2\bigr)\xrightarrow[k\to\infty]{}0.
\end{equation*}
For the harvest term,
\begin{align*}
\left|\int_Q u_kv_{2,k}-\int_Q u^*v_2^*\right|
&\le U_{\max}|Q|^{1/2}\|v_{2,k}-v_2^*\|_2
+\left|\int_Q(u_k-u^*)v_2^*\right|
\xrightarrow[k\to\infty]{}0.
\end{align*}
For the effort cost, weak convergence $u_k\xrightharpoonup[k\to\infty]{L^2(Q)}u^*$ gives
\begin{equation*}
\|u^*\|_2^2\le\liminf_k\|u_k\|_2^2
\quad\Longrightarrow\quad
j(u^*)\le\liminf_k j(u_k).
\end{equation*}
Thus $u^*$ attains the minimum.
\subsection{First and second derivative formulas}\label{app:continuous-derivative-formulas}
For $u,\xi\in L^\infty(Q)$, the tangent $z^\xi=S'(u)\xi$ satisfies the following system. The evolution equations hold in $Q$, the Neumann conditions on $\partial\Omega\times(0,T)$, and the initial conditions in $\Omega$ for $i=1,2$:
\begin{equation}\label{eq:tangentc}
\left\{\begin{aligned}
\partial_tz_1-\mu_1\Delta z_1&=(-m+\alpha v_2)z_1+\alpha v_1z_2,\\
\partial_tz_2-\mu_2\Delta z_2&=-\beta v_2z_1+(r-\beta v_1-qu)z_2-qv_2\xi,\\
\partial_nz_i&=0,\qquad z_i(0)=0.
\end{aligned}\right.
\end{equation}
The \emph{adjoint} $p\coloneqq (p_1,p_2)\in\mathcal W^2$ satisfies the following system. The evolution equations hold in $Q$, the Neumann conditions on $\partial\Omega\times(0,T)$, and the terminal conditions in $\Omega$ for $i=1,2$:
\begin{equation}\label{eq:adjc}
\left\{\begin{aligned}
-\partial_tp_1-\mu_1\Delta p_1&=\gamma_1(v_1-v_1^d)+(-m+\alpha v_2)p_1-\beta v_2p_2,\\[2pt]
-\partial_tp_2-\mu_2\Delta p_2&=\gamma_2(v_2-v_2^d)-\eta qu+\alpha v_1p_1+(r-\beta v_1-qu)p_2,\\[2pt]
\partial_np_i&=0,\qquad p_i(T)=0.
\end{aligned}\right.
\end{equation}
For $z^\xi\coloneqq S'(u)\xi$, $z^\zeta\coloneqq S'(u)\zeta$, with $\xi,\zeta\in L^\infty(Q)$,
\begin{equation}\label{eq:hessc}
\mapdef{j''(u)}{L^\infty(Q)\times L^\infty(Q)}{\R}{(\xi,\zeta)}{
 \lambda\int_Q\xi\zeta+\sum_{i=1}^2\gamma_i\int_Qz_i^\xi z_i^\zeta +\int_Q(\alpha p_1-\beta p_2) \cdot(z_1^\xi z_2^\zeta+z_1^\zeta z_2^\xi) -q\int_Q(\eta+p_2)(\xi z_2^\zeta+\zeta z_2^\xi).}
\end{equation}
\subsection{Differentiability and second variation}
Fix $u\in L^\infty(Q)$ and $R>\|u\|_\infty$, and set $v=S(u)$. Consider increments $\zeta\in L^\infty(Q)$ such that $\|u+\zeta\|_\infty\le R$. Nonnegativity holds for signed controls too. Replacing $r$ by $r+qR$ gives the prey bound. Thus both the states and the reaction Jacobian
\begin{equation*}
\mathsf A(v,u)\coloneqq
\begin{pmatrix}
-m+\alpha v_2&\alpha v_1\\
-\beta v_2&r-\beta v_1-qu
\end{pmatrix}\in L^\infty(Q;\R^{2\times2})
\end{equation*}
are uniformly bounded. Define
\begin{equation*}
\mapdef{\mathscr L_v}{\Y^2}{L^2(0,T;H^1(\Omega)^*)^2}{y}{
 \partial_ty-\diag(\mu_1,\mu_2)\Delta y-\mathsf A(v,u)y}
\end{equation*}
with weak Neumann conditions. Write $\|y(t)\|_p\coloneqq\|y(t)\|_{L^p(\Omega)^2}$, and likewise for $F$. If $\mathscr L_vy=F\in L^\infty(Q)^2$ and $y(0)=0$, variation of constants for the Neumann heat semigroups gives
\begin{equation*}
\|y(t)\|_\infty
\le C\int_0^t\bigl(\|y(s)\|_\infty+\|F(s)\|_\infty\bigr)\,\mathrm{d}s.
\end{equation*}
Let $C_A$ bound $\max_i\sum_{j=1}^{2}|\mathsf A_{ij}|$ on this control ball. The componentwise Neumann heat semigroup is contractive in the maximum norm. Its variation-of-constants formula and Gronwall give
\[
\max_i\|y_i(t)\|_{L^\infty(\Omega)}
\le\int_0^t e^{C_A(t-s)}\max_i\|F_i(s)\|_{L^\infty(\Omega)}\,\mathrm{d}s
\le T e^{C_AT}\max_i\|F_i\|_{L^\infty(Q)}.
\]
The inequalities $\max_i|a_i|\le|a|_{\R^2}\le\sqrt2\max_i|a_i|$ give the corresponding bound in the Euclidean product norm. For almost every $t\in(0,T)$, testing by $y$ gives
\begin{equation*}
\frac12\frac{\mathrm{d}}{\mathrm{d}t}\|y\|_2^2+\sum_{i=1}^2\mu_i\|\nabla y_i\|_2^2
=\underbrace{(\mathsf Ay,y)_\Omega}_{\le C\|y\|_2^2}
+\underbrace{(F,y)_\Omega}_{\le\frac12\|F\|_2^2+\frac12\|y\|_2^2} \le C\|y\|_2^2+\frac12\|F\|_2^2.
\end{equation*}
Integration and Gronwall give
\[
\sup_{0\le t\le T}\|y(t)\|_2^2
+\sum_{i=1}^{2}\mu_i\|\nabla y_i\|_{L^2(Q)}^2
\le C_R\|F\|_{L^2(Q)^2}^2.
\]
For $\phi\in H^1(\Omega)^2$, the equation also gives
\[
|\langle\partial_ty,\phi\rangle|
\le C_R\bigl(\|\nabla y\|_2+\|y\|_2+\|F\|_2\bigr)\|\phi\|_{H^1(\Omega)^2}.
\]
Taking the dual norm and then integrating in time bounds $\partial_ty$ in $L^2(0,T;H^1{}^*)^2$. Together, these estimates give
\begin{equation}\label{eq:pfa_linear_estimate}
\|y\|_{\Y^2}\le C_R\|F\|_{L^\infty(Q)^2}.
\end{equation}
The constant $C_R$ depends only on the data, $T$ and $R$.
For $\xi\in L^\infty(Q)$, the tangent satisfies the following evolution equation in $Q$, initial condition in $\Omega$, and norm bound:
\begin{equation}\label{eq:pfa_tangent_source}
\mathscr L_v z^\xi=\binom{0}{-qv_2\xi},\qquad
z^\xi(0)=0,\qquad \|z^\xi\|_{\Y^2}\le C_R\|\xi\|_\infty.
\end{equation}
For $\zeta\in L^\infty(Q)$, set $\delta v\coloneqq S(u+\zeta)-S(u)$. The state difference solves an equation with bounded coefficients and source $(0,-q\zeta v_2(u+\zeta))^T$, so $\|\delta v\|_{\Y^2}\le C_R\|\zeta\|_\infty$. Subtracting the tangent equation gives the following evolution equation in $Q$ and initial condition in $\Omega$:
\begin{equation*}
\mathscr L_v(\delta v-z^\zeta)=\binom{\alpha\delta v_1\delta v_2}{-\beta\delta v_1\delta v_2-q\zeta\delta v_2},\qquad
(\delta v-z^\zeta)(0)=0.
\end{equation*}
This source is $O(\|\zeta\|_\infty^2)$ in $L^\infty(Q)^2$, so
\begin{equation*}
\|S(u+\zeta)-S(u)-z^\zeta\|_{\Y^2}
\byeq{eq:pfa_linear_estimate}{\le}C_R\|\zeta\|_\infty^2.
\end{equation*}
This proves $S'(u)\xi=z^\xi$.
For $\xi,\zeta\in L^\infty(Q)$, define $B_{\xi\zeta},D_{\xi\zeta}\in L^\infty(Q)$ by
\begin{equation}\label{eq:pfa_second_sources}
B_{\xi\zeta}\coloneqq z_1^\xi z_2^\zeta+z_1^\zeta z_2^\xi,
\qquad D_{\xi\zeta}\coloneqq\xi z_2^\zeta+\zeta z_2^\xi.
\end{equation}
Then $w\coloneqq S''(u)[\xi,\zeta]\in\Y^2$ satisfies the following system. The evolution equations hold in $Q$, the Neumann conditions on $\partial\Omega\times(0,T)$, and the initial conditions in $\Omega$ for $i=1,2$:
\begin{equation}\label{eq:secondstate_explicit}
\left\{\begin{aligned}
\partial_tw_1-\mu_1\Delta w_1&=(-m+\alpha v_2)w_1+\alpha v_1w_2+\alpha B_{\xi\zeta},\\
\partial_tw_2-\mu_2\Delta w_2&=-\beta v_2w_1+(r-\beta v_1-qu)w_2-\beta B_{\xi\zeta}-qD_{\xi\zeta},\\[2pt]
\partial_nw_i&=0,\qquad w_i(0)=0,
\end{aligned}\right.
\end{equation}
with $\|w\|_{\Y^2}\le C_R\|\xi\|_\infty\|\zeta\|_\infty$. To identify this as the second derivative, set $\widetilde z^\xi\coloneqq S'(u+\zeta)\xi$, $\delta z^\xi\coloneqq\widetilde z^\xi-z^\xi$. Since
\begin{equation*}
\delta\mathsf A\coloneqq
\begin{pmatrix}
\alpha\delta v_2&\alpha\delta v_1\\
-\beta\delta v_2&-\beta\delta v_1-q\zeta
\end{pmatrix},
\end{equation*}
subtraction gives, in $\Omega\times(0,T)$,
\begin{equation*}
\mathscr L_v\delta z^\xi
=\delta\mathsf A\,\widetilde z^\xi+\binom{0}{-q\delta v_2\xi},
\end{equation*}
thus
\begin{equation*}
\|\delta z^\xi\|_{\Y^2}\le C_R\|\xi\|_\infty\|\zeta\|_\infty.
\end{equation*}
Set $r^\zeta=\delta v-z^\zeta$. Subtracting \eqref{eq:secondstate_explicit} gives the remainder equation
\[
\mathscr L_v(\delta z^\xi-w)
=\delta\mathsf A\,\delta z^\xi
+\begin{pmatrix}
\alpha(r_2^\zeta z_1^\xi+r_1^\zeta z_2^\xi)\\
-\beta(r_2^\zeta z_1^\xi+r_1^\zeta z_2^\xi)-q\xi r_2^\zeta
\end{pmatrix},
\qquad(\delta z^\xi-w)(0)=0.
\]
The two parts of the source satisfy
\begin{equation*}
\|\delta\mathsf A\,\delta z^\xi\|_\infty
\le C_R\|\zeta\|_\infty\|\delta z^\xi\|_\infty
\le C_R\|\xi\|_\infty\|\zeta\|_\infty^2,\qquad \|r^\zeta\|_\infty(\|z^\xi\|_\infty+\|\xi\|_\infty)
\le C_R\|\xi\|_\infty\|\zeta\|_\infty^2.
\end{equation*}
Consequently,
\begin{equation*}
\|\widetilde z^\xi-z^\xi-w\|_{\Y^2}
\byeq{eq:pfa_linear_estimate}{\le}C_R\|\xi\|_\infty\|\zeta\|_\infty^2.
\end{equation*}
For two base controls $u,\widetilde u$ in this ball, subtracting their tangent equations and applying \eqref{eq:pfa_linear_estimate} gives
\[
\|S'(\widetilde u)\xi-S'(u)\xi\|_{\Y^2}
\le C_R\|\widetilde u-u\|_\infty\|\xi\|_\infty.
\]
In the second-variation source, differences of products are expanded as
\[
\widetilde z_1^\xi\widetilde z_2^\zeta-z_1^\xi z_2^\zeta
=(\widetilde z_1^\xi-z_1^\xi)\widetilde z_2^\zeta
+z_1^\xi(\widetilde z_2^\zeta-z_2^\zeta).
\]
Each term is bounded by $C_R\|\widetilde u-u\|_\infty\|\xi\|_\infty\|\zeta\|_\infty$. The change in the linear operator contributes
$(\mathsf A(S(\widetilde u),\widetilde u)-\mathsf A(S(u),u))S''(\widetilde u)[\xi,\zeta]$, with the same bound. Hence
\[
\|S''(\widetilde u)[\xi,\zeta]-S''(u)[\xi,\zeta]\|_{\Y^2}
\le C_R\|\widetilde u-u\|_\infty\|\xi\|_\infty\|\zeta\|_\infty.
\]
Taking suprema over unit directions proves continuity of both derivatives in operator norm. Thus $S:L^\infty(Q)\to\Y^2$ is $C^2$ on bounded control sets.
Set
\begin{equation*}
\ell_v\coloneqq\binom{\gamma_1(v_1-v_1^d)}{\gamma_2(v_2-v_2^d)-\eta qu}
\in L^2(Q)^2.
\end{equation*}
The backward linear problem \eqref{eq:adjc} has a weak Neumann solution $p\in\mathcal W^2$. For $y=z^\xi$ or $w$, integration by parts gives
\begin{equation}\label{eq:pfa_duality}
\int_Qp\cdot\mathscr L_vy =\int_0^T\bigg[\langle\partial_ty,p\rangle +\sum_{i=1}^2\mu_i(\nabla y_i,\nabla p_i)_\Omega-(\mathsf Ay,p)_\Omega\bigg]\mathrm{d}t =\left[(p,y)_\Omega\right]_0^T+\int_Q\ell_v\cdot y.
\end{equation}
Here $\langle\cdot,\cdot\rangle$ is the componentwise dual pairing between $(H^1(\Omega)^*)^2$ and $H^1(\Omega)^2$. The trace term vanishes because $p(T)=y(0)=0$, and the chain rule gives
\[
j'(u)\xi =\int_Q\ell_v\cdot z^\xi+\int_Q(\lambda u-\eta qv_2)\xi \byeq{eq:pfa_duality}{=}\int_Qp\cdot\mathscr L_vz^\xi +\int_Q(\lambda u-\eta qv_2)\xi \byeq{eq:pfa_tangent_source}{=}\int_Q[\lambda u-qv_2(\eta+p_2)]\xi,
\]
proving \eqref{eq:gradc}. At second order,
\begin{equation*}
j''(u)[\xi,\zeta] ={}\lambda\int_Q\xi\zeta+\sum_{i=1}^2\gamma_i\int_Qz_i^\xi z_i^\zeta -\eta q\int_QD_{\xi\zeta}+\int_Q\ell_v\cdot w,
\end{equation*}
and
\[
\int_Q\ell_v\cdot w \byeq{eq:pfa_duality}{=}\int_Qp\cdot\mathscr L_vw \byeq{eq:secondstate_explicit}{=}\int_Q(\alpha p_1-\beta p_2)B_{\xi\zeta} -q\int_Qp_2D_{\xi\zeta}.
\]
Substitution gives \eqref{eq:hessc}. The bounds $\xi,\zeta\in L^\infty(Q)$, $z^\xi,z^\zeta\in L^\infty(Q)^2$ and $p\in L^2(Q)^2$ ensure that every term is integrable.
\subsection{Variational inequality and pointwise projection}
Let $u^*$ be a local minimizer and take $u\in\Uad$. The segment $u^*+t(u-u^*)$ is admissible for $0\le t\le1$. Thus
\[
0\le\lim_{t\to0^+}\frac{j(u^*+t(u-u^*))-j(u^*)}{t}
=\int_Q g(u^*)(u-u^*).
\]
Set $a=qS_2(u^*)(\eta+p_2(u^*))/\lambda\in L^2(Q)$ and $w=\proj(a)\in\Uad$. The scalar projection satisfies, for every $b\in[0,U_{\max}]$,
\[
(w-a)(b-w)\ge0.
\]
Apply this pointwise with $b=u^*$ and use $g(u^*)=\lambda(u^*-a)$. Testing the variational inequality with $w$ gives
\begin{equation*}
0\le\int_Q(u^*-a)(w-u^*) =-\|u^*-w\|_{L^2(Q)}^2
+\underbrace{\int_Q(w-a)(w-u^*)}_{\le0}
\le-\|u^*-w\|_{L^2(Q)}^2.
\end{equation*}
Therefore $u^*=w$, which is the pointwise projection formula. Conversely, the scalar projection inequality integrates to the variational inequality. The pointwise conditions are
\[
\begin{cases}
g(u^*)\ge0 & \text{where }u^*=0,\\
g(u^*)=0 & \text{where }0<u^*<U_{\max},\\
g(u^*)\le0 & \text{where }u^*=U_{\max}.
\end{cases}
\]
\section{Discrete bounds}\label{app:discrete}
\subsection{Positivity and population balance}
\begin{proof}
For $D\coloneqq\diag(d_0,\ldots,d_M)\in\R^{(M+1)\times(M+1)}$, $d_j\ge0$, $\mu>0$ and $\tau>0$, set
\begin{equation*}
K\coloneqq I-\tau\mu L_h+\tau D.
\end{equation*}
Its off-diagonal entries are nonpositive, and $L_h\one=0$ gives
\begin{equation*}
K_{jj}-\sum_{\substack{0\le k\le M\\ k\ne j}}|K_{jk}|=1+\tau d_j>0.
\end{equation*}
If $y\in\R^{M+1}$, $Ky\ge0$, $y_j=\min_{0\le k\le M}y_k<0$, then
\begin{equation*}
(L_hy)_j=\sum_{\substack{0\le k\le M\\ k\ne j}}(L_h)_{jk}(y_k-y_j)\ge0,
\end{equation*}
so
\begin{equation*}
(Ky)_j=(1+\tau d_j)y_j-\tau\mu(L_hy)_j<0,
\end{equation*}
a contradiction. Applying the argument to $y$ and $-y$ when $Ky=0$ proves invertibility. Taking $y=K^{-1}b$ with $b\ge0$ then gives $K^{-1}\ge0$. We obtain
\begin{equation}\label{eq:pfa_inverse_bound}
K^{-1}\ge0,\qquad K\one\ge\one,\qquad
\|K^{-1}\|_\infty=\|K^{-1}\one\|_\infty\le1.
\end{equation}
The same inverse is contractive in the weighted $L^2$ norm. Indeed, if $Ky=b$, weighted summation by parts gives
\begin{equation*}
\|y\|_h^2+\underbrace{\tau\mu|y|_{1,h}^2}_{\ge0}
+\underbrace{\tau\sum_{j=0}^M w_jd_jy_j^2}_{\ge0}
=\underbrace{\langle b,y\rangle_h}_{\le\|b\|_h\|y\|_h}\quad\Longrightarrow\quad \|K^{-1}b\|_h\le\|b\|_h.
\end{equation*}
Also $B\one=(1+\tau m)\one$ gives $\|B^{-1}\|_\infty=(1+\tau m)^{-1}$, so solving either linear system in \eqref{eq:positive} preserves nonnegativity. Weighted symmetry gives $\ip{L_hy}{\one}=\ip{y}{L_h\one}=0$. Testing with $\one$ yields
\begin{equation*}
\ip{v_1^{n+1}-v_1^n}\one
=\tau\ip{-mv_1^{n+1}+\alpha v_1^n\circ v_2^{n+1}}\one,\qquad
\ip{v_2^{n+1}-v_2^n}\one
=\tau\ip{rv_2^n-(\beta v_1^n+qu^n)\circ v_2^{n+1}}\one.
\end{equation*}
Multiplication by $\beta$ and $\alpha$ gives the transfer cancellation
\[
\tau\alpha\beta\langle v_1^n\circ v_2^{n+1},\one\rangle_h
-\tau\alpha\beta\langle v_1^n\circ v_2^{n+1},\one\rangle_h=0.
\]
Consequently,
\begin{equation*}
\mathcal M_h^{n+1}-\mathcal M_h^n
=-\tau\beta m\langle v_1^{n+1},\one\rangle_h
+\tau\alpha r\langle v_2^n,\one\rangle_h
-\tau\alpha q\langle u^n,v_2^{n+1}\rangle_h,
\end{equation*}
and
\begin{equation*}
\mathcal M_h^N
+\tau\beta m\sum_{n=0}^{N-1}\langle v_1^{n+1},\one\rangle_h
+\tau\alpha q\sum_{n=0}^{N-1}\langle u^n,v_2^{n+1}\rangle_h
=\mathcal M_h^0+\tau\alpha r\sum_{n=0}^{N-1}\langle v_2^n,\one\rangle_h.
\end{equation*}
Dropping the nonpositive terms in the one-step identity gives
\begin{equation*}
\mathcal M_h^{n+1}\byeq{eq:balancep}{\le}\mathcal M_h^n+\tau\alpha r\ip{v_2^n}\one
\le(1+\tau r)\mathcal M_h^n.
\end{equation*}
\end{proof}
\subsection{Uniform bounds}
For $y,z\in\R^{M+1}$, weighted summation by parts and the definition of $\Ih$ give
\begin{equation}\label{eq:bilinear_spatial}
a_h(y,z)=\frac1h\sum_{j=0}^{M-1}(y_{j+1}-y_j)(z_{j+1}-z_j) =\int_0^1\partial_x\Ih y\,\partial_x\Ih z\,\mathrm{d}x =-\ip{L_hy}z.
\end{equation}
To see the role of the endpoint weights, expand the left side:
\begin{equation*}
-\langle L_hy,z\rangle_h
=\frac1h\left[(y_0-y_1)z_0
+\sum_{j=1}^{M-1}(2y_j-y_{j-1}-y_{j+1})z_j
+(y_M-y_{M-1})z_M\right] =\frac1h\sum_{j=0}^{M-1}(y_{j+1}-y_j)(z_{j+1}-z_j).
\end{equation*}
On $[x_j,x_{j+1}]$, the reconstructed derivative is $(y_{j+1}-y_j)/h$. Thus
\[
\int_{x_j}^{x_{j+1}}(\Ih y)'(\Ih z)'\,\mathrm{d}x
=\frac{(y_{j+1}-y_j)(z_{j+1}-z_j)}h.
\]
Summing over the intervals proves \eqref{eq:bilinear_spatial}. Its symmetry and nonnegativity give
\[
W_hL_h=L_h^TW_h,\qquad
-y^TW_hL_hy=|y|_{1,h}^2\ge0.
\]
For the positive scheme, define comparison sequences by
\begin{equation}\label{eq:pfa_comparison_sequences}
P_2^{n+1}\coloneqq(1+\tau r)P_2^n,\qquad
P_1^{n+1}\coloneqq\frac{1+\tau\alpha P_2^{n+1}}{1+\tau m}P_1^n,
\end{equation}
with $P_i^0\coloneqq\|v_i^0\|_\infty$. If $0\le v_i^n\le P_i^n\one$, the inverse bounds give
\begin{equation*}
v_2^{n+1}
\le(1+\tau r)P_2^n(C^n)^{-1}\one
\overset{\mathclap{\text{\eqref{eq:pfa_inverse_bound}, \eqref{eq:pfa_comparison_sequences}}}}{\le}P_2^{n+1}\one,
\qquad
v_1^{n+1}
\le P_1^n(1+\tau\alpha P_2^{n+1})B^{-1}\one
\byeq{eq:pfa_comparison_sequences}{=}P_1^{n+1}\one.
\end{equation*}
For $s\ge0$, we have $1+s\le e^s$. Induction therefore gives
\begin{equation*}
P_2^n\byeq{eq:pfa_comparison_sequences}{=}P_2^0(1+\tau r)^n
\le P_2^0e^{rt_n}\le K_2,
\end{equation*}
and
\[
P_1^n\byeq{eq:pfa_comparison_sequences}{\le}P_1^0\prod_{\ell=0}^{n-1}(1+\tau\alpha P_2^{\ell+1}) \le P_1^0(1+\tau\alpha K_2)^n \le P_1^0e^{\alpha K_2T}=K_1.
\]
For IMEX, assume the bounds hold at time $n$ and $\tau\le\tau_0$. The old-state factors satisfy
\begin{equation*}
1+\tau(-m+\alpha v_{2,j}^n)\ge1-\tau m\ge0,\qquad
1+\tau(r-\beta v_{1,j}^n-qu_j^n)\ge1-\tau(\beta K_1+qU_{\max})\ge0.
\end{equation*}
Since $\|A_i^{-1}\|_\infty\le1$,
\begin{equation*}
0\le v_2^{n+1}\le(1+\tau r)\|v_2^n\|_\infty\one,\qquad
0\le v_1^{n+1}\le(1+\tau\alpha K_2)\|v_1^n\|_\infty\one.
\end{equation*}
This closes the induction. Under \eqref{eq:dataV}, the embedding $H^1(0,1)\hookrightarrow C([0,1])$ gives $P_i^0\le\|v_i^0\|_{L^\infty(0,1)}$, so neither $K_i$ nor $\tau_0$ depends on $h$.
For equal diffusivities, let $A\coloneqq I-\tau\mu L_h$ and $w^n\coloneqq\beta v_1^n+\alpha v_2^n\in\R^{M+1}$. The positive equations give
\[
Aw^{n+1} =w^n-\tau\beta mv_1^{n+1}+\tau\alpha rv_2^n -\tau\alpha qu^n\circ v_2^{n+1} \le(1+\tau r)w^n.
\]
Using $A^{-1}\ge0$ and $A^{-1}\one=\one$,
\[
w^{n+1}\byeq{eq:pfa_inverse_bound}{\le}(1+\tau r)\|w^n\|_\infty\one \le(1+\tau r)^{n+1}\|w^0\|_\infty\one.
\]
For nonnegative IMEX trajectories, cancellation instead gives
\[
Aw^{n+1}=w^n-\tau\beta mv_1^n+\tau\alpha rv_2^n -\tau\alpha qu^n\circ v_2^n \le(1+\tau r)w^n.
\]
The same induction proves \eqref{eq:equaldiscrete}.
\section{Discrete Lagrangians and derivative identities}\label{app:derivatives}
\subsection{IMEX}\label{app:imex-derivatives}
Set $p_i^N=0$. For $n=N-1,\ldots,1$,
\begin{equation}\label{eq:adjI}
\left\{\begin{aligned}
A_1p_1^n={}&p_1^{n+1}+\tau\big[\gamma_1(v_1^n-v_1^{d,n})
+(-m+\alpha v_2^n)\circ p_1^{n+1}
-\beta v_2^n\circ p_2^{n+1}\big],\\
A_2p_2^n={}&p_2^{n+1}+\tau\big[\gamma_2(v_2^n-v_2^{d,n})-\eta qu^n
+\alpha v_1^n\circ p_1^{n+1}+(r-\beta v_1^n-qu^n)\circ p_2^{n+1}\big].
\end{aligned}\right.
\end{equation}
The gradient is
\begin{equation}\label{eq:gradI}
g^{\mathrm I,n}=\lambda u^n-\eta qv_2^n-qv_2^n\circ p_2^{n+1}.
\end{equation}
For $d\in\R^{N(M+1)}$, $z\coloneqq DS^{\rm I}_{h,\tau}(u)d$ satisfies $z_i^0=0$ and
\begin{equation}\label{eq:tangentI}
\left\{\begin{aligned}
A_1z_1^{n+1}={}&z_1^n+\tau[(-m+\alpha v_2^n)\circ z_1^n
+\alpha v_1^n\circ z_2^n],\\
A_2z_2^{n+1}={}&z_2^n+\tau\big[-\beta v_2^n\circ z_1^n
+(r-\beta v_1^n-qu^n)\circ z_2^n-qv_2^n\circ d^n\big].
\end{aligned}\right.
\end{equation}
With $\xi\coloneqq Dp(u)d$ obtained from the differentiated adjoint \eqref{eq:xiI}, the Hessian action is
\begin{equation}\label{eq:hessI}
(H^{\mathrm I}d)^n=\lambda d^n-\eta qz_2^n
-q(z_2^n\circ p_2^{n+1}+v_2^n\circ\xi_2^{n+1}).
\end{equation}
To derive these formulas, fix the initial vectors and define the residuals $F_i^n\in\R^{M+1}$ by
\begin{equation*}
F_1^n\coloneqq A_1v_1^{n+1}-v_1^n-\tau v_1^n\circ(-m+\alpha v_2^n),\qquad
F_2^n\coloneqq A_2v_2^{n+1}-v_2^n-\tau v_2^n\circ(r-\beta v_1^n-qu^n).
\end{equation*}
The \emph{Lagrangian} is
\begin{equation}\label{eq:lagI}
\mapdef{\mathcal L^{\rm I}}{\Yht\times\R^{N(M+1)}\times\R^{2N(M+1)}}{\R}{(v,u,p)}{
 J_{h,\tau}^0(v,u)-\sum_{n=0}^{N-1}\sum_{i=1}^2\ip{p_i^{n+1}}{F_i^n}}.
\end{equation}
For arbitrary state variations $\delta v$ and a control variation $d$, differentiation gives
\begin{align*}
DF_1^n[\delta v,d]
&=A_1\delta v_1^{n+1}-\delta v_1^n
-\tau\bigl[(-m+\alpha v_2^n)\circ\delta v_1^n
+\alpha v_1^n\circ\delta v_2^n\bigr],\\
DF_2^n[\delta v,d]
&=A_2\delta v_2^{n+1}-\delta v_2^n
-\tau\bigl[-\beta v_2^n\circ\delta v_1^n
+(r-\beta v_1^n-qu^n)\circ\delta v_2^n
-qv_2^n\circ d^n\bigr].
\end{align*}
Weighted symmetry and componentwise multiplication give
\[
\langle p_i^{n+1},A_i\delta v_i^{n+1}\rangle_h
=\langle A_ip_i^{n+1},\delta v_i^{n+1}\rangle_h,\qquad
\langle p,a\circ\delta v\rangle_h=\langle a\circ p,\delta v\rangle_h.
\]
The time index is shifted by
\[
\sum_{n=0}^{N-1}\langle A_ip_i^{n+1},\delta v_i^{n+1}\rangle_h
=\sum_{k=1}^{N}\langle A_ip_i^k,\delta v_i^k\rangle_h.
\]
With $v_i^0$ fixed, write $D_v\mathcal L^{\rm I}[\delta v]=\sum_{i=1}^2\sum_{k=1}^N\ip{c_i^k}{\delta v_i^k}$, $\delta v_i^k\in\R^{M+1}$. At $k=N$,
\begin{equation*}
c_i^N=-A_ip_i^N=0\quad\Longrightarrow\quad p_i^N=0.
\end{equation*}
For $1\le k<N$, the cost and the residuals at steps $k-1,k$ give
\[
c_1^k={}\tau\gamma_1(v_1^k-v_1^{d,k})-A_1p_1^k+p_1^{k+1} +\tau(-m+\alpha v_2^k)\circ p_1^{k+1}-\tau\beta v_2^k\circ p_2^{k+1},
\]
and
\[
c_2^k={}\tau[\gamma_2(v_2^k-v_2^{d,k})-\eta qu^k]-A_2p_2^k+p_2^{k+1} +\tau\alpha v_1^k\circ p_1^{k+1} +\tau(r-\beta v_1^k-qu^k)\circ p_2^{k+1}.
\]
Setting $c_i^k=0$ gives \eqref{eq:adjI}. With states fixed,
\begin{equation*}
D_uF_1^n[d^n]=0,\qquad D_uF_2^n[d^n]=\tau qv_2^n\circ d^n.
\end{equation*}
Thus
\[
D_u\mathcal L^{\rm I}[d] =\tau\sum_{n=0}^{N-1}\ip{\lambda u^n-\eta qv_2^n}{d^n} -\tau q\sum_{n=0}^{N-1}\ip{p_2^{n+1}}{v_2^n\circ d^n} =\tau\sum_{n=0}^{N-1}\ip{\lambda u^n-\eta qv_2^n-qv_2^n\circ p_2^{n+1}}{d^n}.
\]
This gives \eqref{eq:gradI}.
Differentiating $F_i^n(v(u),u)=0$ gives \eqref{eq:tangentI}, while differentiating $c_i^n=0$ gives $\xi_i^N=0$ and
\begin{equation}\label{eq:xiI}
\left\{\begin{aligned}
A_1\xi_1^n={}&\xi_1^{n+1}+\tau\big[\gamma_1z_1^n+\alpha z_2^n\circ p_1^{n+1}
 +(-m+\alpha v_2^n)\circ\xi_1^{n+1}-\beta z_2^n\circ p_2^{n+1}-\beta v_2^n\circ\xi_2^{n+1}\big],
\\
A_2\xi_2^n={}&\xi_2^{n+1}+\tau\big[\gamma_2z_2^n-\eta qd^n +\alpha z_1^n\circ p_1^{n+1}+\alpha v_1^n\circ\xi_1^{n+1}\\
&+(-\beta z_1^n-qd^n)\circ p_2^{n+1} +(r-\beta v_1^n-qu^n)\circ\xi_2^{n+1}\big].
\end{aligned}\right.
\end{equation}
Then
\begin{equation*}
D(v_2^n\circ p_2^{n+1})[d]
=z_2^n\circ p_2^{n+1}+v_2^n\circ\xi_2^{n+1},
\end{equation*}
so \eqref{eq:gradI} gives \eqref{eq:hessI}.
\subsection{Positive step and endpoint-dependent harvest}
Define $G_i^n\in\R^{M+1}$ by
\begin{equation*}
G_2^n\coloneqq C^nv_2^{n+1}-(1+\tau r)v_2^n,\qquad
G_1^n\coloneqq Bv_1^{n+1}-v_1^n-\tau\alpha v_1^n\circ v_2^{n+1},
\end{equation*}
and set
\begin{equation}\label{eq:lagP}
\mapdef{\mathcal L^{+,\zeta}}{\Yht\times\R^{N(M+1)}\times\R^{2N(M+1)}}{\R}{(v,u,p)}{
 J_{h,\tau}^\zeta(v,u)-\sum_{n=0}^{N-1}\sum_{i=1}^2\ip{p_i^{n+1}}{G_i^n}}.
\end{equation}
For arbitrary variations, the positive residuals give
\begin{align*}
DG_2^n[\delta v,d]
&=C^n\delta v_2^{n+1}-(1+\tau r)\delta v_2^n
+\tau(\beta\delta v_1^n+qd^n)\circ v_2^{n+1},\\
DG_1^n[\delta v,d]
&=B\delta v_1^{n+1}-\delta v_1^n
-\tau\alpha(\delta v_1^n\circ v_2^{n+1}
+v_1^n\circ\delta v_2^{n+1}).
\end{align*}
The harvest variation with respect to the prey trajectory is
\begin{equation*}
-\eta q\tau\sum_{n=0}^{N-1}
\langle u^n,(1-\zeta)\delta v_2^n+\zeta\delta v_2^{n+1}\rangle_h
=-\eta q\tau\sum_{k=1}^{N-1}
\langle(1-\zeta)u^k+\zeta u^{k-1},\delta v_2^k\rangle_h
-\eta q\tau\zeta\langle u^{N-1},\delta v_2^N\rangle_h,
\end{equation*}
because $\delta v_2^0=0$. This identity accounts for both time levels in the adjoint source and for the terminal term. The coefficients of the terminal variations are
\begin{equation*}
c_1^N=-Bp_1^N,\qquad c_2^N=-\tau\zeta\eta qu^{N-1}-C^{N-1}p_2^N
+\tau\alpha v_1^{N-1}\circ p_1^N.
\end{equation*}
Setting these to zero gives \eqref{eq:terminalP}. For $1\le k<N$,
\[
c_1^k={}\tau\gamma_1(v_1^k-v_1^{d,k})-Bp_1^k+p_1^{k+1} +\tau v_2^{k+1}\circ(\alpha p_1^{k+1}-\beta p_2^{k+1}),
\]
and
\[
c_2^k={}\tau\gamma_2(v_2^k-v_2^{d,k}) -\tau\eta q[(1-\zeta)u^k+\zeta u^{k-1}] -C^{k-1}p_2^k+(1+\tau r)p_2^{k+1} +\tau\alpha v_1^{k-1}\circ p_1^k.
\]
The coefficient $c_2^k$ includes harvest at both steps $k$ and $k-1$. Setting $c_i^k=0$ gives \eqref{eq:adjP}. Holding the states fixed, we have
\begin{equation*}
D_uG_2^n[d^n]=\tau qd^n\circ v_2^{n+1},\qquad D_uG_1^n[d^n]=0,
\end{equation*}
so
\[
D_u\mathcal L^{+,\zeta}[d] ={}\tau\sum_{n=0}^{N-1}\ip{\lambda u^n-\eta q[(1-\zeta)v_2^n+\zeta v_2^{n+1}]}{d^n} -\tau q\sum_{n=0}^{N-1}\ip{v_2^{n+1}\circ p_2^{n+1}}{d^n}.
\]
This proves \eqref{eq:gradP}.
For $d\in\R^{N(M+1)}$, the product rule gives
\begin{equation*}
D(C^nv_2^{n+1})[d]=C^nz_2^{n+1}
+\tau(\beta z_1^n+qd^n)\circ v_2^{n+1},
\end{equation*}
and
\begin{equation*}
D(C^{k-1}p_2^k)[d]=C^{k-1}\xi_2^k
+\tau(\beta z_1^{k-1}+qd^{k-1})\circ p_2^k,
\end{equation*}
and
\begin{equation*}
D(v_1^n\circ v_2^{n+1})[d]
=z_1^n\circ v_2^{n+1}+v_1^n\circ z_2^{n+1}.
\end{equation*}
These identities yield \eqref{eq:tangentP} and \eqref{eq:terminalxi}--\eqref{eq:xiP}. Differentiating \eqref{eq:gradP} gives \eqref{eq:hessP}.
\subsection{Positive tangent and tangent-adjoint equations}\label{app:positive-derivatives}
For $d\in\R^{N(M+1)}$, set $z_i^0=0$. Since
\begin{equation*}
DC^n(u)[d]=\tau\diag(\beta z_1^n+qd^n),
\end{equation*}
the tangent satisfies
\begin{equation}\label{eq:tangentP}
\left\{\begin{aligned}
C^n z_2^{n+1}={}&(1+\tau r)z_2^n
-\tau(\beta z_1^n+qd^n)\circ v_2^{n+1},\\
Bz_1^{n+1}={}&z_1^n+\tau\alpha(z_1^n\circ v_2^{n+1}
+v_1^n\circ z_2^{n+1}).
\end{aligned}\right.
\end{equation}
Differentiating the terminal equation gives
\begin{equation}\label{eq:terminalxi}
\xi_1^N=0,\qquad
C^{N-1}\xi_2^N=-\tau\zeta\eta qd^{N-1}
-\tau(\beta z_1^{N-1}+qd^{N-1})\circ p_2^N.
\end{equation}
The backward derivative is
\begin{equation}\label{eq:xiP}
\left\{\begin{aligned}
B\xi_1^k={}&\xi_1^{k+1}+\tau\big[\gamma_1z_1^k
+z_2^{k+1}\circ(\alpha p_1^{k+1}-\beta p_2^{k+1})+v_2^{k+1}\circ(\alpha\xi_1^{k+1}-\beta\xi_2^{k+1})\big],\\
C^{k-1}\xi_2^k={}&(1+\tau r)\xi_2^{k+1}
+\tau\big[\gamma_2z_2^k-\eta q((1-\zeta)d^k+\zeta d^{k-1})\\
&\hspace{5.5em}+\alpha z_1^{k-1}\circ p_1^k+\alpha v_1^{k-1}\circ\xi_1^k
-(\beta z_1^{k-1}+qd^{k-1})\circ p_2^k\big].
\end{aligned}\right.
\end{equation}
Then
\begin{equation}\label{eq:hessP}
(H^{+,\zeta}d)^n={}\lambda d^n-\eta q[(1-\zeta)z_2^n+\zeta z_2^{n+1}] -q(z_2^{n+1}\circ p_2^{n+1}+v_2^{n+1}\circ\xi_2^{n+1}).
\end{equation}
\subsection{Second variations of the discrete costs}
Fix either state scheme and take control directions $d,e\in\R^{N(M+1)}$. Write $z[d]=DS(u)[d]$, $z[e]=DS(u)[e]$ and $s[d,e]=D^2S(u)[d,e]$. All three variations have zero initial value. For IMEX, differentiating \eqref{eq:tangentI} in direction $e$ gives
\[
\left\{\begin{aligned}
A_1s_1^{n+1}
={}&s_1^n+\tau\bigl[(-m+\alpha v_2^n)\circ s_1^n+\alpha v_1^n\circ s_2^n\bigr]+\tau\alpha\bigl(z_1[d]^n\circ z_2[e]^n+z_1[e]^n\circ z_2[d]^n\bigr),\\
A_2s_2^{n+1}
={}&s_2^n+\tau\bigl[-\beta v_2^n\circ s_1^n+(r-\beta v_1^n-qu^n)\circ s_2^n\bigr]\\
&-\tau\beta\bigl(z_1[d]^n\circ z_2[e]^n+z_1[e]^n\circ z_2[d]^n\bigr)
-\tau q\bigl(d^n\circ z_2[e]^n+e^n\circ z_2[d]^n\bigr).
\end{aligned}\right.
\]
For the positive scheme, differentiating \eqref{eq:tangentP} instead gives
\[
\left\{\begin{aligned}
C^ns_2^{n+1}
={}&(1+\tau r)s_2^n-\tau\beta s_1^n\circ v_2^{n+1}
-\tau(\beta z_1[d]^n+qd^n)\circ z_2[e]^{n+1}
-\tau(\beta z_1[e]^n+qe^n)\circ z_2[d]^{n+1},\\
Bs_1^{n+1}
={}&s_1^n+\tau\alpha\bigl(s_1^n\circ v_2^{n+1}+v_1^n\circ s_2^{n+1}\bigr)+\tau\alpha\bigl(z_1[d]^n\circ z_2[e]^{n+1}+z_1[e]^n\circ z_2[d]^{n+1}\bigr).
\end{aligned}\right.
\]
The cost derivatives follow without differentiating the adjoint. With $\zeta=0$ for IMEX, the first derivative is
\begin{align*}
Dj(u)[d]
={}\lambda\langle u,d\rangle_{h,\tau}
&+\tau\sum_{n=0}^{N-1}\sum_{i=1}^2\gamma_i
\langle v_i^n-v_i^{d,n},z_i[d]^n\rangle_h\\
&-\eta q\tau\sum_{n=0}^{N-1}\Bigl[
\langle d^n,(1-\zeta)v_2^n+\zeta v_2^{n+1}\rangle_h
+\langle u^n,(1-\zeta)z_2[d]^n+\zeta z_2[d]^{n+1}\rangle_h\Bigr].
\end{align*}
Differentiating each factor gives
\begin{align*}
D^2j(u)[d,e]
={}\lambda\langle d,e\rangle_{h,\tau}
&+\tau\sum_{n=0}^{N-1}\sum_{i=1}^2\gamma_i
\Bigl[\langle z_i[d]^n,z_i[e]^n\rangle_h
+\langle v_i^n-v_i^{d,n},s_i[d,e]^n\rangle_h\Bigr]\\
&-\eta q\tau\sum_{n=0}^{N-1}
\Bigl[\langle d^n,(1-\zeta)z_2[e]^n+\zeta z_2[e]^{n+1}\rangle_h\\
&\hspace{5em}+\langle e^n,(1-\zeta)z_2[d]^n+\zeta z_2[d]^{n+1}\rangle_h\\
&\hspace{5em}+\langle u^n,(1-\zeta)s_2[d,e]^n+\zeta s_2[d,e]^{n+1}\rangle_h\Bigr].
\end{align*}
Both second-variation systems are symmetric in $d,e$, so uniqueness gives $s[d,e]=s[e,d]$. Thus the bilinear form is symmetric. Differentiating $Dj(u)[e]=\langle g(u),e\rangle_{h,\tau}$ also gives
\[
D^2j(u)[d,e]=\langle Dg(u)[d],e\rangle_{h,\tau}
=\langle H(u)d,e\rangle_{h,\tau},
\]
which identifies the same Hessian as \eqref{eq:hessI} and \eqref{eq:hessP}.
\subsection{Symmetry and local smoothness}
For an invertible matrix $A\in\R^{s\times s}$ and perturbations $E,F\in\R^{s\times s}$, differentiating $AA^{-1}=I_s$ gives
\begin{equation*}
D(A^{-1})[E]=-A^{-1}EA^{-1},\qquad D^2(A^{-1})[E,F]=A^{-1}EA^{-1}FA^{-1} +A^{-1}FA^{-1}EA^{-1}.
\end{equation*}
Each IMEX step combines polynomials with fixed matrix inverses. For the positive scheme, the loss matrices remain invertible near the admissible box, where their inverses depend smoothly on the entries. All three reduced costs are therefore smooth near $\Uk$ on each fixed grid. The weighted gradient and the Euclidean derivatives are related by
\[
Dj(u)[d]=g(u)^T\Wht d=\nabla j(u)^Td,
\qquad \nabla j(u)=\Wht g(u),\qquad \nabla^2j(u)=\Wht H(u).
\]
Since $\Wht$ is constant,
\[
D^2j(u)[d,e]=e^T\Wht Dg(u)d =e^T\Wht H(u)d =D^2j(u)[e,d].
\]
Thus
\begin{equation*}
H(u)^T\Wht=\Wht H(u),\qquad \big(\Wht^{1/2}H(u)\Wht^{-1/2}\big)^T=\Wht^{1/2}H(u)\Wht^{-1/2}.
\end{equation*}
These identities show that $H$ is self-adjoint in the weighted inner product and $\Wht^{1/2}H\Wht^{-1/2}$ is symmetric. This proves Proposition~\ref{prop:symmetry}.
\section{Compactness and variational convergence proofs}\label{app:convergence}
\subsection{Proof of Lemma~\ref{lem:energy}}
Both schemes have the form
\begin{equation}\label{eq:energyform}
D_\tau v_i^n-\mu_iL_hv_i^{n+1}=f_i^n.
\end{equation}
For IMEX, $f_i^n\coloneqq f_i(v^n,u^n)\in\R^{M+1}$. For the positive scheme,
\begin{equation*}
f_1^n\coloneqq-mv_1^{n+1}+\alpha v_1^n\circ v_2^{n+1},\qquad
f_2^n\coloneqq rv_2^n-(\beta v_1^n+qu^n)\circ v_2^{n+1}.
\end{equation*}
For both schemes, Proposition~\ref{prop:uniform} gives
\begin{equation*}
\|f_1^n\|_\infty\le mK_1+\alpha K_1K_2,\qquad
\|f_2^n\|_\infty\le rK_2+(\beta K_1+qU_{\max})K_2.
\end{equation*}
Since $\sum_{j=0}^M w_j=1$, we have
\begin{equation}\label{eq:forcingbound}
\max_{\substack{i\in\{1,2\}\\0\le n<N}}\|f_i^n\|_\infty\le C,
\qquad \tau\sum_{n=0}^{N-1}\normh{f_i^n}^2\le CT.
\end{equation}
These estimates use only the bounds on $u$.
By symmetry of $a_h$ and the definition of $D_\tau$,
\begin{equation}\label{eq:pfa_energy_polarization}
2\tau a_h(v_i^{n+1},D_\tau v_i^n)
=|v_i^{n+1}|_{1,h}^2-|v_i^n|_{1,h}^2+\tau^2|D_\tau v_i^n|_{1,h}^2.
\end{equation}
Test \eqref{eq:energyform} with $2\tau D_\tau v_i^n$. Young's inequality gives
\begin{align*}
2\tau\normh{D_\tau v_i^n}^2+2\tau\mu_i a_h(v_i^{n+1},D_\tau v_i^n)
&=\underbrace{2\tau\ip{f_i^n}{D_\tau v_i^n}}
_{\substack{
\le2\tau\normh{f_i^n}\normh{D_\tau v_i^n}\\
\le\tau\normh{f_i^n}^2+\tau\normh{D_\tau v_i^n}^2
}}.
\end{align*}
Subtract $\tau\normh{D_\tau v_i^n}^2$ from the resulting inequality, substitute \eqref{eq:pfa_energy_polarization}, and sum over $0\le n<\ell$. The spatial energies telescope:
\begin{equation}\label{eq:energystrong}
\tau\sum_{n=0}^{\ell-1}\normh{D_\tau v_i^n}^2
+\mu_i|v_i^\ell|_{1,h}^2+\mu_i\tau^2\sum_{n=0}^{\ell-1}|D_\tau v_i^n|_{1,h}^2
\le\mu_i|v_i^0|_{1,h}^2+\tau\sum_{n=0}^{\ell-1}\normh{f_i^n}^2.
\end{equation}
For nodal initial data, Cauchy--Schwarz on each grid interval gives
\begin{align*}
|v_i^0|_{1,h}^2
&=\sum_{j=0}^{M-1}h\left|\frac1h\int_{x_j}^{x_{j+1}}\partial_xv_i^0\,\mathrm{d}x\right|^2\le\sum_{j=0}^{M-1}\int_{x_j}^{x_{j+1}}|\partial_xv_i^0|^2\,\mathrm{d}x
=\|\partial_xv_i^0\|_2^2.
\end{align*}
With $\|\Ih v_i^\ell\|_2\le K_i$, this gives \eqref{eq:energy}.
The initial-data bound and the forcing estimate give
\[
\tau\sum_{n=0}^{N-1}\|D_\tau v_i^n\|_h^2\le C,\qquad
\max_{0\le\ell\le N}|v_i^\ell|_{1,h}^2\le C.
\]
The constants depend on the initial $H^1$ norm, the uniform reaction bounds, $T$ and $\mu_i$, with no factor involving $h^{-1}$ or $\tau^{-1}$. The difference between the time reconstructions satisfies
\begin{align*}
\|V_i^+-V_i^-\|_{L^2(Q)}^2
&=\tau\sum_{n=0}^{N-1}\normh{v_i^{n+1}-v_i^n}^2
=\tau^3\sum_{n=0}^{N-1}\normh{D_\tau v_i^n}^2
\byeq{eq:energystrong}{\le}C\tau^2.
\end{align*}
For the spatial reconstructions, integrate over both halves of each grid interval:
\begin{equation*}
\int_0^{h/2}s^2\,\mathrm{d}s+\int_{h/2}^h(h-s)^2\,\mathrm{d}s=\frac{h^3}{12}.
\end{equation*}
Thus
\[
\|\Qh y-\Ih y\|_2^2 =\sum_{j=0}^{M-1}\frac{(y_{j+1}-y_j)^2}{h^2}\,\frac{h^3}{12} =\frac h{12}\sum_{j=0}^{M-1}(y_{j+1}-y_j)^2 \byeq{eq:bilinear_spatial}{=}\frac{h^2}{12}a_h(y,y),
\]
which proves \eqref{eq:reconstructionerror}.
\subsection{Proof of Theorem~\ref{th:statecompact}}
For $y\in\R^{M+1}$,
\begin{align*}
\|\Ih y\|_2^2
=\sum_{j=0}^{M-1}\frac h3(y_j^2+y_jy_{j+1}+y_{j+1}^2)
\byeq{eq:bilinear_spatial}{=}\normh y^2-\frac{h^2}{6}a_h(y,y).
\end{align*}
Since $y_j^2+y_jy_{j+1}+y_{j+1}^2\ge\tfrac12(y_j^2+y_{j+1}^2)$, we have
\begin{equation}\label{eq:normequiv}
\tfrac13\normh y^2\le\|\Ih y\|_{L^2(0,1)}^2\le\normh y^2.
\end{equation}
The energy estimate gives
\begin{equation*}
\|\partial_t\widehat V_i\|_{L^2(Q)}^2
=\tau\sum_{n=0}^{N-1}\|\Ih D_\tau v_i^n\|_2^2
\byeq{eq:energystrong}{\le}C,\qquad
\|\widehat V_i\|_{L^\infty(0,T;H^1)}\le C.
\end{equation*}
The compact embedding $H^1(0,1)\hookrightarrow L^2(0,1)$
and time compactness \cite{Simon1987} give a subsequence
that converges strongly in $L^2(Q)$:
$\widehat V_{i,k}\xrightarrow[k\to\infty]{} v_i.$ On $[t_n,t_{n+1})$, set $\widetilde V_i^-(t)\coloneqq\Ih v_i^n$, $\widetilde V_i^+(t)\coloneqq\Ih v_i^{n+1}$. These functions satisfy $\widetilde V_i^\pm\in L^\infty(0,T;\Vh)$ and
\begin{align*}
\|\widehat V_i-\widetilde V_i^-\|_{L^2(Q)}^2
=\frac{\tau^3}{3}\sum_{n=0}^{N-1}\|\Ih D_\tau v_i^n\|_2^2
\byeq{eq:energystrong}{\le}C\tau^2,\qquad
\|\widetilde V_i^\pm-V_i^\pm\|_{L^2(Q)}
\overset{\mathclap{\text{\eqref{eq:reconstructionerror}, \eqref{eq:energy}}}}{\le}Ch.
\end{align*}
The same time bound holds for $\widetilde V_i^+$. Thus
\begin{equation*}
\|V_{i,k}^\pm-\widehat V_{i,k}\|_{L^2(Q)}
\le\underbrace{\|V_{i,k}^\pm-\widetilde V_{i,k}^\pm\|_{L^2(Q)}}_{\le Ch_k}
+\underbrace{\|\widetilde V_{i,k}^\pm-\widehat V_{i,k}\|_{L^2(Q)}}_{\le C\tau_k} \le C(h_k+\tau_k).
\end{equation*}
Every reconstruction therefore has the same strong limit $v_i$, and $\widetilde V_{i,k}^\pm\xrightharpoonup[k\to\infty]{L^2(0,T;H^1(0,1))}v_i$.
For $\phi\in C^\infty(\overline Q)$ with $\phi(\cdot,T)=0$, set $\phi_j^n\coloneqq\phi(x_j,t_n)$. Summation by parts gives
\begin{equation}\label{eq:time_sbp_explicit}
\tau\sum_{n=0}^{N-1}\ip{D_\tau v_i^n}{\phi^{n+1}}
=-\ip{v_i^0}{\phi^0}-\tau\sum_{n=0}^{N-1}\ip{v_i^n}{D_\tau\phi^n}.
\end{equation}
Testing \eqref{eq:energyform} with $\phi^{n+1}$ therefore yields
\begin{align*}
-\tau\sum_{n=0}^{N-1}\ip{v_i^n}{D_\tau\phi^n}
 +\mu_i\tau\sum_{n=0}^{N-1} a_h(v_i^{n+1},\phi^{n+1})
\byeq{eq:time_sbp_explicit}{=}
 \tau\sum_{n=0}^{N-1}\ip{f_i^n}{\phi^{n+1}}+\ip{v_i^0}{\phi^0}.
\end{align*}
The mass-lumped time term is
\[
\tau_k\sum_{n=0}^{N_k-1}
\langle v_{i,k}^n,D_{\tau_k}\phi^n\rangle_{h_k}
=\int_Q V_{i,k}^-\,\mathcal Q_{h_k,\tau_k}(D_{\tau_k}\phi),
\]
where $N_k=T/\tau_k$ is the number of time steps.
Smoothness of $\phi$ gives uniform convergence of the reconstructed difference quotient to $\partial_t\phi$. Therefore
\begin{equation*}
\left|\int_QV_{i,k}^-\,\mathcal Q_{h_k,\tau_k}(D_{\tau_k}\phi)
-\int_Qv_i\partial_t\phi\right|
\le{}\|V_{i,k}^--v_i\|_2
\|\mathcal Q_{h_k,\tau_k}(D_{\tau_k}\phi)\|_2 +\|v_i\|_2
\|\mathcal Q_{h_k,\tau_k}(D_{\tau_k}\phi)-\partial_t\phi\|_2\xrightarrow[k\to\infty]{}0.
\end{equation*}
For the diffusion term, let $\widetilde\phi_k^+(t)=\mathcal I_{h_k}\phi^{n+1}$ on $[t_n,t_{n+1})$. Then
\[
\tau_k\sum_{n=0}^{N_k-1}
a_{h_k}(v_{i,k}^{n+1},\phi^{n+1})
=\int_Q\partial_x\widetilde V_{i,k}^+\,\partial_x\widetilde\phi_k^+.
\]
Since $\partial_x\widetilde V_{i,k}^+$ converges weakly in $L^2(Q)$ and $\partial_x\widetilde\phi_k^+$ converges uniformly to $\partial_x\phi$,
\begin{equation*}
\int_Q\partial_x\widetilde V_{i,k}^+\,\partial_x\widetilde\phi_k^+
-\int_Q\partial_xv_i\,\partial_x\phi =\int_Q\partial_x\widetilde V_{i,k}^+
(\partial_x\widetilde\phi_k^+-\partial_x\phi)
+\int_Q(\partial_x\widetilde V_{i,k}^+-\partial_xv_i)\partial_x\phi\xrightarrow[k\to\infty]{}0.
\end{equation*}
The initial term is
\[
\langle v_i^0,\phi^0\rangle_{h_k}
=\int_0^1\mathcal Q_{h_k}v_i^0\,\mathcal Q_{h_k}\phi^0
\xrightarrow[k\to\infty]{}\int_0^1v_i^0(x)\phi(x,0)\,\mathrm{d}x,
\]
because the nodal reconstruction of the fixed $H^1(0,1)$ initial datum converges in $L^2(0,1)$.
The nonlinear reaction products converge strongly, since
\begin{equation*}
\|V_{1,k}^-V_{2,k}^+-v_1v_2\|_2
\le K_1\|V_{2,k}^+-v_2\|_2+K_2\|V_{1,k}^--v_1\|_2\xrightarrow[k\to\infty]{}0.
\end{equation*}
For $U_k=\mathcal Q_{h_k,\tau_k}u_k$ and $\psi\in L^2(Q)$,
\begin{align*}
\int_Q(U_kV_{2,k}^+-uv_2)\psi
={}\int_QU_k(V_{2,k}^+-v_2)\psi+\int_Q(U_k-u)v_2\psi
\xrightarrow[k\to\infty]{}0.
\end{align*}
For the first term,
\[
\left|\int_QU_k(V_{2,k}^+-v_2)\psi\right|
\le U_{\max}\|V_{2,k}^+-v_2\|_2\|\psi\|_2\xrightarrow[k\to\infty]{}0,
\]
and $v_2\psi\in L^1(Q)$ gives the second limit by weak-$*$ convergence. By \eqref{eq:jump}, using the state at either endpoint of each time step gives the same weak identity:
\begin{equation*}
-\int_Qv_i\partial_t\phi+\mu_i\int_Q\partial_xv_i\partial_x\phi
=\int_Qf_i(v,u)\phi+\int_0^1v_i^0\phi(\cdot,0).
\end{equation*}
The nonnegative limit solves \eqref{eq:pde} with the Neumann conditions and prescribed initial data, so uniqueness gives $v=S(u)$. Every convergent subsequence has this limit. Relative compactness then gives convergence of the whole sequence.
\subsection{Recovery controls and the liminf inequality}
On each space--time cell $D_j^n\coloneqq C_j\times[t_n,t_{n+1})$, Jensen's inequality gives
\[
\|\Qht\avg f\|_2^2 =\sum_{n=0}^{N-1}\sum_{j=0}^M|D_j^n|\left|\frac1{|D_j^n|}\int_{D_j^n}f\right|^2 \le\sum_{n=0}^{N-1}\sum_{j=0}^M\int_{D_j^n}|f|^2 =\|f\|_2^2.
\]
Cell averages preserve $[0,U_{\max}]$. For $f\in L^2(Q)$, choose $f_\varepsilon\in C(\overline Q)$ with $\|f-f_\varepsilon\|_2\le\varepsilon$. Then
\begin{align*}
\|\Qht\avg f-f\|_2
&\le2\|f-f_\varepsilon\|_2
+\|\Qht\avg f_\varepsilon-f_\varepsilon\|_2.
\end{align*}
Uniform continuity bounds the last term by
\[
 |Q|^{1/2}\sup_{\substack{z,z'\in\overline Q\\|z-z'|\le\sqrt{h^2+\tau^2}}}
 |f_\varepsilon(z)-f_\varepsilon(z')|\xrightarrow[h,\tau\to0]{}0.
\]
Let $h,\tau\to0$, then $\varepsilon\to0$.
For bounded controls, strong $L^2$ convergence also implies weak-$*$ convergence. Indeed, for $\psi\in L^1(Q)$ choose $\psi_\varepsilon\in L^2(Q)$ with $\|\psi-\psi_\varepsilon\|_1\le\varepsilon$. This gives
\begin{equation*}
\left|\int_Q(U_k-u)\psi\right|
\le2U_{\max}\|\psi-\psi_\varepsilon\|_1
+\|U_k-u\|_2\|\psi_\varepsilon\|_2.
\end{equation*}
Let $k\to\infty$, then $\varepsilon\to0$.
Here $N=T/\tau_k$ and $M=1/h_k$. With target array $(v_i^{d,n})_{n<N}\in\R^{N(M+1)}$ and $\zeta=0$ for IMEX, separate the effort term from the cost:
\[
A_k(u_k)\coloneqq j_k(u_k)-\frac\lambda2\|U_k\|_2^2 =\frac12\sum_{i=1}^2\gamma_i\|V_{i,k}^--\mathcal Q_{h_k,\tau_k}\big((v_i^{d,n})_{n<N}\big)\|_2^2 -\eta q\int_QU_k\big[(1-\zeta)V_{2,k}^-+\zeta V_{2,k}^+\big].
\]
Define
\[
 \mapdef{A}{\Uad}{\R}{u}{j(u)-\lambda\|u\|_2^2/2}.
\]
For the tracking terms, write $V_{i,k}^d=\mathcal Q_{h_k,\tau_k}((v_i^{d,n})_{n<N})$. Then
\[
\|(V_{i,k}^--V_{i,k}^d)-(v_i-v_i^d)\|_2
\le\|V_{i,k}^--v_i\|_2+\|V_{i,k}^d-v_i^d\|_2\xrightarrow[k\to\infty]{}0.
\]
For either endpoint of the harvest quadrature,
\[
\left|\int_QU_kV_{2,k}^\pm-\int_Quv_2\right|
\le U_{\max}|Q|^{1/2}\|V_{2,k}^\pm-v_2\|_2
+\left|\int_Q(U_k-u)v_2\right|\xrightarrow[k\to\infty]{}0.
\]
Thus the state, target and harvest terms converge, and whenever $U_k\mathrel{\underset{k\to\infty}{\rightharpoonup^{*}}}u$,
\begin{equation}\label{eq:Aconvergence}
A_k(u_k)\xrightarrow[k\to\infty]{} A(u).
\end{equation}
Lower semicontinuity of the effort term then gives
\[
\liminf_k j_k(u_k) \byeq{eq:Aconvergence}{=}A(u)+\frac\lambda2\liminf_k\|U_k\|_2^2 \ge A(u)+\frac\lambda2\|u\|_2^2 =j(u).
\]
For the extended cost, apply the argument to a subsequence attaining the liminf if it is finite. Otherwise, \eqref{eq:liminf} follows.
For $u_k=\mathcal A_{h_k,\tau_k}u$, strong $L^2$ convergence gives weak-$*$ convergence and convergence of the effort term. Then \eqref{eq:Aconvergence} gives
\begin{equation*}
j_k(u_k)=A_k(u_k)+\frac\lambda2\|U_k\|_2^2\xrightarrow[k\to\infty]{} j(u),
\end{equation*}
which proves \eqref{eq:recovery}.
\subsection{Minimum values and strong convergence of minimizers}
Compactness of the discrete box ensures that $m_k\coloneqq\min_{\mathcal U_{\mathrm{ad},h_k,\tau_k}}j_k\in\R$ is attained. Set $m_\star\coloneqq\min_{\Uad}j$. Recovery controls for a global minimizer and the liminf inequality give the two bounds
\begin{equation*}
\limsup_km_k\le m_\star\le\liminf_km_k
\quad\Longrightarrow\quad m_k\xrightarrow[k\to\infty]{} m_\star.
\end{equation*}
For \eqref{eq:epsglobal},
\begin{equation*}
m_k\le j_k(u_k)\le m_k+\varepsilon_k
\quad\Longrightarrow\quad j_k(u_k)\xrightarrow[k\to\infty]{} m_\star.
\end{equation*}
Take a subsequence $U_k\mathrel{\underset{k\to\infty}{\rightharpoonup^{*}}}u^*$. The liminf inequality gives $j(u^*)=m_\star$, and \eqref{eq:Aconvergence} gives
\begin{align*}
\frac\lambda2\|U_k\|_2^2
=j_k(u_k)-A_k(u_k)
\byeq{eq:Aconvergence}{\xrightarrow[k\to\infty]{}}j(u^*)-A(u^*)
=\frac\lambda2\|u^*\|_2^2.
\end{align*}
Since $\lambda>0$ and $U_k\xrightharpoonup[k\to\infty]{L^2(Q)}u^*$,
\begin{equation*}
\|U_k-u^*\|_2^2=\|U_k\|_2^2+\|u^*\|_2^2-2(U_k,u^*)_Q\xrightarrow[k\to\infty]{}0.
\end{equation*}
This proves Theorem~\ref{th:gamma}.
Every subsequence has a further subsequence converging strongly to a global minimizer, so
\[
\inf_{w\in\argmin_{\Uad}j}\|U_k-w\|_{L^2(Q)}\xrightarrow[k\to\infty]{}0.
\]
If the global minimizer $u^*$ is unique, the full sequence converges to it. For two asymptotically minimizing reconstructions $U_k,\widetilde U_k$,
\begin{equation*}
\|U_k-\widetilde U_k\|_2
\le\|U_k-u^*\|_2+\|\widetilde U_k-u^*\|_2\xrightarrow[k\to\infty]{}0,
\end{equation*}
which proves Corollary~\ref{cor:unique}.
\subsection{Proof of Corollary~\ref{cor:local}}
For fine enough grids, minimize $j_k$ over
\[
K_k\coloneqq\{u\in \mathcal U_{\mathrm{ad},h_k,\tau_k}:
\|\mathcal Q_{h_k,\tau_k}u-u^*\|_2\le r\}.
\]
The recovery controls of $u^*$ belong to $K_k$ for large $k$. This compact set is therefore nonempty. Choose a minimizer $u_k$ and write $U_k=\mathcal Q_{h_k,\tau_k}u_k$.
From every subsequence we can extract a further subsequence such that, in $L^\infty(Q)$, $U_k\mathrel{\underset{k\to\infty}{\rightharpoonup^{*}}}u$. The limit is admissible. Weak lower semicontinuity keeps it in the closed $L^2$ ball of radius $r$, while the liminf inequality and recovery controls give
\[
j(u)\le\liminf_k j_k(u_k)
\le\limsup_k j_k(u_k)\le j(u^*).
\]
Strict local minimality gives $u=u^*$, $j_k(u_k)\xrightarrow[k\to\infty]{} j(u^*)$. Then
\begin{align*}
\frac\lambda2\|U_k\|_2^2
=j_k(u_k)-A_k(u_k)
\byeq{eq:Aconvergence}{\xrightarrow[k\to\infty]{}}j(u^*)-A(u^*)
=\frac\lambda2\|u^*\|_2^2.
\end{align*}
Since $\lambda>0$, the norms converge, turning weak convergence into strong convergence. The subsequence argument gives $U_k\xrightarrow[k\to\infty]{}u^*$ for the full sequence. Eventually $\|U_k-u^*\|_2<r$: the ball constraint is locally inactive, and $u_k$ is a local minimizer over the admissible box.
\subsection{Proof of Proposition~\ref{prop:quadrature}}\label{app:proof-prop-quadrature}
\begin{proof}
Subtracting the two costs and using \eqref{eq:jump} gives
\[
|j^{+,1}_{h,\tau}(u)-j^{+,0}_{h,\tau}(u)|
=\eta q\left|\int_Q(\Qht u)(V_2^+-V_2^-)\right| \le\eta q\underbrace{\|\Qht u\|_{L^2(Q)}}_{\le U_{\max}|Q|^{1/2}}
\underbrace{\|V_2^+-V_2^-\|_{L^2(Q)}}_{\le C\tau}
\le C\tau.
\]
The bound is uniform in $u$, so taking infima gives the same bound for the minimum values.
\end{proof}
\subsection{Control comparison in Remark~\ref{rem:quadrature-controls}}\label{app:proof-growth-comparison}
\begin{proof}
Quadratic growth and the cost comparison give
\begin{equation*}
c\normt{u_1-u_0}^2
\le j_0(u_1)-j_0(u_0) =\underbrace{j_0(u_1)-j_1(u_1)}_{\le\delta}
+\underbrace{j_1(u_1)-j_1(u_0)}_{\le\varepsilon}
+\underbrace{j_1(u_0)-j_0(u_0)}_{\le\delta} \le\varepsilon+2\delta.
\end{equation*}
If $2\delta<cr^2$, a minimizer of $j_1$ on $K_r$ satisfies
\[
\normt{u_1-u_0}^2\le 2\delta/c<r^2.
\]
Thus the ball constraint is inactive and $u_1$ is a local minimizer on $\Uk$.
\end{proof}
\section{Accuracy and local optimization proofs}\label{app:local}
\subsection{Accuracy for smooth states}\label{app:smooth-accuracy}
The space $W^{1,\infty}(0,T;X)$ consists of $X$-valued functions whose norm and weak time-derivative norm are essentially bounded.
\begin{proposition}[Accuracy for compatible smooth data]\label{prop:accuracy}
Let $u\in\Uad\cap W^{1,\infty}(0,T;C([0,1]))$. Assume the states, their second time derivatives, and their fourth spatial derivatives are continuous and uniformly bounded. At $x=0,1$, assume that for every $t\in[0,T]$,
\begin{equation}\label{eq:compatibility}
\partial_xv_i=\partial_x^3v_i=0.
\end{equation}
For $u_j^n\coloneqq u(x_j,t_n)$ and nonnegative initial nodal vectors with $O(h^2)$ maximum-norm error, either scheme satisfies, for small $h,\tau$,
\begin{equation}\label{eq:accuracy}
\max_{0\le n\le N}\sum_{i=1}^2\|v_i(t_n)-\Ih v_i^n\|_{L^2(0,1)}
\le C(\tau+h^2).
\end{equation}
\end{proposition}
Condition \eqref{eq:compatibility} removes the first-order endpoint error.
\subsection{Proof of Proposition~\ref{prop:accuracy}}
Set $V_i^n\coloneqq(v_i(x_j,t_n))_{j=0}^M\in\R^{M+1}$. Inserting $V$ into the scheme, subtracting the right side from the left and dividing by $\tau$ gives the residual $r_i^n\in\R^{M+1}$. Taylor's formula yields
\begin{align*}
\frac{v_i(x,t+\tau)-v_i(x,t)}{\tau}
&=\partial_tv_i(x,t)
 +\frac1\tau\int_t^{t+\tau}(t+\tau-s)\partial_t^2v_i(x,s)\,\mathrm{d}s.
\end{align*}
Its remainder is bounded by $\frac\tau2\|\partial_t^2v_i\|_\infty$. At an interior node,
\begin{equation*}
\left|\frac{v_i(x-h,t)-2v_i(x,t)+v_i(x+h,t)}{h^2}
-\partial_x^2v_i(x,t)\right|
\le\frac{h^2}{12}\|\partial_x^4v_i\|_\infty.
\end{equation*}
The derivative norms are over $[0,1]\times[0,T]$. At the endpoints,
\begin{equation*}
\frac{2(v_i(h,t)-v_i(0,t))}{h^2} =\frac2h\partial_xv_i(0,t)+\partial_x^2v_i(0,t) +\frac h3\partial_x^3v_i(0,t)+O(h^2),
\end{equation*}
and
\begin{equation*}
\frac{2(v_i(1-h,t)-v_i(1,t))}{h^2} =-\frac2h\partial_xv_i(1,t)+\partial_x^2v_i(1,t) -\frac h3\partial_x^3v_i(1,t)+O(h^2).
\end{equation*}
The Neumann conditions cancel the $h^{-1}$ terms, and \eqref{eq:compatibility} cancels the $h$ terms. Since $|V_i^{n+1}-V_i^n|\le\tau\|\partial_tv_i\|_\infty$, the reaction quadratures differ by $O(\tau)$. For diffusion at the new time, differentiate \eqref{eq:pde} in time in the distributional sense:
\begin{equation*}
\mu_i\partial_t\partial_x^2v_i=\partial_t^2v_i-\partial_tf_i(v,u),
\end{equation*}
where
\begin{equation*}
\partial_tf_1=(-m+\alpha v_2)\partial_tv_1+\alpha v_1\partial_tv_2,\qquad
\partial_tf_2=(r-\beta v_1-qu)\partial_tv_2 -\beta v_2\partial_tv_1-qv_2\partial_tu.
\end{equation*}
The regularity assumptions bound every term. Hence $\|\partial_x^2v_i(t+\tau)-\partial_x^2v_i(t)\|_\infty\le C\tau$ and
\begin{equation*}
\max_{\substack{i\in\{1,2\}\\0\le n<N}}\|r_i^n\|_\infty\le C(\tau+h^2).
\end{equation*}
Set $e_i^n\coloneqq V_i^n-v_i^n\in\R^{M+1}$ and $E_n\coloneqq\|e_1^n\|_\infty+\|e_2^n\|_\infty$. For IMEX,
\begin{equation*}
A_ie_i^{n+1}=e_i^n+\tau[f_i(V^n,u^n)-f_i(v^n,u^n)]+\tau r_i^n.
\end{equation*}
Lipschitz continuity on the bounded range of $V^n,v^n$ and $\|A_i^{-1}\|_\infty\le1$ give
\begin{equation}\label{eq:error_recursion_explicit}
E_{n+1}\le(1+C\tau)E_n+C\tau(\tau+h^2).
\end{equation}
For the positive scheme, subtract the equations using the computed matrix $C^n$ to get
\begin{equation*}
C^ne_2^{n+1}
=(1+\tau r)e_2^n-\tau\beta e_1^n\circ V_2^{n+1}+\tau r_2^n,\qquad Be_1^{n+1} =e_1^n+\tau\alpha e_1^n\circ V_2^{n+1} +\tau\alpha v_1^n\circ e_2^{n+1}+\tau r_1^n.
\end{equation*}
The inverse bounds give
\begin{equation*}
\|e_2^{n+1}\|_\infty
\le(1+C\tau)\|e_2^n\|_\infty+C\tau\|e_1^n\|_\infty+C\tau(\tau+h^2),\qquad
\|e_1^{n+1}\|_\infty
\le(1+C\tau)\|e_1^n\|_\infty+C\tau\|e_2^{n+1}\|_\infty+C\tau(\tau+h^2).
\end{equation*}
Substituting the first bound into the second and using $\tau^2\le T\tau$ recovers \eqref{eq:error_recursion_explicit}. Iteration gives
\[
E_n \byeq{eq:error_recursion_explicit}{\le}(1+C\tau)^nE_0 +C\tau(\tau+h^2)\sum_{\ell=0}^{n-1}(1+C\tau)^\ell \le e^{CT}\big[E_0+CT(\tau+h^2)\big] \le C(\tau+h^2).
\]
The interpolation error on each interval follows from the linear interpolation remainder:
\[
|v_i(x,t_n)-\Ih V_i^n(x)|
\le\frac12(x-x_j)(x_{j+1}-x)\|\partial_x^2v_i(t_n)\|_\infty
\le\frac{h^2}{8}\|\partial_x^2v_i(t_n)\|_\infty.
\]
Also, $\Ih e_i^n$ is a convex combination of its two endpoint values, so $\|\Ih e_i^n\|_\infty\le\|e_i^n\|_\infty$. Therefore
\begin{align*}
\|v_i(t_n)-\Ih v_i^n\|_2
\le\|v_i(t_n)-\Ih V_i^n\|_2+\|\Ih e_i^n\|_2
\le Ch^2+\|e_i^n\|_\infty.
\end{align*}
This proves \eqref{eq:accuracy}.
\subsection{Proof of Theorem~\ref{th:second}}
For $d\in\mathcal C(u^*)\subset\R^{N(M+1)}$, $u^*+td$ is feasible whenever
\begin{equation*}
0\le t<\min\left(
\left\{\frac{U_{\max}-u_i^*}{d_i}:d_i>0\right\}
\cup\left\{\frac{u_i^*}{-d_i}:d_i<0\right\}\right),
\end{equation*}
with $\min\varnothing\coloneqq+\infty$ and all ratios positive. The fundamental theorem of calculus applied twice gives
\begin{align*}
j(u^*+td)-j(u^*)
&=t\langle g(u^*),d\rangle_{h,\tau}
+t^2\int_0^1(1-s)\langle H(u^*+std)d,d\rangle_{h,\tau}\,\mathrm{d}s.
\end{align*}
Continuity of $H$ makes the difference between this integral and $\frac12\langle H(u^*)d,d\rangle_{h,\tau}$ tend to zero. Since $d$ is critical, Taylor expansion gives
\[
j(u^*+td)-j(u^*) =t\ipt{g(u^*)}d+\frac{t^2}{2}\ipt{H(u^*)d}d+o(t^2) =\frac{t^2}{2}\ipt{H(u^*)d}d+o(t^2).
\]
At a local minimizer the left side is nonnegative, which proves necessity.
For sufficiency, suppose quadratic growth fails. Then feasible $u_k\xrightarrow[k\to\infty]{} u^*$, $u_k\ne u^*$, satisfy
\begin{equation*}
j(u_k)<j(u^*)+k^{-1}\normt{u_k-u^*}^2.
\end{equation*}
Set $s_k\coloneqq\normt{u_k-u^*}>0$, $d_k\coloneqq(u_k-u^*)/s_k\in\R^{N(M+1)}$. On the fixed grid, compactness gives a subsequence with
\begin{equation*}
d_k\xrightarrow[k\to\infty]{} d\in T(u^*),\qquad \normt d=1.
\end{equation*}
Stationarity gives $\ipt{g(u^*)}{d_k}\ge0$, while Taylor expansion and boundedness of $H(u^*)$ give
\[
\ipt{g(u^*)}{d_k} <s_k\left[\frac1k-\frac12\ipt{H(u^*)d_k}{d_k}+o(1)\right] \le Cs_k\xrightarrow[k\to\infty]{}0,
\]
so $d\in\mathcal C(u^*)$. Dividing by $s_k^2$,
\begin{equation*}
\frac1k>
\underbrace{\frac{\ipt{g(u^*)}{d_k}}{s_k}}_{\ge0}
+\frac12\ipt{H(u^*)d_k}{d_k}+o(1).
\end{equation*}
The first term is nonnegative, so taking limits gives
\begin{align*}
0\ge\frac12\ipt{H(u^*)d}d
\ge\frac\kappa2\normt d^2
=\frac\kappa2>0,
\end{align*}
a contradiction. For stationary $u^*$, $g_i(u^*)=0$ on $F$. Under strict complementarity,
\begin{equation*}
\ipt{g(u^*)}d=\sum_{i\in A_0\cup A_U}(\Wht)_{ii}g_i(u^*)d_i.
\end{equation*}
Every summand is nonnegative and vanishes only for $d_i=0$, so the critical cone is the free subspace.
\subsection{Proof of Theorem~\ref{th:fbs}}
For $\theta>0$, set
$y^*\coloneqq u^*-\tfrac{\theta}{\lambda}g(u^*)\in\R^{N(M+1)}$.
Stationarity on $F$ and strict complementarity on $A$ imply
\[
\begin{cases}
y_i^*<0, & \text{if }i\in A_0,\\
y_i^*>U_{\max}, & \text{if }i\in A_U,\\
y_i^*=u_i^*\in(0,U_{\max}), & \text{if }i\in F.
\end{cases}
\]
On the fixed grid, all strict inequalities have a positive margin. Choose $\delta>0$ smaller than every number in
\[
\{-y_i^*:i\in A_0\}\cup\{y_i^*-U_{\max}:i\in A_U\}
\cup\{u_i^*,U_{\max}-u_i^*:i\in F\}.
\]
Continuity gives a neighborhood of $u^*$ on which
\[
\left\|u-\frac\theta\lambda g(u)-y^*\right\|_\infty<\delta.
\]
Projection therefore fixes active coordinates at their bounds and leaves free coordinates unchanged throughout this neighborhood. Its derivative is $P_F$, giving
\begin{equation*}
D\Phi_\theta(u^*)=P_F\left(I-\tfrac{\theta}{\lambda}H(u^*)\right).
\end{equation*}
The map is locally constant when $F=\varnothing$. Otherwise, write $A=\{i_1^A<\cdots<i_{|A|}^A\}$ and define
\[
 \mapdef{E_A}{\R^{|A|}}{\R^{N(M+1)}}{x}{\sum_{k=1}^{|A|} x_ke_{i_k^A}}.
\]
Ordering coordinates as $(F,A)$, with $E_F^TH(u^*)E_A\in\R^{|F|\times |A|}$,
\begin{equation*}
D\Phi_\theta(u^*)\simeq
\begin{pmatrix}
I_F-\tfrac{\theta}{\lambda}H_F&-\tfrac{\theta}{\lambda}E_F^TH(u^*)E_A\\
0&0
\end{pmatrix}.
\end{equation*}
Starting close enough to $u^*$, one sweep identifies the active coordinates. On $u_A=u_A^*$, write $u=u^*+E_Fx$ and choose an open neighborhood $\mathcal O_F\subset\R^{|F|}$ of zero where $u^*+E_Fx\in\mathcal O_{h,\tau}$ and projection leaves free coordinates unchanged. The free coordinates then follow
\[
 \mapdef{\Psi_\theta}{\mathcal O_F}{\R^{|F|}}{x}{x-\tfrac{\theta}{\lambda}E_F^Tg(u^*+E_Fx)},
\]
with $\Psi_\theta(0)=0$ and $D\Psi_\theta(0)=I_F-\tfrac{\theta}{\lambda}H_F$.
Since $H_F^TW_F=W_FH_F$ and $W_F$ is positive diagonal, there is a $W_F$-orthonormal eigenbasis $(q_\ell)_{\ell=1}^{|F|}\subset\R^{|F|}$ with
\begin{equation*}
H_Fq_\ell=\mu_\ell q_\ell,\qquad q_\ell^TW_Fq_m=\delta_{\ell m},
\end{equation*}
and hence
\begin{equation*}
\left(I_F-\tfrac{\theta}{\lambda}H_F\right)q_\ell=\left(1-\tfrac{\theta}{\lambda}\mu_\ell\right)q_\ell,
\end{equation*}
where $\delta_{\ell m}$ is the Kronecker delta, equal to one if $\ell=m$ and zero otherwise. For $e=\sum_{\ell=1}^{|F|} c_\ell q_\ell\in\R^{|F|}$, $c_\ell\in\R$,
\begin{equation*}
\left\|\left(I_F-\tfrac{\theta}{\lambda}H_F\right)e\right\|_{W_F}^2
=\sum_{\ell=1}^{|F|}\left|1-\tfrac{\theta}{\lambda}\mu_\ell\right|^2|c_\ell|^2.
\end{equation*}
Thus the free operator norm is \eqref{eq:rho}. The condition $\rho_\theta<1$ is equivalent to requiring that, for every $\mu\in\sigma(H_F)$, we have $-1<1-\theta\mu/\lambda<1$. Here $\sigma(H_F)$ is the spectrum. Since $\mu_{\min}>0$, this gives
\[
\rho_\theta<1 \quad\iff\quad 0<\theta<2\lambda/\mu_{\max}.
\] For $\rho_\theta<\rho'<1$, continuity gives $\varepsilon>0$ such that, on the ball $\|x\|_{W_F}<\varepsilon$, we have $\|D\Psi_\theta(x)\|_{W_F}\le\rho'$.
For such $x$,
\begin{align*}
\|\Psi_\theta(x)\|_{W_F}
=\left\|\int_0^1D\Psi_\theta(tx)x\,\mathrm{d}t\right\|_{W_F}
\le\rho'\|x\|_{W_F}.
\end{align*}
The update keeps this ball invariant. If the starting control is close enough to $u^*$, continuity puts the first iterate in the ball, and all later iterates satisfy
\begin{equation*}
\normt{u_{k+1}-u^*}\le\rho'\normt{u_k-u^*},
\end{equation*}
which is local $q$-linear convergence.
For $\theta>2\lambda/\mu_{\max}$, the mode with the largest eigenvalue has factor $1-\theta\mu_{\max}/\lambda<-1$. At equality, the factor is $-1$ and higher-order terms determine stability. The smallest spectral factor comes from balancing the factors at the two extreme eigenvalues:
\begin{equation*}
1-\frac{\theta\mu_{\min}}\lambda
=-1+\frac{\theta\mu_{\max}}\lambda
\quad\Longrightarrow\quad
\theta_R=\frac{2\lambda}{\mu_{\min}+\mu_{\max}}.
\end{equation*}
Substitution gives $\rho_{\theta_R}=(\mu_{\max}-\mu_{\min})/(\mu_{\max}+\mu_{\min})$. Restricting $\theta\le1$ gives $\min\{1,\theta_R\}$.
\subsection{Ritz residual and spectral distance}\label{app:proof-ritz-residual}
\begin{proof}
Choose an orthonormal eigenbasis $(e_\ell)_{\ell=1}^{|F|}$ of $\widetilde H_F$, with eigenvalues $\mu_\ell$. For $y=\sum_{\ell=1}^{|F|} a_\ell e_\ell$, $\sum_{\ell=1}^{|F|} a_\ell^2=1$,
\begin{align*}
\|\widetilde H_Fy-\widehat\mu y\|_2^2
=\sum_{\ell=1}^{|F|} a_\ell^2(\mu_\ell-\widehat\mu)^2
\ge\min_{\mu\in\sigma(\widetilde H_F)}|\widehat\mu-\mu|^2.
\end{align*}
Thus
\[
\min_{\mu\in\sigma(\widetilde H_F)}|\widehat\mu-\mu|
\le\|\widetilde H_Fy-\widehat\mu y\|_2.
\]
\end{proof}
\section{Harvest-quadrature proofs}\label{app:harvest-quadrature-proofs}
\subsection{Removal identities}\label{app:proof-removal-identities}
\begin{proof}
The removal step gives
\begin{equation*}
y^n-y^{n+1}\byeq{eq:pure-harvest-new}{=}\tau q u^n y^{n+1},\qquad
y^k=y^0\prod_{n=0}^{k-1}(1+\tau q u^n)^{-1}.
\end{equation*}
Summing over time gives
\begin{equation*}
\tau q\sum_{n=0}^{N-1} u^ny^{n+1}
 =\sum_{n=0}^{N-1}(y^n-y^{n+1})=y^0-y^N,
\end{equation*}
and
\begin{equation*}
j_1(u)-j_0(u)
 =\eta\tau q\sum_{n=0}^{N-1} u^n(y^n-y^{n+1})
 =\eta\tau^2q^2\sum_{n=0}^{N-1}(u^n)^2y^{n+1}.
\end{equation*}
These are \eqref{eq:pure-telescope-new} and \eqref{eq:pure-defect-new}.
\end{proof}
\subsection{Proof of Proposition~\ref{prop:pure-curvature-new}}\label{app:proof-prop-pure-curvature-new}
\begin{proof}
Differentiating the product formula for $y^N$ gives, for $0\le m,n<N$,
\begin{equation*}
\partial_{u^n}y^N=-\frac{\tau qy^N}{1+\tau qu^n},\qquad
\partial_{u^m}\frac1{1+\tau qu^n}=-\frac{\tau q\delta_{mn}}{(1+\tau qu^n)^2}.
\end{equation*}
Here $\delta_{mn}$ is the Kronecker delta. Applying the product rule once more gives
\begin{align*}
\partial_{u^m}\partial_{u^n}y^N
&=-\tau q\left(\frac{\partial_{u^m}y^N}{1+\tau qu^n}
-\frac{\tau q\delta_{mn}y^N}{(1+\tau qu^n)^2}\right)=\frac{\tau^2q^2y^N(1+\delta_{mn})}{(1+\tau qu^m)(1+\tau qu^n)}.
\end{align*}
The weighted gradient is $g_1(u)_n=\tau^{-1}\partial_{u^n}j_1(u)$. Thus \eqref{eq:pure-telescope-new} gives
\begin{equation*}
g_1(u)_n=\lambda u^n-\frac{\eta qy^N}{1+\tau qu^n},\qquad
(H_1(u))_{mn}=\lambda\delta_{mn}
+\frac{\eta\tau q^2y^N(1+\delta_{mn})}{(1+\tau qu^m)(1+\tau qu^n)}.
\end{equation*}
For $\xi\in\R^N$, the two terms in $1+\delta_{mn}$ yield
\begin{align*}
\langle H_1(u)\xi,\xi\rangle_\tau
&=\lambda\|\xi\|_\tau^2+\eta\tau^2q^2y^N\left[
\left(\sum_{n=0}^{N-1}\frac{\xi^n}{1+\tau qu^n}\right)^2
+\sum_{n=0}^{N-1}\left(\frac{\xi^n}{1+\tau qu^n}\right)^2\right]\ge\lambda\|\xi\|_\tau^2.
\end{align*}
Permuting the entries leaves the unique minimizer unchanged: $u=s\mathbf1$. Along this line,
\[
 \frac{\mathrm{d}}{\mathrm{d}s}j_1(s\mathbf1)
 =N\tau\bigl[\lambda s-\eta qy^0(1+\tau q s)^{-N-1}\bigr].
\]
The bracketed expression is strictly increasing. Its nonnegative root determines $s$, capped at $U_{\max}$. When $\eta=0$, this gives $s=0$.
To compute the old-state Hessian at zero, use that $y^n$ depends only on $u^0,\ldots,u^{n-1}$:
\[
\partial_{u^m}y^n(0)=
\begin{cases}
-\tau qy^0,&m<n,\\
0,&m\ge n.
\end{cases}
\]
Twice differentiating $-\eta\tau q\sum_{n=0}^{N-1}u^ny^n$ and evaluating at zero gives, for $m\ne\ell$,
\[
\partial_{u^m}\partial_{u^\ell}j_0(0)
=-\eta\tau q\bigl(\partial_{u^\ell}y^m(0)
+\partial_{u^m}y^\ell(0)\bigr)
=\eta\tau^2q^2y^0.
\]
For $m=\ell$, the reward contributes zero and the effort contributes $\lambda\tau$. Dividing by the time weight $\tau$ yields
\[
H_0(0)=\lambda I+\eta\tau q^2y^0(\one\one^T-I),\qquad
H_1(0)=\lambda I+\eta\tau q^2y^0(\one\one^T+I).
\]
Thus, at zero,
\[
 H_1(0)-H_0(0)\byeq{eq:pure-defect-new}{=}2\eta\tau q^2y^0I.
\]
For $\xi\ne0$ with $\one^T\xi=0$, multiplication by the rank-one matrix gives $\one\one^T\xi=0$. Hence
\[
\frac{\langle H_0(0)\xi,\xi\rangle_\tau}{\|\xi\|_\tau^2}
=\lambda-\eta\tau q^2y^0,\qquad
\frac{\langle H_1(0)\xi,\xi\rangle_\tau}{\|\xi\|_\tau^2}
=\lambda+\eta\tau q^2y^0.
\] Negative curvature persists by continuity at $u=s\mathbf1$ for small $0<s<U_{\max}$. Both $u\pm\varepsilon\xi$ are then feasible for small $\varepsilon>0$, so $j_0$ is nonconvex on $[0,U_{\max}]^N$.
\end{proof}
\subsection{Refinement bound in the pure-removal model}\label{app:proof-reduced-refinement}
\begin{proof}
For $\xi,\chi\in\R^N$, set $z[\xi]^k\coloneqq Dy^k(u)[\xi]$,
$s[\xi,\chi]^k\coloneqq D^2y^k(u)[\xi,\chi]$. Differentiating the product formula for $y^k$ gives
\begin{equation*}
z[\xi]^k=-\tau qy^k\sum_{n=0}^{k-1}\frac{\xi^n}{1+\tau qu^n},
\end{equation*}
and
\begin{equation*}
s[\xi,\chi]^k=\tau^2q^2y^k\left[
\left(\sum_{n=0}^{k-1}\frac{\xi^n}{1+\tau qu^n}\right)
\left(\sum_{n=0}^{k-1}\frac{\chi^n}{1+\tau qu^n}\right)
+\sum_{n=0}^{k-1}\frac{\xi^n\chi^n}{(1+\tau qu^n)^2}\right].
\end{equation*}
Since $u^n\ge0$ and $k\tau\le T$, Cauchy--Schwarz gives
\begin{align*}
\tau\sum_{n=0}^{k-1}\frac{|\xi^n|}{1+\tau qu^n}
&\le\sqrt{k\tau}\left(\tau\sum_{n=0}^{k-1}|\xi^n|^2\right)^{1/2}
\le\sqrt T\|\xi\|_\tau,\qquad
\tau^2\sum_{n=0}^{k-1}\frac{|\xi^n\chi^n|}{(1+\tau qu^n)^2}\le\tau\|\xi\|_\tau\|\chi\|_\tau.
\end{align*}
Combining these estimates with $0<y^k\le y^0$ gives
\begin{equation}\label{eq:pure-first-sensitivity-bound-restyled}
\max_{0\le k\le N}|z[\xi]^k|\le qy^0\sqrt T\|\xi\|_\tau,
\end{equation}
and
\begin{equation}\label{eq:pure-sensitivity-bounds-restyled}
\max_{0\le k\le N}|s[\xi,\chi]^k|
\le q^2y^0(T+\tau)\|\xi\|_\tau\|\chi\|_\tau.
\end{equation}
Differentiating each factor in the defect twice gives
\begin{align*}
D^2(j_1-j_0)(u)[\xi,\chi]
&\byeq{eq:pure-defect-new}{=}\eta\tau^2q^2\sum_{n=0}^{N-1}
\bigl[2\xi^n\chi^ny^{n+1}+2u^n\xi^nz[\chi]^{n+1}+2u^n\chi^nz[\xi]^{n+1}+(u^n)^2s[\xi,\chi]^{n+1}\bigr].
\end{align*}
Use $u^n\le U_{\max}$, $\tau\sum_{n=0}^{N-1}|\xi^n|\le\sqrt T\|\xi\|_\tau$, and \eqref{eq:pure-first-sensitivity-bound-restyled}--\eqref{eq:pure-sensitivity-bounds-restyled}. The first term, the two cross terms and the last term give, respectively,
\begin{align*}
|D^2(j_1-j_0)(u)[\xi,\chi]|
&\le\eta\tau q^2y^0
\bigl[2+4qU_{\max}T+q^2U_{\max}^2T(T+\tau)\bigr]
\|\xi\|_\tau\|\chi\|_\tau.
\end{align*}
Taking the supremum over unit directions and using $\tau\le T$ gives
\begin{equation*}
\|H_1(u)-H_0(u)\|_{\mathcal L(\R^N)}
\le\eta\tau q^2y^0\bigl(2+4qU_{\max}T+2q^2U_{\max}^2T^2\bigr).
\end{equation*}
For fixed parameters, the bound falls below $\lambda$ on all sufficiently fine time grids, making $j_0$ strictly convex.
\end{proof}
\subsection{Proof of Proposition~\ref{prop:quadrature-induced-local-minimum}}\label{app:proof-prop-quadrature-induced-local-minimum}
\begin{proof}
Multiplying $j_0$ by $\tau q^2/\lambda$ gives
\[
 F_c(x,z)=\tfrac12(x^2+z^2)-c\left(x+\frac{z}{1+x}\right).
\]
At $(0,c)$, $\partial_xF_c=c(c-1)>0$. This derivative remains positive on a small rectangle about $(0,c)$, so $F_c(x,z)>F_c(0,z)$ for every feasible point there with $x>0$. On the edge $x=0$,
\[
 F_c(0,z)=F_c(0,c)+\tfrac12(z-c)^2.
\]
Thus $(0,c)$ is a strict constrained local minimizer, including when $c=\tau q U_{\max}$. The point $(c,c/(1+c))$ is feasible and gives
\[
 F_c\left(c,\frac{c}{1+c}\right)-F_c(0,c)=-\frac{c^2}{2(1+c)^2}.
\]
Proposition~\ref{prop:pure-curvature-new} gives the flux assertion.
\end{proof}
\subsection{Proof of Proposition~\ref{prop:budget}}\label{app:proof-prop-budget}
\begin{proof}
Rewrite the old-state reward using
\begin{equation*}
a+\frac{E-a}{1+\tau qa}=E-\frac{\tau qa(E-a)}{1+\tau qa}.
\end{equation*}
The quotient rule gives
\begin{equation*}
\frac{\mathrm{d}}{\mathrm{d}a}\frac{a(E-a)}{1+\tau qa}
=\frac{(E-2a)(1+\tau qa)-\tau qa(E-a)}{(1+\tau qa)^2}
=\frac{E-2a-\tau qa^2}{(1+\tau qa)^2}.
\end{equation*}
Differentiating again, and using $a\le E$ and the assumed bound on $\lambda$, yields
\begin{align*}
f_0''(a)
&=2\lambda\tau-\frac{2\eta\tau^2q^2y^0(1+\tau qE)}{(1+\tau qa)^3}\le2\lambda\tau-\frac{2\eta\tau^2q^2y^0}{(1+\tau qE)^2}<0.
\end{align*}
For the flux reward, expand the denominator as
\begin{equation*}
(1+\tau qa)(1+\tau q(E-a))=1+\tau qE+\tau^2q^2a(E-a).
\end{equation*}
Its first and second derivatives are $\tau^2q^2(E-2a)$ and $-2\tau^2q^2$. The reciprocal rule therefore gives
\begin{align*}
f_1''(a)
&=2\lambda\tau+\frac{2\eta\tau^2q^2y^0}{[1+\tau qE+\tau^2q^2a(E-a)]^2}+\frac{2\eta\tau^4q^4y^0(E-2a)^2}{[1+\tau qE+\tau^2q^2a(E-a)]^3}>0.
\end{align*}
Since $f_1(a)=f_1(E-a)$, strict convexity puts its minimum at $a=E/2$. The strictly concave function $f_0$ has equal minima at the two endpoints.
For the continuous control $u_a(t)\coloneqq a$ on $[0,\tau)$,
$u_a(t)\coloneqq E-a$ on $[\tau,2\tau]$,
\begin{equation*}
\int_0^{2\tau}u_a\,\mathrm{d}t=\tau E,\qquad
j_{\rm c}(u_a)
 =\frac{\lambda\tau}{2}[a^2+(E-a)^2]
 -\eta y^0(1-e^{-\tau qE}).
\end{equation*}
The effort term alone depends on $a$, with its minimum at $E/2$.
\end{proof}
\subsection{Harvest quadrature in the coupled system}\label{app:proof-coupled-sign}
\begin{proof}
In \eqref{eq:coupled-defect-new}, the harvest and predation terms are
nonnegative, the growth term is nonpositive, and
\begin{align*}
 -\mu_2\langle u^n,L_hv_2^{n+1}\rangle_h
 &=\frac{\mu_2}{h}\sum_{j=0}^{M-1}
 (u_{j+1}^n-u_j^n)(v_{2,j+1}^{n+1}-v_{2,j}^{n+1})
\end{align*}
can have either sign, as can the total defect. To see this, take one step with constant data
$v_1^0=a\one$, $v_2^0=b\one$, $u^0=c\one$, where
$a\ge0$, $b>0$, $0<c\le U_{\max}$. Since $L_h\one=0$,
$\langle\one,\one\rangle_h=1$, the prey step gives
\begin{equation*}
v_2^1=\frac{(1+\tau r)b}{1+\tau(\beta a+qc)}\one,
\end{equation*}
so
\begin{align*}
D(u)=\eta q\tau c\left(
 b-\frac{(1+\tau r)b}{1+\tau(\beta a+qc)}\right)
 =\frac{\eta q\tau^2cb(\beta a+qc-r)}{1+\tau(\beta a+qc)}.
\end{align*}
For $\eta>0$, choose $0<c<\min\{U_{\max},r/q\}$. Then $D(u)<0$ when
$a=0$, while $D(u)>0$ when $a>r/\beta$. For $\eta=0$, $D=0$.
\end{proof}
\begin{proof}[Proof of Proposition~\ref{prop:coupled-local-minimum}]
Write $x=u^0$, $z=u^1$ and $y(x)=v_2^1$. Since $\tau=1$, the first prey step and the old-state cost are
\begin{align*}
 \bigl(\tfrac32 I-\mu_2L_h+\diag(x)\bigr)y(x)&=\tfrac32\one,\\
 j_0(x,z)&=\tfrac12\bigl(\normh{x}^2+\normh{z}^2\bigr)
                -2\ip{x}{\one}-2\ip{z}{y(x)}.
\end{align*}
Set $A=\tfrac32 I-\mu_2 L_h$. It is self-adjoint in $\langle\cdot,\cdot\rangle_h$ and $A^{-1}\one=\tfrac23\one$. Thus $y(0)=\one$ and $Dy(0)[d]=-A^{-1}d$. For any first-step direction $d$,
\[
 D_x j_0(0,2\one)[d]=-2\ip{\one}{d}+4\ip{\one}{A^{-1}d}
                       =\tfrac23\ip{\one}{d}.
\]
Every component of the first-step gradient is positive on a sufficiently small neighborhood of $(0,2\one)$, since the grid is finite and its quadrature weights are positive. Integrating along the segment from $(0,z)$ to $(x,z)$ gives $j_0(x,z)>j_0(0,z)$ for every nearby feasible $x\ne0$. On the face $x=0$,
\[
 j_0(0,z)=-2+\tfrac12\normh{z-2\one}^2.
\]
This proves strict local minimality in the full box. At $x=2\one$, $y(x)=\tfrac37\one$, giving $j_0(2\one,\one)=-33/14$.

To test the flux cost, keep $x=0$ and set $z=t\one$. The first step leaves the initial populations unchanged, so the second prey step gives $v_2^2=3\one/(3+2t)$. Hence
\[
 j_1(0,t\one)=\tfrac12t^2-\frac{6t}{3+2t},\qquad
 \left.\frac{\mathrm d}{\mathrm dt}j_1(0,t\one)\right|_{t=2}=\frac{80}{49}.
\]
The direction $(0,-\one)$ is feasible at $(0,2\one)$ and is a descent direction for $j_1$.
\end{proof}
\subsection{Exact global certificate for the two-step problem}\label{app:proof-global-certificate}
\begin{proof}
For $F_0$ in Section~\ref{sec:reduced-certificate}, $\partial_zF_0=z-2/(1+x)$ and $\partial_{zz}F_0=1$. Thus, at each $x\in[0,3]$, the unique minimizing value is $z=2/(1+x)\in[1/2,2]$. Substitution reduces the problem to
\[
 \mapdef{\varphi_0}{[0,3]}{\R}{x}{
 F_0\left(x,\frac2{1+x}\right)
 =\tfrac12x^2-2x-\frac2{(1+x)^2}}.
\]
Set
\[
 \mapdef{P}{\R}{\R}{x}{x^4+x^3-3x^2-5x+2}.
\]
Then
\[
 \varphi_0'(x)=x-2+\frac4{(1+x)^3}=\frac{P(x)}{(1+x)^3}.
\]
The rational Sturm sequence $P_0\coloneqq P$, $P_1\coloneqq P'$,
$P_{k+1}\coloneqq-\operatorname{rem}(P_{k-1},P_k)$, where $\operatorname{rem}(A,B)$ is the remainder of polynomial division of $A$ by $B$, is
\begin{equation*}
P_1(x)=4x^3+3x^2-6x-5,\qquad P_2(x)=(27x^2+54x-37)/16,\qquad
P_3(x)=(320-256x)/27,\qquad P_4(x)=-1163/256.
\end{equation*}
Define
\[
 \mapdef{V}{\R}{\{0,1,2,3,4\}}{a}{
 \text{sign changes in }(P_0(a),\ldots,P_4(a))},
\]
with zeros omitted. We have $V(0)=3$ and $V(3)=1$, so there are exactly two interior stationary points. Rational root isolation places one in each interval
\begin{equation*}
I_1\coloneqq[0.3408975428038,\,0.3408975428039],\qquad
I_2\coloneqq[1.8220156900313,\,1.8220156900314].
\end{equation*}
For $0\le a\le x\le b$, monotonicity of each term gives
\[
 \frac{a^2}{2}-2b-\frac{2}{(1+a)^2}
 \le\varphi_0(x)\le
 \frac{b^2}{2}-2a-\frac{2}{(1+b)^2}.
\]
On $I_1,I_2$, these bounds give
\begin{equation*}
\varphi_0(I_1)\subset[-1.7360327224,\,-1.7360327223],\qquad
\varphi_0(I_2)\subset[-2.2352980504,\,-2.2352980503].
\end{equation*}
All interval endpoints are exact rationals. The Sturm variations at the isolating endpoints satisfy
\[
V(0.3408975428038)=3,\quad V(0.3408975428039)=2,
\qquad V(1.8220156900313)=2,\quad V(1.8220156900314)=1.
\]
Each interval therefore contains exactly one root, and $V(0)-V(3)=2$ excludes any other interior critical point. Comparing these bounds with $\varphi_0(0)=-2$ and $\varphi_0(3)=-13/8$ proves that the stationary point in $I_2$ is the unique global minimizer. Its control is approximately $(1.822,0.7087)$, with $F_0\approx-2.235$. Appendix~\ref{app:proof-prop-quadrature-induced-local-minimum} proves that $(0,2)$ is a strict local minimum.
For $F_1$, strict convexity and symmetry give $x=z=t$. Then
\[
 t-\frac2{(1+t)^3}=0
 \quad\Longleftrightarrow\quad t(1+t)^3=2.
\]
Rational root isolation gives \[
t\in[0.5436890126920,\,0.5436890126921],
\]
so $t\approx0.5437$ and $F_1\approx-0.8651$.
For the continuous cost $j_c$ on $K_c$ from Section~\ref{sec:reduced-certificate}, take $u\in K_c$ and set $s\coloneqq\frac12\int_0^2u$. Jensen's inequality gives
$\frac12\int_0^2u^2\ge s^2$, with equality only for $u\equiv s$.
Define
\[
 \mapdef{f_c}{[0,3]}{\R}{s}{s^2-2+2e^{-2s}}.
\]
Then $f_c'(s)=2s-4e^{-2s}, f_c''(s)=2+8e^{-2s}>0.$ The unique minimizer is $s\approx0.6011$, with cost approximately $-1.038$. For two-step controls, $j_c(x,z)=\frac12(x^2+z^2)-2+2e^{-(x+z)}$. The old-state and flux global controls give continuous costs of about $0.07021$ and $-1.030$, respectively.
\end{proof}
\section{Uniform derivative estimates and stability proofs}\label{app:uniform-derivative-proofs}
\subsection{Proof of Theorem~\ref{th:quadrature-C2-new}}\label{app:proof-th-quadrature-C2-new}
\begin{proof}
On a fixed grid, the state map $S^+$ is smooth near $\Uk$. For $a,b\in\R^{N(M+1)}$, write its first and second variations as
\[
 z[a]\coloneqq DS^+(u)[a],\qquad s[a,b]\coloneqq D^2S^+(u)[a,b].
\]
Both variations lie in $\R^{2(N+1)(M+1)}$ and start at zero. With spatial products taken componentwise, the tangent equations are
\begin{align*}
 D_\tau z_2^n-\mu_2L_hz_2^{n+1}
 &={}rz_2^n-(\beta v_1^n+qu^n)z_2^{n+1}
 -\beta v_2^{n+1}z_1^n-qv_2^{n+1}a^n,\\
 D_\tau z_1^n-\mu_1L_hz_1^{n+1}
 &={}-mz_1^{n+1}+\alpha v_2^{n+1}z_1^n
 +\alpha v_1^nz_2^{n+1}.
\end{align*}
Write the common linear part as
\[
 D_\tau w^n-\operatorname{blockdiag}(\mu_1L_h,\mu_2L_h)w^{n+1}
 =P_nw^n+Q_nw^{n+1}+F^n,
\]
where $w^n,F^n\in\R^{2(M+1)}$. In the component order $(w_1,w_2)$, the block matrices are
\[
P_n=\begin{pmatrix}
\alpha\diag(v_2^{n+1})&0\\
-\beta\diag(v_2^{n+1})&rI
\end{pmatrix},\qquad
Q_n=\begin{pmatrix}
-mI&\alpha\diag(v_1^n)\\
0&-\diag(\beta v_1^n+qu^n)
\end{pmatrix}.
\]
Each diagonal block acts by componentwise multiplication. Since $0\le v_i^n\le K_i\mathbf1$ and $0\le u^n\le U_{\max}\mathbf1$, the weighted operator norms satisfy
\[
\|P_n\|\le(\alpha+\beta)K_2+r\le A,\qquad
\|Q_n\|\le m+(\alpha+\beta)K_1+qU_{\max}\le A,
\]
with $A$ from Remark~\ref{rem:explicit-constants-new}. Set $E_n\coloneqq\|w_1^n\|_h^2+\|w_2^n\|_h^2$. Test with $2\tau w^{n+1}$ and use $2\langle x-y,x\rangle=\|x\|^2-\|y\|^2+\|x-y\|^2$. 

The energy identity is
\begin{equation*}
E_{n+1}-E_n+\|w^{n+1}-w^n\|_h^2
+2\tau\sum_{i=1}^{2}\mu_i|w_i^{n+1}|_{1,h}^2 =\underbrace{2\tau\langle P_nw^n,w^{n+1}\rangle_h}
_{\substack{
\le 2\tau A\sqrt{E_nE_{n+1}}\\
\le \tau A(E_n+E_{n+1})
}}
+\underbrace{2\tau\langle Q_nw^{n+1},w^{n+1}\rangle_h}
_{\le 2\tau AE_{n+1}}
+\underbrace{2\tau\langle F^n,w^{n+1}\rangle_h}
_{\le \tau(\|F^n\|_h^2+E_{n+1})}.
\end{equation*}

Dropping the nonnegative terms on the left gives
\begin{align*}
 E_{n+1}-E_n
 &\le\tau\bigl[A E_n+(3A+1)E_{n+1}+\|F^n\|_h^2\bigr].
\end{align*}
Rearranging this inequality gives
\[
[1-(3A+1)\tau]E_{n+1}\le(1+A\tau)E_n+\tau\|F^n\|_h^2.
\]
For $\tau\le[2(3A+1)]^{-1}$, the denominator is at least $1/2$, and
\[
\frac{1+A\tau}{1-(3A+1)\tau}
=1+\frac{(4A+1)\tau}{1-(3A+1)\tau}
\le1+2(4A+1)\tau.
\]
Consequently,
\[
 E_{n+1}\le[1+2(4A+1)\tau]E_n+2\tau\|F^n\|_h^2.
\]
Iteration bounds $\max_{0\le n\le N}E_n$ uniformly in terms of $E_0$ and $\|F\|_{h,\tau}^2$, which also controls $P_nw^n+Q_nw^{n+1}+F^n$ in the space--time norm.
For a scalar heat step $D_\tau w^n-\mu L_hw^{n+1}=f^n$, test with
$2\tau D_\tau w^n$. Summation by parts gives
\[
2\tau\|D_\tau w^n\|_h^2
+\mu\bigl(|w^{n+1}|_{1,h}^2-|w^n|_{1,h}^2+|w^{n+1}-w^n|_{1,h}^2\bigr) =\underbrace{2\tau\langle f^n,D_\tau w^n\rangle_h}
_{\substack{
\le2\tau\|f^n\|_h\|D_\tau w^n\|_h\\
\le\tau\|f^n\|_h^2+\tau\|D_\tau w^n\|_h^2
}}.
\]
Thus, for $1\le k\le N$,
\begin{equation}\label{eq:heat-endpoint-restyled}
\begin{aligned}
 \tau\sum_{n=0}^{k-1}\|D_\tau w^n\|_h^2+\mu|w^k|_{1,h}^2
 +\mu\sum_{n=0}^{k-1}|w^{n+1}-w^n|_{1,h}^2
 \le\tau\sum_{n=0}^{k-1}\|f^n\|_h^2+\mu|w^0|_{1,h}^2.
\end{aligned}
\end{equation}
Combine the time-derivative bound at $k=N$ with the gradient bound at each $k$ to get
\begin{align*}
 \tau\sum_{n=0}^{N-1}\|D_\tau w^n\|_h^2+\mu\max_{0\le n\le N}|w^n|_{1,h}^2
 &\byeq{eq:heat-endpoint-restyled}{\le}
 2\left(\tau\sum_{n=0}^{N-1}\|f^n\|_h^2+\mu|w^0|_{1,h}^2\right).
\end{align*}
For any $f\in H^1(0,1)$ and $x\in[0,1]$, the continuous representative satisfies
\[
f(x)=\int_0^1f(y)\,\mathrm{d}y+\int_0^1(f(x)-f(y))\,\mathrm{d}y,
\qquad
|f(x)-f(y)|\le\int_{\min(x,y)}^{\max(x,y)}|f'(s)|\,\mathrm{d}s.
\]
Thus $|f(x)|\le\|f\|_{L^1}+\|f'\|_{L^1}\le\|f\|_{L^2}+\|f'\|_{L^2}$. Applying this to the nodal interpolant $\Ih w\in H^1(0,1)$ gives
\begin{align*}
 \|w\|_\infty
 &=\|\Ih w\|_{L^\infty}
 \le\|\Ih w\|_{L^2}+\|(\Ih w)'\|_{L^2}
 \le\|w\|_h+|w|_{1,h}.
\end{align*}
These mesh-independent estimates apply first to the bounded state reactions, then to the tangent forcing $F_2^n=-qv_2^{n+1}a^n$, $F_1^n=0$, giving
\begin{equation}\label{eq:state-MR-new}
 \|D_\tau v\|_{h,\tau}\le C,
\end{equation}
and
\begin{equation}\label{eq:tangent-MR-new}
 \|D_\tau z[a]\|_{h,\tau}+\max_{0\le n\le N}\|z[a]^n\|_\infty
 \le C\|a\|_{h,\tau}.
\end{equation}
The second variation has the same linear part, with forcing
\begin{equation*}
F_2^n
 \coloneqq{}-(\beta z_1[a]^n+qa^n)z_2[b]^{n+1}
 -(\beta z_1[b]^n+qb^n)z_2[a]^{n+1},
\qquad
F_1^n
 \coloneqq{}\alpha\bigl(z_1[a]^nz_2[b]^{n+1}
 +z_1[b]^nz_2[a]^{n+1}\bigr).
\end{equation*}
The product inequality $\|xy\|_h\le\|x\|_\infty\|y\|_h$ gives, for example,
\begin{align*}
 \|a\circ z_2[b]^+\|_{h,\tau}
 &\le\max_{0\le n<N}\|z_2[b]^{n+1}\|_\infty\|a\|_{h,\tau}
 \byeq{eq:tangent-MR-new}{\le}C\|a\|_{h,\tau}\|b\|_{h,\tau}.
\end{align*}
Applying this inequality to each product gives $\|F\|_{h,\tau}\le C\|a\|_{h,\tau}\|b\|_{h,\tau}$, so
\begin{equation}\label{eq:second-MR-new}
 \|D_\tau s[a,b]\|_{h,\tau}
 \le C\|a\|_{h,\tau}\|b\|_{h,\tau}.
\end{equation}
Only the stated data enter the constants. The directions enter through their weighted $L^2$ norms.
Differentiating \eqref{eq:defect-exact-new} gives
\begin{equation*}
D(j_1-j_0)(u)[a]
 =-\eta q\tau^2\sum_{n=0}^{N-1}\langle a^n,D_\tau v_2^n\rangle_h
 -\eta q\tau^2\sum_{n=0}^{N-1}\langle u^n,D_\tau z_2[a]^n\rangle_h,
\end{equation*}
and
\[
D^2(j_1-j_0)(u)[a,b] =-\eta q\tau^2\sum_{n=0}^{N-1}\langle a^n,D_\tau z_2[b]^n\rangle_h -\eta q\tau^2\sum_{n=0}^{N-1}\langle b^n,D_\tau z_2[a]^n\rangle_h -\eta q\tau^2\sum_{n=0}^{N-1}\langle u^n,D_\tau s_2[a,b]^n\rangle_h.
\]
Since $\|u\|_{h,\tau}\le U_{\max}\sqrt T$, Cauchy--Schwarz gives
\begin{equation*}
|D(u)|
 \le\eta q\tau\|u\|_{h,\tau}\|D_\tau v_2\|_{h,\tau}
 \byeq{eq:state-MR-new}{\le}C\tau,
\end{equation*}
and
\begin{align*}
|D(j_1-j_0)(u)[a]|
 \le\eta q\tau\bigl(
 \|a\|_{h,\tau}\|D_\tau v_2\|_{h,\tau}
 +\|u\|_{h,\tau}\|D_\tau z_2[a]\|_{h,\tau}\bigr)
 \overset{\mathclap{\text{\eqref{eq:state-MR-new}, \eqref{eq:tangent-MR-new}}}}{\le}
 C\tau\|a\|_{h,\tau},
\end{align*}
and
\begin{align*}
|D^2(j_1-j_0)(u)[a,b]|
 \le\eta q\tau\bigl(
 \|a\|_{h,\tau}\|D_\tau z_2[b]\|_{h,\tau}
 +\|b\|_{h,\tau}\|D_\tau z_2[a]\|_{h,\tau}
 +\|u\|_{h,\tau}\|D_\tau s_2[a,b]\|_{h,\tau}\bigr)
 \overset{\mathclap{\text{\eqref{eq:tangent-MR-new}, \eqref{eq:second-MR-new}}}}{\le}
 C\tau\|a\|_{h,\tau}\|b\|_{h,\tau}.
\end{align*}
The bound on $D$ gives \eqref{eq:C0-new}, and Riesz representation turns the derivative bounds into \eqref{eq:C1-new}--\eqref{eq:C2-new}.
\end{proof}
\subsection{Explicit constants}\label{app:proof-explicit-constants}
\begin{remark}[Explicit constants]
\label{rem:explicit-constants-new}
Let $K_i$ be the uniform state bounds. Set
\[
A\coloneqq m+r+(\alpha+\beta)(K_1+K_2)+qU_{\max},\qquad
R\coloneqq\sqrt2\exp((4A+1)T).
\]
Assume $\tau\le\bar\tau\coloneqq[2(3A+1)]^{-1}$. The first variation is controlled in the discrete $L^2$ and maximum norms by
\begin{equation*}
C_z\coloneqq\sqrt T RqK_2,\quad
C_{z,\infty}\coloneqq\left(R+\frac{1+2A\sqrt T R}{\sqrt{\min\{\mu_1,\mu_2\}}}\right)qK_2.
\end{equation*}
For the second variation in the discrete $L^2$ norm, use
\[
C_s\coloneqq2\sqrt T R C_{z,\infty}\bigl((\alpha+\beta)C_z+q\bigr).
\]
The proof below gives these sensitivity bounds and
\[
\|H_1(u)-H_0(u)\|_{\mathcal L(\R^{N(M+1)})}
\le\eta q(1+2A\sqrt T R)\left(2qK_2+\frac{U_{\max}C_s}{R}\right)\tau.
\]
\end{remark}
\begin{proof}
For the same variation system with $w^0=0$, iterate the energy recurrence. For $1\le k\le N$,
\begin{align*}
E_k
&\le2\tau\sum_{n=0}^{k-1}[1+2(4A+1)\tau]^{k-1-n}\|F^n\|_h^2
\le2e^{2(4A+1)T}\|F\|_{h,\tau}^2=R^2\|F\|_{h,\tau}^2.
\end{align*}
Hence $\max_{0\le n\le N}\|w^n\|_h\le R\|F\|_{h,\tau}$. The heat-equation forcing is $f^n=P_nw^n+Q_nw^{n+1}+F^n$. It satisfies
\begin{equation*}
\|f\|_{h,\tau}
\le A\left(
\underbrace{\|w^-\|_{h,\tau}}_{\le\sqrt T R\|F\|_{h,\tau}}
+\underbrace{\|w^+\|_{h,\tau}}_{\le\sqrt T R\|F\|_{h,\tau}}
\right)+\|F\|_{h,\tau} \le(1+2A\sqrt T R)\|F\|_{h,\tau}.
\end{equation*}
Applying \eqref{eq:heat-endpoint-restyled} to each component gives
\begin{equation*}
\|D_\tau w\|_{h,\tau}\le(1+2A\sqrt T R)\|F\|_{h,\tau},\qquad
\max_{0\le n\le N}|w^n|_{1,h}
\le\frac{1+2A\sqrt T R}{\sqrt{\min\{\mu_1,\mu_2\}}}\|F\|_{h,\tau}.
\end{equation*}
For pairs, $|w|_{1,h}^2=|w_1|_{1,h}^2+|w_2|_{1,h}^2$. The tangent forcing has norm at most $qK_2\|a\|_{h,\tau}$. Thus, for either time reconstruction,
\begin{equation*}
\|z_i[a]^\pm\|_{h,\tau}\le C_z\|a\|_{h,\tau},\qquad
\max_{0\le n\le N}\|z_i[a]^n\|_\infty\le C_{z,\infty}\|a\|_{h,\tau}.
\end{equation*}
For the second-variation forcing, bound each product by the spatial supremum norm of one factor and the space--time norm of the other. The two components satisfy
\[
\|F_1\|_{h,\tau}\le2\alpha C_{z,\infty}C_z\|a\|_{h,\tau}\|b\|_{h,\tau},\qquad
\|F_2\|_{h,\tau}\le2C_{z,\infty}(\beta C_z+q)\|a\|_{h,\tau}\|b\|_{h,\tau}.
\]
Since the norm of the pair is at most the sum of the component norms,
\begin{equation*}
\|F\|_{h,\tau}\le2C_{z,\infty}\bigl((\alpha+\beta)C_z+q\bigr)
\|a\|_{h,\tau}\|b\|_{h,\tau}.
\end{equation*}
The energy bound applied to the second variation then gives
\begin{equation*}
\|s_i[a,b]^\pm\|_{h,\tau}\le2\sqrt T R C_{z,\infty}\bigl((\alpha+\beta)C_z+q\bigr)\|a\|_{h,\tau}\|b\|_{h,\tau}
=C_s\|a\|_{h,\tau}\|b\|_{h,\tau}.
\end{equation*}
Applying the time-derivative estimate to the tangent and second-variation equations gives
\[
\|D_\tau z[a]\|_{h,\tau}
\le(1+2A\sqrt T R)qK_2\|a\|_{h,\tau},\qquad \|D_\tau s[a,b]\|_{h,\tau}
\le(1+2A\sqrt T R)\frac{C_s}{\sqrt T R}\|a\|_{h,\tau}\|b\|_{h,\tau}.
\]
Substitute these bounds into the formula for $D^2(j_1-j_0)$ in Appendix~\ref{app:proof-th-quadrature-C2-new} and use $\|u\|_{h,\tau}\le U_{\max}\sqrt T$:
\[
\begin{aligned}
|D^2(j_1-j_0)(u)[a,b]|
&\le\eta q\tau\bigl(\|a\|_{h,\tau}\|D_\tau z_2[b]\|_{h,\tau}
+\|b\|_{h,\tau}\|D_\tau z_2[a]\|_{h,\tau}
+\|u\|_{h,\tau}\|D_\tau s_2[a,b]\|_{h,\tau}\bigr)\\
&\le\eta q(1+2A\sqrt T R)\left(2qK_2+\frac{U_{\max}C_s}{R}\right)\tau\|a\|_{h,\tau}\|b\|_{h,\tau}.
\end{aligned}
\]
For a linear operator $L$ on the weighted control space,
\[
\|L\|_{h,\tau}
=\sup_{\|a\|_{h,\tau}=1}\|La\|_{h,\tau}
=\sup_{\|a\|_{h,\tau}=\|b\|_{h,\tau}=1}|\langle La,b\rangle_{h,\tau}|.
\]
With $L=H_1(u)-H_0(u)$, the bilinear form is $D^2(j_1-j_0)(u)[a,b]$. Taking the supremum proves the Hessian bound.
\end{proof}
\subsection{Proof of Corollary~\ref{cor:effort-dominance-new}}\label{app:proof-cor-effort-dominance-new}
\begin{proof}
For $a\in\R^{N(M+1)}$, use $\gamma_i\ge0$ to bound the tracking contribution:
\begin{equation*}
\underbrace{\gamma_i\|z_i[a]^-\|_{h,\tau}^2}_{\ge0}
+\underbrace{\gamma_i\langle v_i^--v_i^d,s_i[a,a]^-\rangle_{h,\tau}}
_{\ge-\gamma_i(K_i\sqrt T+D_i)C_s\|a\|_{h,\tau}^2} \ge-\gamma_i(K_i\sqrt T+D_i)C_s\|a\|_{h,\tau}^2.
\end{equation*}
For $\zeta\in\{0,1\}$, the harvest contribution satisfies
\begin{align*}
&\bigl|{}-2\eta q\langle a,(1-\zeta)z_2[a]^-+\zeta z_2[a]^+\rangle_{h,\tau}
-\eta q\langle u,(1-\zeta)s_2[a,a]^-+\zeta s_2[a,a]^+\rangle_{h,\tau}\bigr|\\
&\le\underbrace{2\eta q\bigl|\langle a,(1-\zeta)z_2[a]^-+\zeta z_2[a]^+\rangle_{h,\tau}\bigr|}
_{\le2\eta qC_z\|a\|_{h,\tau}^2}
+\underbrace{\eta q\bigl|\langle u,(1-\zeta)s_2[a,a]^-+\zeta s_2[a,a]^+\rangle_{h,\tau}\bigr|}
_{\le\eta qU_{\max}\sqrt T C_s\|a\|_{h,\tau}^2}
\le(2\eta qC_z+\eta qU_{\max}\sqrt T C_s)\|a\|_{h,\tau}^2.
\end{align*}
Adding the effort contribution gives
\[
 \langle H_\zeta(u)a,a\rangle_{h,\tau}
 \ge(\lambda-B)\|a\|_{h,\tau}^2.
\]
For $u,v\in\Uk$, integration along the feasible segment gives
\begin{equation*}
j_\zeta(v)-j_\zeta(u)-\langle g_\zeta(u),v-u\rangle_{h,\tau}
=\int_0^1(1-t)\langle H_\zeta(u+t(v-u))(v-u),v-u\rangle_{h,\tau}\,\mathrm{d}t \ge\frac{\lambda-B}{2}\|v-u\|_{h,\tau}^2.
\end{equation*}
A minimizer exists by compactness of $\Uk$. If $u$ and $v$ were two minimizers, their variational inequalities would give
\[
0\ge\langle g_\zeta(v)-g_\zeta(u),v-u\rangle_{h,\tau}
\ge(\lambda-B)\|v-u\|_{h,\tau}^2.
\]
Since $\lambda>B$, this forces $u=v$.
\end{proof}
\subsection{Proof of Corollary~\ref{cor:control-stability-new}}\label{app:proof-cor-control-stability-new}
\begin{proof}
The set $K_r$ is nonempty, convex and compact. For every $u\in K_r$ and $a\in\R^{N(M+1)}$,
\[
 \langle H_1(u)a,a\rangle_{h,\tau}
 \byeq{eq:C2-new}{\ge}(\kappa-C\tau)\|a\|_{h,\tau}^2.
\]
Because $\kappa-C\tau>0$, $j_1$ has a unique minimizer $u_1$ on $K_r$.
Set $d\coloneqq u_1-u_0$. The segment joining the two controls lies in $K_r$, so
\begin{align*}
 \langle g_0(u_1)-g_0(u_0),d\rangle_{h,\tau}
 &=\int_0^1\langle H_0(u_0+td)d,d\rangle_{h,\tau}\,\mathrm{d}t
 \ge\kappa\|d\|_{h,\tau}^2.
\end{align*}
Testing each variational inequality with the other control gives
$\langle g_0(u_0),d\rangle_{h,\tau}\ge0$,
$\langle g_1(u_1),d\rangle_{h,\tau}\le0$. Thus
\begin{align*}
\kappa\|d\|_{h,\tau}^2
&\le\langle g_0(u_1)-g_0(u_0),d\rangle_{h,\tau}
=\langle g_0(u_1)-g_1(u_1),d\rangle_{h,\tau}
+\underbrace{\langle g_1(u_1),d\rangle_{h,\tau}}_{\le0}
+\underbrace{-\langle g_0(u_0),d\rangle_{h,\tau}}_{\le0}\\
&\le\langle g_0(u_1)-g_1(u_1),d\rangle_{h,\tau}
\byeq{eq:C1-new}{\le}C\tau\|d\|_{h,\tau}.
\end{align*}
This proves \eqref{eq:O-tau-stability-new} and gives $\|d\|_{h,\tau}<r$. For any $u\in\Uk$, the point $u_1+t(u-u_1)$ stays in $K_r$ for small $t>0$. Testing along this segment and dividing by $t$ gives the full box inequality at $u_1$. Strong convexity then makes $u_1$ a strict local minimizer and the only box-stationary point in $K_r$. Apply the same argument on the whole box to prove the global statement.
\end{proof}
\subsection{Spectral perturbation at a common control}\label{app:proof-spectral-perturbation}
\begin{proof}
For any $a\in\R^{N(M+1)}$,
\[
 |\langle(H_1(u)-H_0(u))a,a\rangle_{h,\tau}|
 \byeq{eq:C2-new}{\le}C\tau\|a\|_{h,\tau}^2.
\]
For $a\ne0$, division by $\|a\|_{h,\tau}^2$ gives the Rayleigh-quotient bounds
\begin{equation*}
\frac{\langle H_0(u)a,a\rangle_{h,\tau}}{\|a\|_{h,\tau}^2}-C\tau
\le\frac{\langle H_1(u)a,a\rangle_{h,\tau}}{\|a\|_{h,\tau}^2}
\le\frac{\langle H_0(u)a,a\rangle_{h,\tau}}{\|a\|_{h,\tau}^2}+C\tau.
\end{equation*}
The matrices $\widetilde H_\zeta=\Wht^{1/2}H_\zeta \Wht^{-1/2}$ are Euclidean symmetric. Moreover, for $x=\Wht^{1/2}a$,
\[
\frac{x^T\widetilde H_\zeta x}{x^Tx}
=\frac{\langle H_\zeta a,a\rangle_{h,\tau}}{\|a\|_{h,\tau}^2}.
\]
For the eigenvalues ordered increasingly, the min--max formula is
\[
\lambda_k(H_\zeta)
=\min_{\substack{E\subset\R^{N(M+1)}\\\dim E=k}}
\ \max_{a\in E\setminus\{0\}}
\frac{\langle H_\zeta a,a\rangle_{h,\tau}}{\|a\|_{h,\tau}^2}.
\]
Taking the maximum and then the minimum in the two Rayleigh inequalities yields
\[
\lambda_k(H_0)-C\tau\le\lambda_k(H_1)\le\lambda_k(H_0)+C\tau.
\]
For a fixed subspace with an orthonormal basis matrix $E$, its Euclidean compressions satisfy
\[
\|E^T(\widetilde H_1-\widetilde H_0)E\|_2
\le\|\widetilde H_1-\widetilde H_0\|_2\le C\tau,
\]
so the same argument applies to the compressed operators.
\end{proof}
\section{Production--destruction proofs}\label{app:network-proofs}
\subsection{Weighted transfer cancellation}\label{app:proof-network-transfer}
\begin{proof}
The weighted column sum is
\begin{align*}
 \sum_{i=1}^s c_iw_aK_{(i,a),(j,a)}^n
 &=c_jw_a\sum_{\substack{1\le i\le s\\i\ne j}}k_{ij,a}^n
 -c_jw_a\sum_{\substack{1\le i\le s\\i\ne j}}k_{ij,a}^n=0.
\end{align*}
\end{proof}
\subsection{Proof of Theorem~\ref{th:pd-network}}\label{app:proof-th-pd-network}
\begin{proof}
Set $Q\coloneqq\operatorname{diag}(\rho)^{-1}(A^n)^T\operatorname{diag}(\rho)$.
Index the coordinates by $1\le\kappa,\nu\le sd$. For $\kappa\ne\nu$, we have $Q_{\kappa\nu}\le0$. For each row $\kappa$,
\[
 \sum_{\nu=1}^{sd} Q_{\kappa\nu}=1+\tau(m_\kappa^n+h_\kappa^n)>0.
\]
Thus $Q$ is a strictly row-diagonally dominant Z-matrix, meaning that its off-diagonal entries are nonpositive. If $Qx=f\ge0$
and $x_\kappa=\min_{1\le\nu\le sd}x_\nu<0$, then
\begin{align*}
 f_\kappa
 =Q_{\kappa\kappa}x_\kappa
 +\sum_{\substack{1\le\nu\le sd\\\nu\ne\kappa}}Q_{\kappa\nu}x_\nu\le\left(\sum_{\nu=1}^{sd} Q_{\kappa\nu}\right)x_\kappa<0.
\end{align*}
This contradicts the assumption. Applying the minimum and maximum arguments to $Qx=0$ gives $x=0$, so $Q$ is invertible with $Q^{-1}\ge0$. Transpose and undo the positive diagonal scaling to get $(A^n)^{-1}\ge0$: the step has a unique nonnegative solution.
Multiplying \eqref{eq:pd-step} by $\rho^T$ gives \eqref{eq:pd-mass}.
Since
\[
 \rho^TA^n=\rho^T+\tau\rho^T\operatorname{diag}(m^n+h^n),
\]
right multiplication by $(A^n)^{-1}\ge0$ gives
$\rho^T(A^n)^{-1}\le\rho^T$. Thus
\[
 \rho^T|(A^n)^{-1}z|
 \le\rho^T(A^n)^{-1}|z|\le\rho^T|z|.
\]
This proves \eqref{eq:pd-resolvent}. The production bound gives
\begin{equation*}
\mathcal M(y^{n+1})
 \byeq{eq:pd-mass}{\le}(1+\tau\sigma)\mathcal M(y^n)+\tau b_n,
\end{equation*}
and iteration yields
\begin{align*}
\mathcal M(y^n)
 \le(1+\tau\sigma)^n\mathcal M(y^0)
 +\tau\sum_{r=0}^{n-1}(1+\tau\sigma)^{n-1-r}b_r
 \le e^{\sigma t_n}\left(
 \mathcal M(y^0)+\tau\sum_{r=0}^{n-1}b_r\right).
\end{align*}
For $p^n=0$, summation gives
\begin{align*}
\underbrace{\mathcal M(y^N)}_{\ge0}
+\underbrace{\tau\sum_{n=0}^{N-1}\rho^T\operatorname{diag}(m^n)y^{n+1}}_{\ge0}
+\tau\sum_{n=0}^{N-1}\rho^T\operatorname{diag}(h^n)y^{n+1}
&\byeq{eq:pd-mass}{=}\mathcal M(y^0).
\end{align*}
Remove the nonnegative final mass and mortality terms to obtain the harvest bound.
\end{proof}
\subsection{Exact fixed-grid derivatives}\label{app:network-derivative-formulas}
Assume $G,p,m,h$ are $C^2$ on an open neighborhood $\mathcal N\subset\R^{sd}\times\R^{d_u}$ of the computed pairs of states and controls, with $p,m,h$ taking values in $\R^{sd}$ and
\[
 \mapdef{G}{\mathcal N}{\R^{sd\times sd}}{(y,u)}{
 \operatorname{blockdiag}_{i=1}^s(\mu_iL_h)+K(y,u)}.
\]
Define
\[
 \mapdef{F_n}{\{(y^-,y^+,u):(y^-,u)\in\mathcal N,\ y^+\in\R^{sd}\}}{\R^{sd}}{(y^-,y^+,u)}{
 A(y^-,u)y^+-y^--\tau p(y^-,u)}.
\]
At $y^-=y^n$, $y^+=y^{n+1}$, $u=u^n$, the matrix $A(y^n,u^n)=A^n$ is invertible, so the implicit function theorem gives a $C^2$ state map on an open neighborhood $\mathcal O_{\rm net}\subset(\R^{d_u})^N$ of the computed control:
\[
 \mapdef{S_{\rm net}}{\mathcal O_{\rm net}}{(\R^{sd})^{N+1}}{u}{
 (y^n)_{n=0}^N\text{ solving \eqref{eq:pd-step}}}.
\]
Write the residual derivatives as
\[
\partial_{y^+}F_n=A^n,\qquad
C_n\coloneqq\partial_{y^-}F_n\in\R^{sd\times sd},\qquad
B_n\coloneqq\partial_uF_n\in\R^{sd\times d_u}.
\]
For directions $d,e\in(\R^{d_u})^N$, the tangents $z[d],z[e]\in(\R^{sd})^{N+1}$ and second variation $w[d,e]$ satisfy
\[
A^nz[d]^{n+1}=-C_nz[d]^n-B_nd^n,\qquad z[d]^0=0,
\]
and
\[
A^nw[d,e]^{n+1}
=-C_nw[d,e]^n-D^2F_n[(z[d]^n,z[d]^{n+1},d^n),(z[e]^n,z[e]^{n+1},e^n)],
\qquad w[d,e]^0=0.
\]
Take a $C^2$ cost $J:\mathcal O_J\to\R$ on an open neighborhood $\mathcal O_J\subset(\R^{sd})^{N+1}\times(\R^{d_u})^N$ of the trajectory--control pair. Choose the neighborhoods small enough that $(y^n,u^n)\in\mathcal N$ for $n=0,\ldots,N-1$ and $(S_{\rm net}(u),u)\in\mathcal O_J$ for $u\in\mathcal O_{\rm net}$. Define
\[
 \mapdef{j}{\mathcal O_{\rm net}}{\R}{u}{J(S_{\rm net}(u),u)}.
\]
The Euclidean adjoints $a^{n+1}\in\R^{sd}$ satisfy
\begin{equation}\label{eq:pd-terminal}
 (A^{N-1})^Ta^N=-\partial_{y^N}J,
\end{equation}
and, for $1\le k<N$,
\[
 (A^{k-1})^Ta^k=-\partial_{y^k}J-C_k^Ta^{k+1}.
\]
The Euclidean reduced gradient blocks and Hessian bilinear form are
\begin{equation*}
\partial_{u^n}j=\partial_{u^n}J+B_n^Ta^{n+1},
\end{equation*}
and
\begin{equation*}
D^2j(u)[d,e]
=D^2J[(z[d],d),(z[e],e)]
+\sum_{n=0}^{N-1}(a^{n+1})^TD^2F_n[(z[d]^n,z[d]^{n+1},d^n),(z[e]^n,z[e]^{n+1},e^n)].
\end{equation*}
Appendix~\ref{app:proof-network-derivatives} derives these identities. If only the harvest reward depends on $y^N$, then
\[
 \partial_{y^N}J=-\zeta\eta\tau\operatorname{diag}(h^{N-1})\ell.
\]
Thus $a^N=0$ for the old-state rule, while the flux reward gives
\[
 (A^{N-1})^Ta^N=\eta\tau\operatorname{diag}(h^{N-1})\ell.
\]

\subsection{Residual derivatives and adjoint identities}\label{app:proof-network-derivatives}
\begin{proof}
The step matrix is
\[
A(y,u)=I-\tau G(y,u)+\tau\diag(m(y,u)+h(y,u)).
\]
For directions $\delta=(\xi,b)$ and $\varepsilon=(\chi,c)$ in the state--control variables, its derivatives are
\begin{align*}
DA(y,u)[\delta]
&=-\tau DG(y,u)[\delta]
+\tau\diag\bigl(Dm(y,u)[\delta]+Dh(y,u)[\delta]\bigr),\\
D^2A(y,u)[\delta,\varepsilon]
&=-\tau D^2G(y,u)[\delta,\varepsilon]
+\tau\diag\bigl(D^2m(y,u)[\delta,\varepsilon]+D^2h(y,u)[\delta,\varepsilon]\bigr).
\end{align*}
The diffusion block is constant, so only the transfer matrix $K(y,u)$ contributes to $DG$ and $D^2G$. Differentiating
$F_n=A(y^-,u)y^+-y^--\tau p(y^-,u)$ once gives
\[
DF_n[\xi^-,\xi^+,b]
=A^n\xi^++DA(y^n,u^n)[\xi^-,b]y^{n+1}
-\xi^--\tau Dp(y^n,u^n)[\xi^-,b].
\]
The residual derivatives with respect to the old state and control are
\[
 \mapdef{C_n}{\R^{sd}}{\R^{sd}}{\xi}{
 \begin{aligned}[t]
 &D_yA(y^n,u^n)[\xi]y^{n+1}-\xi-\tau D_yp(y^n,u^n)[\xi],
 \end{aligned}}
\]
and
\[
 \mapdef{B_n}{\R^{d_u}}{\R^{sd}}{b}{
 D_uA(y^n,u^n)[b]y^{n+1}-\tau D_up(y^n,u^n)[b]}.
\]
For directions $(\xi^-,\xi^+,b),(\chi^-,\chi^+,c)\in\R^{sd}\times\R^{sd}\times\R^{d_u}$, evaluating the derivatives of $A,p$ at $(y^n,u^n)$ gives
\[
\begin{aligned}
D^2F_n[(\xi^-,\xi^+,b),(\chi^-,\chi^+,c)]
&=D^2A[(\xi^-,b),(\chi^-,c)]y^{n+1}\\
&\quad+DA[(\xi^-,b)]\chi^++DA[(\chi^-,c)]\xi^+-\tau D^2p[(\xi^-,b),(\chi^-,c)].
\end{aligned}
\]
For the cost and adjoints of Appendix~\ref{app:network-derivative-formulas}, set
\[
 \mapdef{\mathcal L}{\mathcal O_J\times(\R^{sd})^N}{\R}{(y,u,a)}{
 J(y,u)+\sum_{n=0}^{N-1}(a^{n+1})^TF_n}.
\]
The terminal state appears only in $J$ and the final residual, so
\begin{equation*}
\partial_{y^N}\mathcal L=\partial_{y^N}J+(A^{N-1})^Ta^N.
\end{equation*}
For $1\le k<N$, the state $y^k$ appears in the two adjacent residuals. Hence
\begin{equation*}
\partial_{y^k}\mathcal L=\partial_{y^k}J+(A^{k-1})^Ta^k+C_k^Ta^{k+1}.
\end{equation*}
Setting these derivatives to zero gives \eqref{eq:pd-terminal} and the backward adjoint equation. For $0\le n<N$, the control derivative is
\begin{equation*}
\partial_{u^n}\mathcal L=\partial_{u^n}J+B_n^Ta^{n+1}.
\end{equation*}
Fix the adjoints at $u$. For every nearby control $\widetilde u$, the state equations give
$j(\widetilde u)=\mathcal L(S_{\rm net}(\widetilde u),\widetilde u,a)$.
The chain rule therefore gives
\begin{align*}
Dj(u)[d]
&=D_u\mathcal L[d]+D_y\mathcal L[z[d]]
=\sum_{n=0}^{N-1}(\partial_{u^n}J+B_n^Ta^{n+1})^Td^n.
\end{align*}
The state term vanishes because $\partial_{y^k}\mathcal L=0$ for $1\le k\le N$ and $z[d]^0=0$. A second application of the chain rule, with the adjoints still fixed, gives
\begin{equation*}
D^2j(u)[d,e]
=D^2_{(y,u)}\mathcal L[(z[d],d),(z[e],e)]+D_y\mathcal L[w[d,e]].
\end{equation*}
The last term is zero, since $w[d,e]^0=0$. Expanding the remaining second derivative of $\mathcal L$ gives the Hessian bilinear form.
\end{proof}

\setlength{\parskip}{2.5pt plus 1pt minus 1pt}
\setlength{\textfloatsep}{8pt plus 2pt minus 2pt}
\setlength{\floatsep}{6pt plus 2pt minus 2pt}
\setlength{\intextsep}{6pt plus 2pt minus 2pt}
\renewcommand{\floatpagefraction}{0.85}
\section{Numerical methods and supplementary checks}\label{app:numerical-details}
\subsection{Algorithm details}\label{app:algorithm-details}
At iteration $k=0,1,\ldots$, write $u_k\in\Uk$ for the control, $g_k\coloneqq g(u_k)$ for its weighted gradient and $d_k\in\R^{N(M+1)}$ for the search direction.\\
\emph{FBS and Anderson.} FBS uses \eqref{eq:fbs}. For damped, regularized Type-II Anderson \cite{Walker2011}, set $f_k\coloneqq \Phi_\theta(u_k)-u_k$, $m_k\coloneqq \min\{k,10\}$, and
\begin{equation*}
\Delta U_k\coloneqq [u_{k-m_k+1}-u_{k-m_k},\ldots,u_k-u_{k-1}],\qquad
\Delta F_k\coloneqq [f_{k-m_k+1}-f_{k-m_k},\ldots,f_k-f_{k-1}].
\end{equation*}
Both matrices lie in $\R^{N(M+1)\times m_k}$, and we solve for $a_k\in\R^{m_k}$:
\begin{equation*}
a_k\coloneqq \argmin_{a\in\R^{m_k}}\{\|f_k-\Delta F_ka\|_2^2+\delta\|a\|_2^2\}.
\end{equation*}
For $m_k>0$, differentiation of the least-squares objective gives the normal equations
\[
(\Delta F_k^T\Delta F_k+\delta I_{m_k})a_k=\Delta F_k^Tf_k.
\]
For $b\ne0$, $b^T(\Delta F_k^T\Delta F_k+\delta I)b
=\|\Delta F_kb\|_2^2+\delta\|b\|_2^2>0$, so $a_k$ is uniquely defined. Omit the correction when $m_k=0$:
\begin{equation*}
u_{k+1}\coloneqq \proj\big(\Phi_\theta(u_k)-\omega(\Delta U_k+\Delta F_k)a_k\big).
\end{equation*}
We use $\omega=0.8$, $\delta=10^{-10}$ and no safeguard.\\
\emph{SPG.} Spectral projected gradient \cite{Birgin2000} uses a scalar step $a_k>0$ and
\begin{equation*}
d_k\coloneqq \proj(u_k-a_kg_k)-u_k,\qquad u_{k+1}\coloneqq u_k+\sigma_kd_k.
\end{equation*}
Because $\Wht$ is positive diagonal, coordinate clipping is also the weighted projection onto the box. Its variational inequality, with $u_k$ as the comparison point, gives
\[
\langle d_k+a_kg_k,-d_k\rangle_{h,\tau}\ge0
\quad\Longrightarrow\quad
\langle g_k,d_k\rangle_{h,\tau}\le-\frac1{a_k}\|d_k\|_{h,\tau}^2.
\]
Thus $d_k\ne0$ is a descent direction. Moreover, $u_k+\sigma d_k$ stays in the box for $0\le\sigma\le1$. The expansion
\[
j(u_k+\sigma d_k)=j(u_k)+\sigma\langle g_k,d_k\rangle_{h,\tau}+o(\sigma)
\]
ensures acceptance for sufficiently small positive $\sigma$. Take the first $\sigma_k$ in $1,1/2,1/4,\ldots$ such that
\[
 j(u_k+\sigma_kd_k)\le
 \max_{\max\{0,k-19\}\le i\le k}j(u_i)
 +10^{-4}\sigma_k\ipt{g_k}{d_k}.
\]
For $s_k\coloneqq u_{k+1}-u_k$, $y_k\coloneqq g_{k+1}-g_k\in\R^{N(M+1)}$, the next step $a_{k+1}$ alternates
\begin{equation*}
\frac{\ipt{s_k}{s_k}}{\ipt{s_k}{y_k}},\qquad
\frac{\ipt{s_k}{y_k}}{\ipt{y_k}{y_k}},
\end{equation*}
clipped to $[10^{-8},10^4]$, with the first ratio used for even $k$. When $\langle s_k,y_k\rangle_{h,\tau}>0$, these ratios follow from the two scalar secant fits
\[
\argmin_{b\in\R}\|bs_k-y_k\|_{h,\tau}^2
=\frac{\langle s_k,y_k\rangle_{h,\tau}}{\langle s_k,s_k\rangle_{h,\tau}},\qquad
\argmin_{a\in\R}\|s_k-ay_k\|_{h,\tau}^2
=\frac{\langle s_k,y_k\rangle_{h,\tau}}{\langle y_k,y_k\rangle_{h,\tau}}.
\]
The first step is the reciprocal of the fitted Hessian scalar $b$, while the second fits the inverse directly. We start with $a_0=0.5$ and reset to $0.5$ whenever $\ipt{s_k}{y_k}\le0$. Rejected trials need only a state solve.\\
\emph{L-BFGS-B} \cite{Byrd1995,Zhu1997} uses $\Wht g$, a memory of 38, \texttt{ftol=0}, at most 60 line-search trials, and a Euclidean projected-gradient tolerance of $10^{-7}\min_{\substack{0\le n<N\\0\le j\le M}}(\tau w_j)$. The iteration limits are 1500 for reference runs and 3000 for refinement and multistart runs.
\subsection{Derivative, Hessian, and balance tests}\label{app:numerical-checks}
\begin{figure}[t]\centering
\includegraphics[width=.485\textwidth]{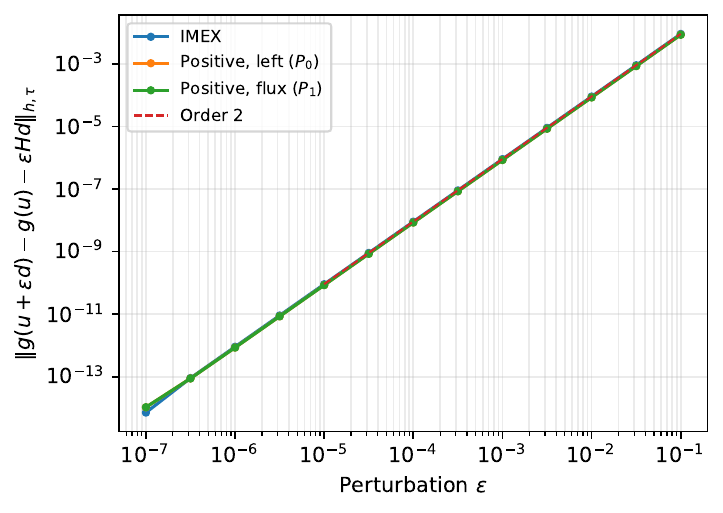}\hfill
\includegraphics[width=.485\textwidth]{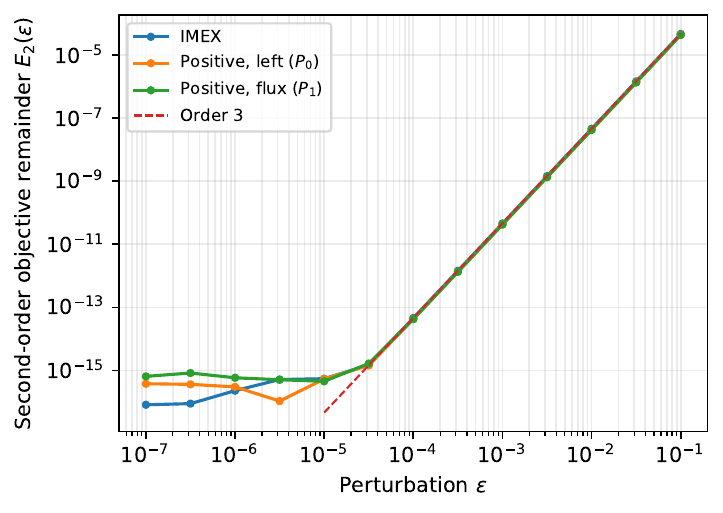}
\caption{Hessian validation: order 2 for $E_g$ (left), order 3 for $E_2$ before cancellation (right). The flux test includes the nonzero terminal adjoint.}\label{fig:taylor}
\end{figure}

Use the reference data of Section~\ref{sec:numerical}.
At $u\equiv0.7$, take $\widetilde d,d,e\in\R^{N(M+1)}$ with
\[
 \widetilde d_j^n\coloneqq0.35\cos(2\pi x_j)\cos(\pi t_n/T)+0.1\xi_j^n,
 \quad d\coloneqq\widetilde d/\|\widetilde d\|_\infty,
\]
where $\xi_j^n$ are independent standard normal samples. For $e$, draw a second independent standard normal array and normalize it in the infinity norm. Set
\begin{equation*}
E_1(\varepsilon)\coloneqq|j(u+\varepsilon d)-j(u)-\varepsilon\ipt{g(u)}d|,\qquad
E_g(\varepsilon)\coloneqq\normt{g(u+\varepsilon d)-g(u)-\varepsilon H(u)d},
\end{equation*}
and
\begin{equation*}
E_2(\varepsilon)\coloneqq\left|j(u+\varepsilon d)-j(u)-\varepsilon\ipt{g(u)}d
-\frac{\varepsilon^2}{2}\ipt{H(u)d}d\right|.
\end{equation*}
The state maps and reduced costs are smooth near $u$. For small $\varepsilon$, the fundamental theorem of calculus gives
\begin{align*}
j(u+\varepsilon d)-j(u)-\varepsilon\langle g(u),d\rangle_{h,\tau}
&=\varepsilon^2\int_0^1(1-t)
\langle H(u+t\varepsilon d)d,d\rangle_{h,\tau}\,\mathrm{d}t,\\
g(u+\varepsilon d)-g(u)-\varepsilon H(u)d
&=\varepsilon\int_0^1[H(u+t\varepsilon d)-H(u)]d\,\mathrm{d}t.
\end{align*}
Subtracting the second-order term in the first identity gives
\[
E_2(\varepsilon)
=\varepsilon^2\left|\int_0^1(1-t)
\langle[H(u+t\varepsilon d)-H(u)]d,d\rangle_{h,\tau}\,\mathrm{d}t\right|.
\]
Let $L_H$ be a Lipschitz constant for $H$ on this neighborhood in the weighted operator norm. Then
\begin{equation*}
E_1(\varepsilon)\le\frac{\varepsilon^2}{2}
\sup_{0\le t\le1}\|H(u+t\varepsilon d)\|_{h,\tau}\|d\|_{h,\tau}^2,\qquad E_g(\varepsilon)\le\frac{L_H\varepsilon^2}{2}\|d\|_{h,\tau}^2,\qquad E_2(\varepsilon)\le\frac{L_H|\varepsilon|^3}{6}\|d\|_{h,\tau}^3.
\end{equation*}
These identities predict orders 2, 2 and 3 when the leading terms do not vanish. The measured orders in Table~\ref{tab:derivatives} come from fitting $\log E$ against $\log\varepsilon$ on $[10^{-3},10^{-1}]$. Smaller steps enter the cancellation regime (Fig.~\ref{fig:taylor}). The \emph{symmetry defect} is
\begin{equation*}
\frac{|\ipt{H(u)d}e-\ipt d{H(u)e}|}{\max\{1,|\ipt{H(u)d}e|,|\ipt d{H(u)e}|\}}.
\end{equation*}
\begin{table}[t]\centering\small
\caption{Derivative checks: $E_g$ tests the full Hessian action, $E_2$ the second variation.}\label{tab:derivatives}
\begin{tabular}{lrrrr}
\toprule
Problem & $E_1$ order & $E_g$ order & $E_2$ order & Symmetry defect \\
\midrule
I & 2.008 & 2.001 & 3.003 & $8.67\times10^{-19}$ \\
$P_0$ & 2.010 & 2.003 & 3.004 & $8.67\times10^{-19}$ \\
$P_1$ & 2.009 & 2.002 & 3.004 & $4.34\times10^{-19}$ \\
\bottomrule
\end{tabular}
\end{table}
With $u\equiv2$, the positive solver keeps the states nonnegative on every tested time grid, down to $N=4$, whereas IMEX gives negative states at $N=4,6,8$. Balance is checked on the raw states, without clipping.
\subsection{Algorithm comparison and rejected runs}\label{app:algorithm-comparison}
\begin{figure}[!htbp]\centering
\includegraphics[width=.325\textwidth]{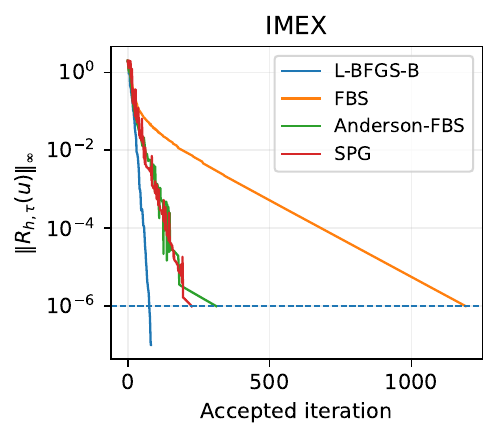}\hfill
\includegraphics[width=.325\textwidth]{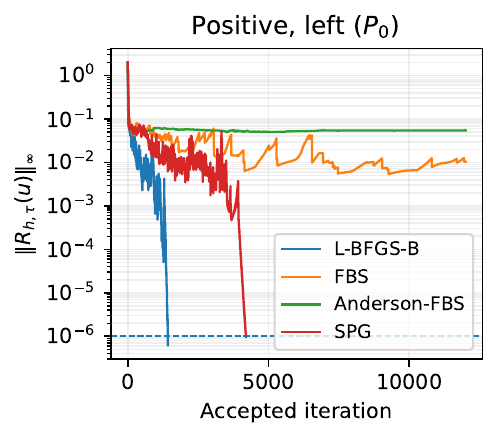}\hfill
\includegraphics[width=.325\textwidth]{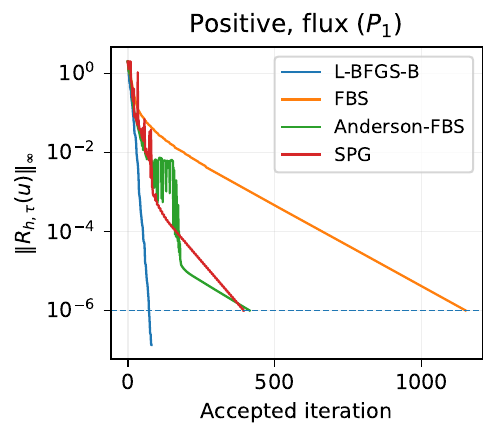}
\caption{Residual histories on separate iteration scales. Dashed lines mark the common tolerance. Failed $P_0$ runs are included.}\label{fig:algorithms}
\end{figure}

All methods start from zero. FBS uses $\theta=0.0093$ and Anderson-FBS $\theta=0.01$, both with a 12000-iteration limit. SPG allows 10000 iterations. We use the same relaxation parameters for all three problems. Anderson changes the iteration map and, for IMEX, converges slightly above the unaccelerated threshold from Theorem~\ref{th:fbs}.
Table~\ref{tab:algorithms} reports final residuals. Figure~\ref{fig:algorithms} shows the residual histories for each method.
\begin{table}[!htbp]\centering\small
\caption{Algorithm comparison. Acceptance requires the independent check $\|R\|_\infty\le10^{-6}$. Counts include terminal verification.}\label{tab:algorithms}
\begin{tabular}{llrrrrrl}
\toprule
Problem & Method & Objective & Iter. & State & Adjoint & $\|R\|_\infty$ & Accepted \\
\midrule
I & L-BFGS-B & -1.5491 & 83 & 103 & 103 & $9.95\times10^{-8}$ & yes \\
I & FBS & -1.5491 & 1191 & 1192 & 1192 & $9.96\times10^{-7}$ & yes \\
I & Anderson-FBS & -1.5491 & 312 & 313 & 313 & $9.94\times10^{-7}$ & yes \\
I & SPG & -1.5491 & 226 & 228 & 227 & $9.86\times10^{-7}$ & yes \\
$P_0$ & L-BFGS-B & -1.5949 & 1427 & 1468 & 1468 & $6.25\times10^{-7}$ & yes \\
$P_0$ & FBS & -1.5944 & 12000 & 12001 & 12001 & $1.03\times10^{-2}$ & no \\
$P_0$ & Anderson-FBS & -1.572 & 12000 & 12001 & 12001 & $5.47\times10^{-2}$ & no \\
$P_0$ & SPG & -1.5947 & 4198 & 4203 & 4199 & $9.70\times10^{-7}$ & yes \\
$P_1$ & L-BFGS-B & -1.521 & 82 & 101 & 101 & $1.32\times10^{-7}$ & yes \\
$P_1$ & FBS & -1.521 & 1150 & 1151 & 1151 & $9.96\times10^{-7}$ & yes \\
$P_1$ & Anderson-FBS & -1.521 & 415 & 416 & 416 & $9.95\times10^{-7}$ & yes \\
$P_1$ & SPG & -1.521 & 394 & 396 & 395 & $9.88\times10^{-7}$ & yes \\
\bottomrule
\end{tabular}
\end{table}
\clearpage
\subsection{State self-convergence}\label{app:state-self-convergence}
With the data of Section~\ref{sec:numerical}, fix the control
\begin{equation*}
u(x,t)\coloneqq0.7+0.25\sin(\pi x)\cos(2\pi t/T).
\end{equation*}
Both tests compare the two terminal states at the same physical time $T$. On a fixed spatial grid with $M=160$, the difference between runs with $N$ and $2N$ time steps is
\[
E_v^{\rm time}(N)
=\left(\sum_{i=1}^2\|v_i^{N,(M,N)}-v_i^{2N,(M,2N)}\|_h^2\right)^{1/2},
\]
where $(M,N)$ records the number of spatial and temporal intervals. For the spatial test, $N=960$ is fixed. We compare the states on grids with $M$ and $2M$ intervals at their common nodes,
\[
E_v^{\rm space}(M)
=\left(\sum_{i=1}^2\sum_{j=0}^M w_j
\left|v_{i,j}^{N,(M,N)}-v_{i,2j}^{N,(2M,N)}\right|^2\right)^{1/2}.
\]
For either difference, write $\Delta=\tau$ or $h$ and define the observed order by
\[
p_{\rm obs}(\Delta)=\log_2\!\left(\frac{E_v(\Delta)}{E_v(\Delta/2)}\right).
\]
If $E_v(\Delta)=c\Delta^p+o(\Delta^p)$ with $c>0$, then $p_{\rm obs}(\Delta)\to p$. For this control, both state schemes give observed orders close to one under time refinement and close to two under space refinement (Table~\ref{tab:fixed}). The two positive-scheme objectives share a state solver, so their state differences are identical. The Gaussian initial profiles in this test have nonzero first derivatives at the endpoints. Proposition~\ref{prop:accuracy} assumes the boundary conditions in \eqref{eq:compatibility} even at $t=0$ and bounds the error over the whole trajectory. Table~\ref{tab:fixed} reports self-convergence at the final time for the data used in the experiments.
\begin{table}[!htbp]\centering\small
\caption{Differences between terminal states on successive grids at fixed control. Time refinement uses $M=160$ and space refinement uses $N=960$.}\label{tab:fixed}
\begin{tabular}{lrrrrr}
\toprule
Study & $N$ or $M$ & IMEX difference & Order & Positive difference & Order \\
\midrule
time & 60 & $5.745\times10^{-2}$ & 1.090 & $1.469\times10^{-1}$ & 0.993 \\
 & 120 & $2.698\times10^{-2}$ & 1.043 & $7.378\times10^{-2}$ & 1.007 \\
 & 240 & $1.310\times10^{-2}$ & 1.021 & $3.671\times10^{-2}$ & 1.006 \\
 & 480 & $6.455\times10^{-3}$ & -- & $1.829\times10^{-2}$ & -- \\
space & 20 & $1.165\times10^{-2}$ & 2.023 & $1.162\times10^{-2}$ & 2.022 \\
 & 40 & $2.866\times10^{-3}$ & 2.006 & $2.860\times10^{-3}$ & 2.006 \\
 & 80 & $7.135\times10^{-4}$ & 2.002 & $7.121\times10^{-4}$ & 2.002 \\
 & 160 & $1.782\times10^{-4}$ & -- & $1.778\times10^{-4}$ & -- \\
\bottomrule
\end{tabular}
\end{table}

\end{document}